\pdfoutput=1 
\documentclass[12pt,a4paper]{article}\usepackage[top=18mm,left=12mm]{geometry}
\usepackage{amsmath,amssymb}
\usepackage{float}\floatstyle{plaintop}\restylefloat{table}
\usepackage[tableposition=top]{caption}
\usepackage{pstricks,fancyvrb,multirow,lineno,url}
\usepackage[pdftex]{graphicx}
\usepackage{listings}\usepackage{paralist, longtable}
\newtheorem{theorem}{Theorem}[section]
\newtheorem{remark}[theorem]{Remark}
\usepackage{etoolbox}\AtBeginEnvironment{remark}{\renewcommand{\baselinestretch}{1}}
\usepackage{titlesec}
\titlespacing*{\section}{0pt}{1.5ex plus 1ex}{1.5ex plus .5ex}
\titlespacing*{\subsection}{0pt}{1.5ex plus 0.5ex}{1.25ex plus .25ex}
\titlespacing*{\subsubsection}{0pt}{1.5ex plus 0.5ex}{0.8ex plus .0ex}
\newtheorem{lemma}[theorem]{Lemma}
\newtheorem{corollary}[theorem]{Corollary}
\newtheorem{solution}[theorem]{Solution}
\newtheorem{exercise}[theorem]{Exercise}
\newtheorem{definition}[theorem]{Definition}
\def\bexe{\begin{exercise}}\def\eexe{\eex\end{exercise}}
\def\bsol{\begin{solution}}\def\esol{\eex\end{solution}}
\def\bexa{\begin{example}}\def\eexa{\end{example}}
\def\brem{\begin{remark}}\def\erem{\end{remark}}
\def\bthm{\begin{theorem}}\def\ethm{\end{theorem}}
\def\blem{\begin{lemma}}\def\elem{\end{lemma}}

\def\bcor{\begin{corollary}}\def\ecor{\end{corollary}}
\def\bdefi{\begin{definition}}\def\edefi{\end{definition}}
\newcommand{\IDEA}{\textbf{Idea of the Proof.} }

\newcommand{\abs}[1]{\lvert#1\rvert}\newcommand{\norm}[1]{\lVert#1\rVert}

\def\bmip{\begin{minipage}{\textwidth}}\def\emip{\end{minipage}}
\def\huga#1{\begin{gather} #1 \end{gather}}
\def\hugast#1{\begin{gather*} #1 \end{gather*}}
\def\hual#1{\begin{align} #1 \end{align}}
\def\hualst#1{\begin{align*} #1 \end{align*}}
\def\hueq#1{\begin{equation} #1 \end{equation}}

\newcommand{\btab}[2]{\begin{tabular}{#1}#2\end{tabular}}

\newcommand{\R}{{\mathbb R}}
\newcommand{\C}{{\mathbb C}}\newcommand{\N}{{\mathbb N}}
\newcommand{\Z}{{\mathbb Z}}

\def\CO{{\cal O}}\def\CS{{\cal S}}
\def\CM{{\cal M}}\def\CN{{\cal N}}

\def\diag{{\rm diag}}

\def\vx{\vec x}

\def\ga{\gamma}\def\om{\omega}
\def\ds{\displaystyle}
\def\vt{\vartheta}\def\pa{{\partial}}\def\lam{\lambda}
\newcommand{\bi}{\begin{itemize}}\newcommand{\ei}{\end{itemize}}
\newcommand{\ben}{\begin{enumerate}}\newcommand{\een}{\end{enumerate}}
\newcommand{\bce}{\begin{center}}\newcommand{\ece}{\end{center}}
\newcommand{\bci}{\begin{compactitem}}\newcommand{\eci}{\end{compactitem}}
\newcommand{\bcen}{\begin{compactenum}}\newcommand{\ecen}{\end{compactenum}}
\newcommand{\bcena}{\begin{compactenum}[(a)]}
\newcommand{\reff}[1]{(\ref{#1})}
\newcommand{\ov}[1]{{\overline {#1}}}

\newcommand{\spr}[1]{\left\langle #1 \right\rangle}
\newcommand{\hs}[1]{{\hspace{#1}}}\newcommand{\vs}[1]{{\vspace{#1}}}

\def\eps{\varepsilon}
\def\ra{\rightarrow}

\newcommand{\barr}{\begin{array}}\newcommand{\earr}{\end{array}}
\newcommand{\bpm}{\begin{pmatrix}}\newcommand{\epm}{\end{pmatrix}}
\newcommand{\bsm}{\left(\begin{smallmatrix}}
\newcommand{\esm}{\end{smallmatrix}\right)}
\newcommand{\ba}{\begin{array}}\newcommand{\ea}{\end{array}}
\def\dd{\, {\rm d}}\def\ri{{\rm i}}

\def\er{{\rm e}}

\def\om{\omega}\def\Om{\Omega}

\def\bphi{\varphi}\def\del{\delta}

\def\nab{\nabla}\def\eex{\hfill\mbox{$\rfloor$}}

\def\sign{{\rm sign}}
\def\Del{\Delta}

\def\sig{\sigma}
\def\al{\alpha}

\def\Ga{\Gamma}\def\Lam{\Lambda}
\def\kap{\kappa}\def\Sig{\Sigma}

\def\bd{\begin{displaymath}} \def\ed{\end{displaymath}}
\def\ba{\begin{array}} \def\ea{\end{array}}
  
\def\eps{\varepsilon}

\def\ind{{\rm ind}}

\def\pdep{{\tt pde2path}}
\def\auto{{\sc Auto}}\def\coco{{\sc Coco}}

\def\auto{{\tt AUTO}}

\def\mlab{{\sc Matlab}}

\def\dhome{/hh/path/pde2path/demos/acsuite}

\newlength{\tew}

\def\ig{\includegraphics}

\def\wt{\tilde{w}}

\renewcommand{\arraystretch}{1.05}\renewcommand{\baselinestretch}{1.0}
\def\medskip{}\def\bigskip{}
\usepackage{setspace}
\def\sm{\small}\def\rb{\raisebox}
\def\bsub{\begin{subequations}}
\def\esub{\end{subequations}}
\def\Mf{{M^{{\text{full}}}}}\def\MV{{M}}%_{{\rm V}}}}
\def\Xhem{X_{\text{2HS}}}

\usepackage{hyperref}
\def\taskip{\renewcommand{\arraystretch}{1}\renewcommand{\baselinestretch}{1}}
\def\teskip{\renewcommand{\arraystretch}{1.1}\renewcommand{\baselinestretch}{1.1}}
\def\hulst#1#2{\taskip\lstinputlisting[#1]{#2}\teskip}
\def\hutab#1{\taskip\begin{\table}#1\end{table}\teskip}
 
\newfloat{Algorithm}{ht}{lop}[section]
\def\slfig{\opt{hu,ho}{\begin{figure}[ht]}\opt{sl}{\begin{figure}[H]}}
\def\slfigH{\opt{hu,ho}{\begin{figure}[H]}\opt{sl}{\begin{figure}[H]}}
\def\sltab{\opt{hu,ho}{\begin{table}[ht]}\opt{sl}{\begin{table}[H]}}
\def\gptool{{\tt gptoolbox}}
\def\wt{\widetilde}\def\srot{s_{\text{rot}}}\def\delm{\delta_{\text{mesh}}}

\definecolor{codegreen}{rgb}{0,0.6,0}
\definecolor{codegray}{rgb}{0.5,0.5,0.5}
\definecolor{codepurple}{rgb}{0.58,0,0.82}
\definecolor{backcolour}{rgb}{0.95,0.95,0.92}
 
\lstdefinestyle{mystyle}{
    backgroundcolor=\color{backcolour},   
    commentstyle=\color{codegreen},
    keywordstyle=\color{black},
    numberstyle=\small\color{codegray},
    stringstyle=\color{codepurple},
    basicstyle=\footnotesize\ttfamily,
    breakatwhitespace=false,         
    breaklines=true,                 
    captionpos=b,                    
    keepspaces=true,                 
    numbers=left,                    
    numbersep=5pt,                  
    showspaces=false,                
    showstringspaces=false,
    showtabs=false,                  
    tabsize=2, 
stepnumber=5,
language=matlab,  
  xleftmargin=1mm,
}
 
\begin{document}
\text{}\vspace{0mm} 
\bce\Large
Differential geometric bifurcation problems in \pdep\ -- algorithms and 
tutorial examples\\[4mm]
\normalsize 
Alexander Meiners, Hannes Uecker\\[2mm]
\footnotesize
Institut f\"ur Mathematik, Universit\"at Oldenburg, D26111 Oldenburg, 
alexander.meiners@uni-oldenburg.de, hannes.uecker@uni-oldenburg.de\\[3mm]
\normalsize
\today
\ece 
\begin{abstract} We describe how some differential geometric 
bifurcation problems can be treated with the \mlab\ continuation and 
bifurcation toolbox \pdep. The basic setup 
consists in solving the PDEs  for the normal 
displacement of an immersed surface $X\subset\R^3$ and subsequent update of $X$ 
in each continuation step, combined with bifurcation detection and localization, 
followed by possible branch switching.  
Examples treated include some minimal surfaces such as 
Enneper's surface and a Schwarz--P--family, some non--zero 
constant mean curvature surfaces such as 
liquid bridges and nodoids, and some 4th order biomembrane models. 
In all of these we find interesting symmetry breaking bifurcations. 
Some of these are (semi)analytically known and thus are used as benchmarks.
\end{abstract}
\tableofcontents

\section{Introduction} \label{isec}
There are various algorithms and toolboxes for 
the numerical continuation %(in parameters) 
and bifurcation analysis for solutions of partial differential equations 
(PDEs), for instance \auto\ \cite{auto}, %(wrt to PDEs focussing on 1D boundary value problems), 
\coco\ \cite{cocobook}, %(Continuation Core, which in principle allows great flexibility by delegating the PDE definition/discretization to the user), 
{\tt BifurcationKit.jl}\ \cite{veltz20}, and \pdep\ 
\cite{p2pbook, p2phome}. In its standard setup, 
\pdep\ is for PDEs for 
functions $u:\Om\times\Lambda\to \R^N$, 
where $\Om\subset\R^d$ is a fixed domain, $d=1,2$, or 3, $N\in\N$, and 
$\Lam\subset \R^p$ is a parameter domain, or for time--dependent 
functions $u:I\times\Om\times\Lambda\to\R^N$, $I\subset\R$, which then includes the 
continuation and bifurcation of time periodic orbits. Essentially, %As far as we know, 
this also applies to {\tt BifurcationKit.jl}, while wrt PDEs \auto\ 
originally focusses on 1D boundary value problems, and \coco\ in principle allows great 
flexibility by delegating the PDE definition/discretization to the user. 

However, none of these packages seem directly applicable to differential 
geometric PDEs in parametric form, which deal directly 
with  manifolds, e.g., surfaces in 2D, which are not graphs over 
a fixed domain. There are well established numerical methods 
for the discretization of such PDEs, for instance the surface finite 
element method (surface FEM) \cite{DzE13}, but 
there seem to be few algorithms or packages 
which combine these with continuation and bifurcation. 
Two notable exceptions are 
the algorithm from \cite{bru18}, and the {\tt SurfaceEvolver}, for 
which bifurcation aspects are for instance discussed in \cite{Bra96}.

\begin{figure}[t]
\centering 
\btab{l}{
\btab{ll}{{\sm (a)}&{\sm (b)}\\
\hs{0mm}\ig[width=0.23\tew,height=40mm]{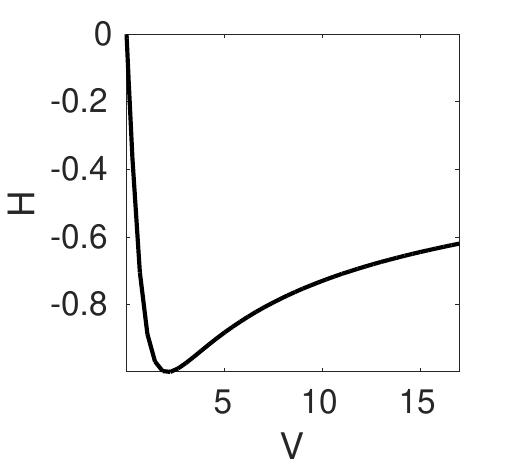}
\hs{0mm}\ig[width=0.23\tew]{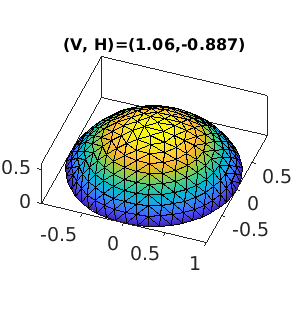}\ig[width=0.23\tew]{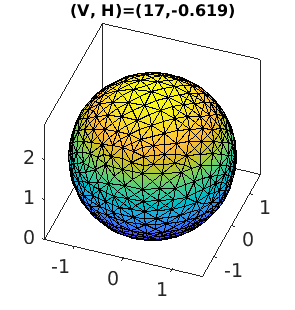}
&\ig[height=0.3\tew,height=44mm]{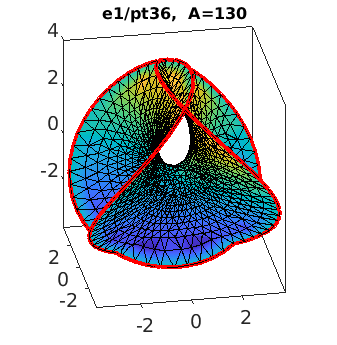}}\\
\btab{lll}{{\sm (c)}&{\sm (d)}&{\sm (e)}\\
\hs{0mm}\rb{6mm}{\ig[width=0.35\tew]{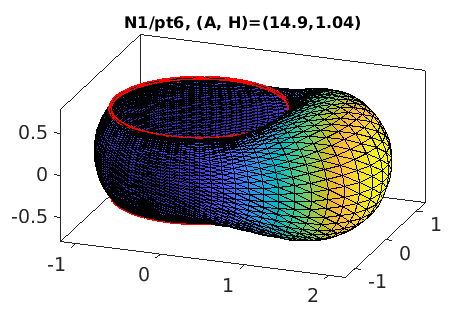}}
&\ig[width=0.34\tew]{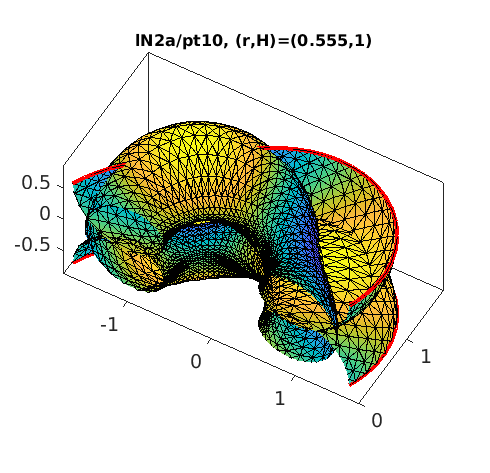}
&\rb{10mm}{\ig[width=0.24\tew]{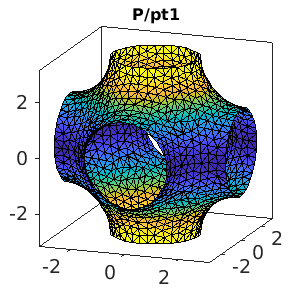}}
}
}
\vs{-8mm}
\caption{{\small Preview of solutions (solution branches) we compute. 
(a) $H$ over $V$ for spherical caps, 
and sample solutions, \S\ref{scsec}. The colors indicate $u$ in the 
last step, yellow$>$blue, and thus besides giving visual structure 
to $X$ indicate the ``direction'' of the continuation. 
$H$ is negative since here $N$ is the outer 
normal. (b) Enneper minimal surface (a bounded part, with the boundary 
shown in red), \S\ref{ensec}. 
(c) A liquid bridge between two circles, with excess volume and 
hence after a symmetry breaking bifurcation, \S\ref{snsec}. (d) A nodoid with 
periodic BCs, cut open for plotting, \S\ref{nodnumpBC}. 
(e) A Schwarz P surface, \S\ref{SPsec}. Samples (b)--(e) are each 
again from branches of solutions of the respective 
problems, see Figures \ref{ef1}, \ref{nodf0},\ref{nodf4}, and \ref{tpsf1} 
for the bifurcation diagrams.
\label{f1}}}
\end{figure}

Here we present an extension of \pdep\ aimed at geometric PDE 
bifurcation problems. We focus on constant 
mean curvature (CMC) surfaces, which are not necessarily graphs, 
and with, e.g., the mean curvature, or the area or enclosed 
volume as the primary bifurcation parameter. 
See Fig.\,\ref{f1} for a preview of the type of solutions we compute. 
For $X$ a two dimensional surface %(possibly with boundary) 
immersed in $\R^3$, we for instance want to study the parameter 
dependent problem 
\begin{subequations}\label{CMC}
\hual{
H(X)-H_0&=0,\\ 
V(X)-V_0&=0, 
}
\end{subequations}
possibly with boundary conditions (BCs) in (a), 
where $H(X)$ is the mean curvature at each point of $X$, 
and $V(X)$ is the (properly defined) 
volume enclosed by $X$. The system \reff{CMC} is obtained for 
minimizing the area $A(X)$ under the volume constraint $V(X)=V_0$, i.e., 
as the Euler--Lagrange equations for minimizing the energy 
\huga{\label{cmcen}
E(X)=A(X)+H_0(V(X)-V_0), 
}
and $V_0\in\R$ typically plays the role 
of an ``external continuation parameter'', while $H_0$, which for instance describes a 
spatially constant pressure difference for interfaces between fluids, is ``free''.  
%Naturally, these roles of $H_0$ and $V_0$ can 
%also be interchanged (continuation of $V$ as a function of $H$). 

Following \cite{bru18}, our setting for \reff{CMC} and 
generalizations is as follows. 
Let $X_0$ be a surface satisfying \eqref{CMC} for some $V_0$ and $H_0$, 
and define a new surface via $X=X_0+u\,N_0$, $u:X_0\to\R$ with suitable 
boundary conditions, where 
$N_0:X_0\to\mathbb S^2$ is (a choice of) the unit normal vector field of $X_0$. 
Then (\ref{CMC}) reads 
\def\Hti{\wt H}\def\Vti{\wt V}
\bsub\label{feq}
\huga{
	G(u,\Hti)=H(X)-\Hti=0 \text{ (with boundary conditions if applicable)}, 
}
which is a quasilinear elliptic equation for $u$ coupled to the 
volume constraint 
\huga{q(u)=V(X)-\Vti. } 
\esub 
Thus, 
\huga{\label{n1}
\text{after solving \reff{feq} for $u,\Hti,\Vti$ 
we can update $X_0=X_0+uN_0$, $H_0=\Hti, V_0=\Vti$, and repeat.}
} 
 
We generally compute (approximate), e.g., the mean curvature $H$ 
from a surface FEM discretization of $X$, see \S\ref{ddgsec}. 
Alternatively, we may assume a parametrization $\phi:\Omega\ra\R^3$ 
of $X$ over some bounded domain $\Om\subset\R^2$, and compute, e.g., 
$H$ from a classical 2D FEM mesh in $\Om$, although this 
is generally more complicated and less robust than using surface meshes, 
and we mainly review it  for completeness in App.~\ref{sc2sec}. 
Both approaches can and usually must be combined 
with adaptive mesh refinement and coarsening as $X$ changes. 
Both can also be applied to other 
geometric PDEs, also of higher order, for instance  
fourth order biomembrane model, see \S\ref{wsec}.  
In this case, the analog of (\ref{feq}a) 
can be rewritten as a system of (2nd order) 
PDEs for a vector valued $u$, and the same ideas apply.

\begin{figure}[ht]
\btab{p{0.78\tew}l}{{\small 
\bci 
\item Make a directory, e.g., {\tt myp2p} anywhere in your path. 
Download \pdep\ from \cite{p2phome} to {\tt myp2p} and unpack, 
which gives you the \pdep\ home directory {\tt myp2p/pde2path}. 
\item In \mlab, change into \pdep/ and
call {\tt setpdepath} to make the libraries available. 
(We also recommend {\tt ilupack} \cite{ilupack}, which is not 
used here but otherwise in \pdep\ for large scale problems). 
\item Test, i.e: change directory into {\tt pde2path/demos/geomtut/spcap1}, 
load the script {\tt cmds1.m} into the 
editor (i.e., type {\tt edit cmds1.m} at the 
command line). To get an understanding what command does what, 
we then recommend to run {\tt cmds1.m} in ``cell--mode'', i.e., to 
proceed ``cell--by--cell''.  
\item Find a demo that is closest to the problem 
you want to study; copy this demo directory to a new directory 
{\tt myproblem/} or any other name 
(we recommend not as a subdirectory of {\tt pde2path} 
but {\em somewhere else}, for instance in a subdirectory {\tt myproblems} 
of {\tt myp2p}. In {\tt myproblem}, modify the relevant files 
(usually at least {\tt *init} and {\tt cmds}) and explore. 
\eci }
&
\hs{2mm}\rb{-58mm}{\ig[width=0.24\tew]{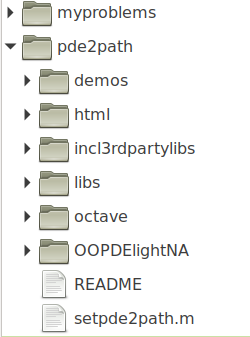}}
}
\vs{-4mm}
\caption{\small Installation of \pdep\ (version 3.1), 
and a ``typical'' directory structure of {\tt myp2p}.
\label{dtf}}
\end{figure}

Our work comes with a number of demos which are subdirectories
of {\tt pde2path/demos/geomtut}, see Table \ref{gtab1}. 
The rather large number of demos is aimed at showing versatility, 
and, more importantly, is due to our own needs for extensive testing, 
in particular  
of various BCs, and of different mesh handling strategies. 
Table \ref{atab} summarizes acronyms and notation 
used throughout, and Fig.\ref{dtf} explains the basic installation steps  
for \pdep. See also \cite{qsrcb} for a quick overview of installation 
and usage of \pdep, and of all demos coming with \pdep, and 
\cite{actutb} or \cite[Chapters 5 and 6]{p2pbook} for getting started 
with \pdep\ via simple classical PDEs. 

\taskip
\begin{table}[ht]\caption{{\sm Demo directories 
in  {\tt pde2path/demos/geomtut}. The first two and last three are 
rather introductory and not dealing with bifurcations. \label{gtab1}}}
\centering %\bce
\vs{-2mm}
{\small 
\begin{tabular}{ll}%p{0.14\tew}|p{0.82\tew}} 
directory&remarks\\
\hline
{\tt spcap1}&Spherical caps via surface meshes, introductory demo. \\
{\tt bdcurve}&Experiments on minimal surfaces with different boundary curves.\\
{\tt enneper}&Bifurcation from Enneper's surface, closely related 
to {\tt bdcurve}. \\
%{\tt libri}&Liquid bridges\\
{\tt nodDBC}&Nodoids with Dirichlet BCs, including so called liquid bridges.\\
{\tt nodpBC}&Nodoids with periodic BCs.\\
{\tt TPS}&Triply Periodic Surfaces, here Schwarz P.\\
{\tt biocyl}&Helfrich cylinders with clamped BCs 
as an example of a 4th order problem.\\
{\tt biocaps}&Disk type solutions as a variant of {\tt biocyl}.\\
\hline
{\tt spheres}&Continuation of spheres, and tests for VPMCF, \S\ref{spsec}.\\
{\tt hemispheres}&Continuation of hemispheres on a supporting plane, 
and VPMCF, \S\ref{hssec}.\\
{\tt spcap2}&Spherical caps via 2D FEM in the preimage, \S\ref{sc2sec}.\\
\end{tabular}
}
%\ece\vs{-4mm}
\end{table}
\teskip

\taskip
\begin{table}[ht]\caption{Notations and acronyms; for given 
$X_0$, quantities of 
$X=X_0+uN_0$ will also be considered as functions of $u$, e.g., 
$A(u)=A(X_0+uN_0)$. }\label{atab}
\centering\vs{-2mm}
{\small 
\begin{tabular}{p{0.16\tew}|p{0.32\tew}|p{0.16\tew}|p{0.3\tew}} 
%\hline
$X$&surface immersed in $\R^3$&$N=N(X)$&surface unit normal\\
$A{=}A(X){=}A(u)$&{area of $X$, resp.~of $X{=}X_0{+}uN_0$}&
$V=V(X)$&(algebraic) volume, e.g., \reff{vol}\\
$H=H(X)$&mean curvature, e.g., \reff{diff-mean}&$K=K(X)$&Gaussian curvature \\
$G(u,\lam)=0$&{generic form of a PDE such as (\ref{feq}a), 
$\lam$ as a generic parameter}&
ind$(X)$&index, i.e., number of unstable eigenvalues of linearization\\
$L=\pa_u H(u)$&Jacobi op.~(with BCs)&$q(u,\lam)=0$&{generic constraint such as (\ref{feq}b)}\\
%$Z_2$ (symmetry)&discrete symmetry group, e.g., $z\mapsto -z$ for the $z$--coordinate&$S^1$ (symmetry)&continuous symmetry group of rotations.\\
\hline 
BC&boundary condition&DBC/NBC&Dirichlet/Neumann BC\\
pBC&periodic BC&PC&phase condition\\
BP/FP&branch/fold point&CMC&constant mean curvature\\%FP&fold point\\
TPS&triply periodic surface&
TPMS&triply periodic minimal surface\\
MCF&mean curvature flow&VPMCF&volume preserving MCF
\end{tabular}
}
%\ece\vs{-4mm}
\end{table}
\teskip

\brem\label{dynrem}{\rm Here we focus on stationary problems of type 
\reff{CMC}, which give critical points of the 
volume preserving mean curvature flow (VPMCF). 
A time $t$ dependent 2D manifold 
$X(t)\subset\R^3$ deforms by mean curvature flow (MCF) if, 
assuming the correct sign for $H$, i.e., $H>0$ for $X$ bounding a 
convex body and $N$ the inner normal, 
\huga{\label{mcf}
\dot X=-H(X)N. 
}
This is the $L^2$ 
gradient flow for the area functional $A(X)$, and can be considered 
as a quasilinear parabolic PDE, at least on short times. 
For closed and compact $X$ there always is 
finite time blowup (generically by shrinking to a ``spherical point''), and 
we refer to \cite{mante11} for an introduction to this huge field, 
which inter alia heavily relies on maximum (comparison) principles. 
%, i.e., ``avoidance principles''. 

The VPMCF reads  
\huga{\label{vpmcf}
\dot X=-(H(X)-\ov H)N, \quad 
\ov H=\frac 1{A(X)} \int_X H\dd S, 
} 
and for closed $X$ conserves the enclosed volume $V(X)$. 
For non--closed $X$ one typically 
studies Neumann type BCs on ``support planes'', see, e.g., 
\cite{hart13} %, AGM15},  
and in most cases the analysis is done near axisymmetric states 
such as spheres, spherical caps, and cylinders. 
In general, the existence and regularity theory for \reff{vpmcf} is 
less well understood than for \reff{mcf} 
due to the lack of general maximum principles for \reff{vpmcf}. 

Our notion of stability of solutions of \reff{CMC} 
(indicated by thick lines in bifurcation diagrams, while 
branches of unstable solutions are drawn as thinner lines) 
refers to \reff{vpmcf} if we have an active volume constraint such as 
(\ref{feq}b), and to \reff{mcf} if not, with the exception of the fourth 
order problems in \S\ref{wsec}, see Remark \ref{hstabrem}. 
Moreover, by ind$(X)$ we denote 
the number of unstable eigenvalues of the linearization of (the 
discretization of) \reff{feq}, including the constraints (if active). 

We also provide very basic setups to numerically integrate \reff{mcf} and 
\reff{vpmcf}  
by explicit Euler stepping. This often has to be combined 
with mesh adaptation, and in this case $A$ does not necessarily decrease 
monotonously for MCF. Moreover, our VPMCF typically conserves $V$ only 
up to $0.5\%$ error. Thus, both are not necessarily efficient or highly 
accurate, 
but can be used to generate initial guesses for the continuation of 
steady states of \reff{CMC}. See \S\ref{scsec}, \S\ref{ensec} (MCF) 
and \S\ref{spsec}, \S\ref{hssec} (VPMCF) 
for examples, and, e.g., \cite{BNP10,BGN20,BGN22} for much more sophisticated 
numerical algorithms for geometric flows including 
\reff{mcf} and \reff{vpmcf}, and 
detailed discussion. 
}
\eex
\erem 

The remainder of this note is organized as follows: 
In \S\ref{dasec} we review some differential geometric background, and 
the \pdep\ data structures and functions to deal with 
geometric PDEs. The central \S\ref{exsec} discusses the main demos, 
and in %the short 
\S\ref{sosec} we give a summary, and an outlook on ongoing 
and future work. In \S\ref{hssec0} we comment on 
the further demos {\tt spheres} and {\tt hemispheres}, 
which do not show bifurcations but deal with VPMCF and Neumann type free BCs, 
and present a classical FEM setup for spherical caps. 
See also \cite{sig} for supplementary information (movies) 
on some of the rather complicated bifurcation diagrams we obtain.

\section{Geometric background, and data structures}
\label{dasec} 
\subsection{Differential geometry}\label{geom}
We briefly review the geometric PDE setup, and recommend \cite{Des04,Tapp16,
UY17} for further background, among many others. 

Throughout, let $\Sig$ be a 2D connected compact orientable manifold, with coordinates $x,y$, and 
possibly with boundary 
$\pa \Sig$, and for some $\al\in(0,1)$ immersed by  $X\in C^{2,\al}(\Sig,\R^3)$. 
By pulling back the standard metric of $\R^3$ we obtain the first and second 
fundamental forms on $\Sig$ expressed via $X$ as 
	\hual{ \label{fundfform}
g&=\bpm g_{11}& g_{12}\\g_{12}& g_{22} \epm=
\bpm \norm{X_x}^2&\spr{X_x,X_y}\\ 
\spr{X_x,X_y}&\norm{X_y}^2\epm, \qquad
h= \bpm h_{11} &h_{12}\\ h_{21}& h_{22} \epm =\bpm
\spr{X_{xx},N} 
& \spr{X_{xy},N}\\\spr{X_{xy},N}& \spr{X_{yy},N}\epm,
	}
with unit normal $N$, which we consider as a field on $\Sig$, or locally on $X$, which 
will be clear from the context. The mean curvature $H$ then is 
\hueq{\label{h1} 
H=\frac 1 2 \frac{h_{11}g_{22}-2h_{12}g_{12}+h_{22}g_{11}}{g_{11}g_{22}-g_{12}^2}, 
}
which is the mean of the minimal and maximal normal curvatures 
$\kap_1$ and $\kap_2$. The Gaussian curvature is 
\hueq{\label{k1} 
K=\kap_1\kap_2. 
}
The sign of $H$ depends on the orientation of $X$, i.e., on the choice of $N$. 
A sphere has positive $H$ iff 
$N$ is the inner normal. The Gaussian curvature does not depend on $N$ or 
any isometry of $\Sig$ (Gau\ss' Theorema egregium). 
%In abuse of notation 
%we also refer to geometric quantities induced by the pullback metric 
%$g$ as quantities on $X$.

A generalization of the directional derivative of a 
function $f$ to vector fields or tensors is the covariant 
derivative $\nabla_Z$ for some vector field $Z$ on $X$. 
For a vector field $Y$, the 
covariant derivative in the $j$'th coordinate direction is
 defined as $\nab_j Y_i:=\frac{\pa Y_i}{\pa x_j}+\Ga_{jk}^iY_k$, 
and for a 1-form $\om$ we have 
$\nab_j\om_i:=\frac{\pa \om_i}{\pa x_j}-\Ga_{jk}^i \om_k$,  
with the Christoffel symbols 
$\Ga_{jk}^i=\frac 1 2
g^{il}(\pa_{x_j}g_{kl}+\pa_{x_k}g_{jl}-\pa_{x_l}g_{jk})$, 
where $g^{ij}$ are the entries of $g^{-1}$, and where we use Einstein's 
summation convention, i.e., summation over repeated indices. 
The covariant derivative is linear in the first argument, giving
a general definition of $\nabla_Z Y$ with some vector field $Z$, and 
if $f$ is a function on $X$, then 
\huga{\label{coZ} 
\nabla_Z f =\spr{g\nabla f,Z}_{\R^2}. %g(\nabla f,Z). 
} 
Throughout we are dealing with surfaces ($2$ dimensional manifold
immersed into $\R^3$), hence the gradient $\nabla$ 
is the {\em surface gradient}, i.e., the usual 
gradient $\nabla_{\R^d}$ in $\R^3$ projected onto the tangent space, 
\huga{\label{sgrad}
\nabla f =\nabla_{\R^3} f-\spr{\nabla_{\R^3}f,N}N, 
}
which later will be needed to (formulate and) implement phase
conditions, and, e.g., Neumann type BCs. 
%Later we will impose some BCs including conditions on the surface%gradient. 
This also
gives the {\em Laplace Beltrami operator} via 
\hugast{
\Delta f=g^{ij}\nabla_i\nabla_j f, 
}
which then also applies to general tensors. 
Using the Gau\ss--Weingarten relation 
$\frac{\pa^2 X}{\pa x_i\pa x_j}=\Ga_{ij}^k\frac{\pa X}{\pa x_k}+h_{ij} N$
we obtain 
\hugast{
\Delta X=g^{ij}\nab_i\nab_j X = g^{ij}\left(\frac{\pa^2 X}{\pa
x_i\pa x_j}-\Ga^i_{jk}\frac{\pa X}{\pa x_k}\right)=g^{ij}h_{ij}
N=2H(X)N=2\vec{H}(X),
} where $\vec{H}(X)$ is called the mean curvature
vector,  and 
\huga{\label{diff-mean} H(X)=\frac 1 2\spr{\Del X,N}. } 

The area of $X$ is 
\hual{
A(X)&=\int_X\dd S, \label{area}
} 
and, based on Gau\ss' divergence theorem, the (algebraic) volume is 
\hueq{
V(X)=\frac{1}{3} \int_X \spr{X,N}\dd S\label{vol}.
}  
If $X$ is a closed manifold bounding $\Om\subset\R^3$, i.e.,
$\pa\Om=X$, and $N$ the outer normal, then $V(X)=|\Om|$ is the physical
volume. If $X$ is not closed, then we typically need 
to add a third of the flux of $\vec{x}$ through the open ends 
to $V(X)$ (see the examples below).  

We denote the set of all immersed surfaces with the same boundary $\ga$ by
\huga{
\text{$\mathcal N_\ga=\{X:\text{$X$ is an immersed surface as above and }
\pa X=\ga\}$.}
} 
The following Lemma states that all immersions $Y\in \CN_\ga$ 
close to $X$ are graphs over $X$ determined by a 
function $u$ as $Y=X+uN$, which justifies our numerical approach \reff{n1}. 
The condition that $Y$ has the same boundary as 
$X$ in general cannot be dropped, as obviously motions of $\pa X$
tangential to $X$ cannot be captured in the form $X+uN$; see
\S\ref{bc1sec} for an illustration.  
\blem \label{l-normvar} \cite{KPP17}.  For 
$X\in C^{2,\alpha}(\Sig,\R^3)$ with boundary $\pa X=\ga$ there exists a
neighborhood $U\subset C^{2,\alpha}(\Sig,\R^3)$ of $X$ such that for all
$Y\in U\cap \mathcal N_\ga$ there exists a diffeomorphism
$\phi:\Sig\rightarrow\Sig$ and a $u\in C^{2,\alpha}(\Sig)$ such that
\hueq{Y\circ\phi=X +u\,N.}
\elem 

Assume that a CMC surface $X_0$ with boundary $\pa X_0{=}\ga$ 
and volume $V(X_0){=}V_0$ belongs to a family of CMC surfaces $X_t$, 
$t\in(-\eps,\eps)$ for some $\eps>0$. 
For example, the spherical caps $S_t$ from Fig.\,\ref{f1}(a) 
with the boundary $\ga=\{(x,y,0)\in\R^3: x^2+y^2=1\}$ are a 
family of CMC immersions fully described by the height 
$t\in\R$. By Lemma
\ref{l-normvar}, the parameter $t$ uniquely defines $u$ in a small 
neighborhood, i.e., $X_t=X+u\,N$, and the system of equations for $u$
reads \hual{\label{cmc-con} 
H(u)-H_0=0, 
} 
for some $H_0\in\R$, where we abbreviate $H(u)=H(X+u\,N)$, etc. 
%Following \cite[\S2.1]{Lop13} 
If we consider variational vector fields at $X_0$ in the form 
$\ds\psi=\frac{\pa X_t}{\pa t}\vert_{t=0}=uN$, and additionally 
assume that $X_t\in \CN_\ga$, then necessarily 
\huga{\label{DBC1}
u|_{\pa X}=0, \text{ Dirichlet boundary conditions (DBCs).}
}  
Such $X_t$ are called {\em admissible variations} in \cite[\S2.1]{Lop13}, 
and we have the following results on derivatives of $A$ and $V$. 
\blem\label{derlem}\cite[\S2.1]{Lop13} 
For an admissible one parameter variation $X_t$ of
$X\in C^{2,\alpha}(\Sig)$ and variational vector 
fields $\psi=\frac{\pa X}{\pa t}\big |_{t=0}=uN$ the 
functions $t\mapsto A(t)=A(X_t)$ and $t\mapsto V(t)=V(X_t)$ 
are smooth, and 
\hual{ \label{f-area}
&V'(0)=\int_{X_0} u\dd S, \quad 
A'(0)=-2\int_{X_0} H_0 u \dd S,\text{ and }
A''(0)= -\int_{X_0} (\Delta u
+\|S_0\|^2u)u \dd S,  
}
where $\norm{S_0}^2=4H_0^2-2K_0$ with the Gaussian
curvature $K_0$. Thus 
\hual{\label{f-mean} 
\frac{\dd}{\dd t}H(X_t)\Big\vert_{t=0}=-\Delta u-\norm{S_0}^2 u, 
} 
and the directional derivative \eqref{f-mean} is given by 
the self-adjoint Fredholm operator $L$ on $L^2(X_0)$ with 
\hueq{\label{jac}L=\pa_u H(0)=-\Del
 -\norm{S_0}^2,\quad\text{with DBCs.}
}
\elem

\brem\label{georem}{\rm 
a) The operator in \eqref{jac} without BCs is called {\em Jacobi
operator}, and a 
nontrivial kernel function is called a {\em Jacobi field} on $X=X_0$. 
An immersion $X$ with a Jacobi field satisfying the BCs is called 
{\em degenerate}. 
The Fredholm property allows the use of the 
Crandall-Rabinowitz bifurcation result \cite{CR71}: Given a $C^1$
branch $(t_0-\eps,t_0+\eps)\ni t\mapsto X_t$, if $X_t$ is
non--degenerate for $t\in(t_0-\eps,t_0)\cup(t_0,t_0+\eps)$, and if at
$t_0$ a {\em simple} eigenvalue $t\mapsto \mu_0(t)$ crosses
transversally, i.e., $\mu(t_0)=0$, $\mu_0'(t_0)\ne 0$, then a branch
$\wt X_t$ bifurcates at $t_0$.  

See also \cite{KPP17} for a formulation
via Morse indices $\ind(X_t)$=number of negative eigenvalues of $L$, 
counted with multiplicity, used to 
find bifurcation points in families of nodoids, which we shall
numerically corroborate in \S\ref{nodnumDBC}. An equivariant version
can be found in \cite[Theorem 5.4]{KPS18}, applied to bifurcations of
triply periodic minimal surfaces, 
for which linearizations always have a trivial 5 dimensional
kernel due to translations and rotations, see \S\ref{tpssec} for
numerical illustration. See also \cite{GoS2002,hoyle,kiel2012} and 
\cite[Chapters 2 and 3]{p2pbook} for 
general discussion of  Crandall--Rabinowitz 
type results, and  of Krasnoselski type results (odd multiplicity 
of critical eigenvalues, based on degree theory), 
including equivariant versions. 

b) Besides the (zero) DBCs \reff{DBC1} corresponding to a fixed
boundary $\pa X=\ga$, we shall consider so called free boundaries of 
Neumann type. This means that $\pa X\subset \Ga$, where
$\Ga\subset\R^3$ is a fixed 2D support manifold (e.g., a plane), and
that $X$ intersects $\Ga$ orthogonally. Following \cite{Ros08}, we
summarize the second derivative of $A$ in this case as follows: if $h_\Ga$ is
the second fundamental form of $\Ga$, and $\psi=uN$, then $N|_{\pa X}$
is tangent to $\Ga$, such that $h_{\Ga}(N,N)$ is well defined, and 
\hual{
\label{appb} 
A''(0)&= 
-\int_{X_0} (\Delta u +\|S_0\|^2u)u \dd S-\int_{\pa
X_0}h_\Ga(N,N)u^2\dd s. }
Note that the term $\int_{\pa X_0}\ldots$ in \reff{appb} vanishes 
if $\Ga$ is a plane and hence $h_\Ga\equiv 0$.  

c) The formulas \reff{f-area}--\reff{jac} translate to our 
discrete computational  setting in a straightforward way, see
\S\ref{ddgssec}. However, given some $X_0$, we compute $X=X_0+uN_0$ via
Newton loops for iterates $u_n$ with 
derivatives (of $V,A$ and $H$) 
evaluated at $X_n=X_0+u_n N_0$, and hence the formulas are 
accordingly adjusted.  
}\eex\erem 

\subsection{Default data and initialization of a \pdep\ struct {\tt p}}\label{dsec} 
Before explaining the modifications needed for the geometric problems, 
we briefly review the standard setup of \pdep, 
and as usual assume that all problem data is contained 
in the \mlab\ struct {\tt p} as in {\tt p}roblem. 
In the standard FEM setting this includes 
the object {\tt p.pdeo} (with sub--objects {\tt fem} and {\tt grid}), 
which provides methods to generate 
FEM meshes, to code BCs, and to assemble FEM matrices {\tt M} (mass matrix) 
and {\tt K} (e.g., Laplacian), or directly a rhs $G$. 
Typical initializations and first continuation steps 
in the FEM setup  (for semilinear problems) then run as follows, where steps 1,2 and 5 are usually combined into an init-function.
\bcen 
\item Call {\tt p=stanparam()} to initialize most fields in {\tt p} with 
default values (see source of {\tt stanparam.m} for default fields and values). 
\item Call a {\tt pdeo} constructor, for instance 
{\tt p.pdeo=stanpdeo2D(p,vararg)}, where here and 
in the following {\tt vararg} stands for variable arguments.  
\item In a function {\tt oosetfemops} (in the current directory), use 
{\tt p.pdeo.assema} to generate a mass matrix {\tt p.mat.M}  
and a stiffness matrix {\tt p.mat.K} (typically corresponding to $-\Delta$), 
and possibly further FEM matrices, e.g., for BCs. 
\item Use {\tt p.mat.M} and {\tt p.mat.K} in a function {\tt r=sG(p,u)} to 
encode the PDE, and optionally the Jacobian in {\tt Gu=sGjac(p,u)} 
(usually recommended, but numerical Jacobians are also supported). The input 
argument {\tt u} contains the ``PDE unknowns'' $u$ and the parameters 
appended at the end. If required by the problem, similarly 
create a function {\tt q=qf(p,u)} for the constraints as in (\ref{feq}b), 
and a function {\tt qu=qder(p,u)} for the derivatives of {\tt qf}. 
\item Initialize {\tt p.u} with a first solution (or a solution guess, 
to be corrected in a Newton loop). 
\item Call {\tt p=cont(p)} to (attempt to) continue the initial solution 
in some parameter, including bifurcation detection, localization, and 
saving to disk. 
\item Call {\tt p=swibra(dir,bpt,newdir)} to attempt branch switching 
at branch point {\tt bpt} in directory {\tt dir}; subsequently, call 
{\tt p=cont(p)} again, with saving in {\tt newdir}.
\item Perform further tasks such as fold or branch--point continuation;
  use {\tt plotbra(dir,pt,vararg)} to plot bifurcation diagrams, and
  {\tt plotsol(dir,pt,vararg)} to plot sample solutions. 
\ecen 
 
\brem\label{fuharem}{\rm The rhs, Jacobian, and some further 
functions needed/used to run \pdep\ on a problem {\tt p}, are 
interfaced via function handles in {\tt p.fuha}. For instance, 
you can give the function encoding your rhs $G$ any name, e.g., {\tt myrhs}, 
with signature {\tt res=myrhs(p,u)}, 
and then set {\tt p.fuha.sG=@myrhs}, but you can also 
simply keep the ``standard names'' {\tt sG} and {\tt sGjac} and 
encode these in the respective problem directory. 
For many handles in {\tt p.fuha} there are 
standard choices which we seldomly modify, e.g., 
{\tt p.fuha.headfu=@stanheadfu} (the header for printouts). 
 Functions for which the 
``default choice'' is more likely to be modified include, e.g.,%
\bci 
\item {\tt p.fuha.outfu=@stanbra}, signature 
{\tt out=stanbra(p,u)},  branch output; 
\item {\tt p.fuha.lss=@lss}, signature 
{\tt [x,p]=lss(A,b,p)}, linear systems solver $x=A^{-1}b$.\\  
Other options include, e.g., 
{\tt lssbel} (bordered elimination) and {\tt lssAMG} (preconditioned 
GMRES using {\tt ilupack} \cite{ilupack}). 
\eci 
During continuation, the current solution is plotted 
via {\tt plotsol(p)}, and similarly for a posteriori plotting 
(from disk). The behavior of {\tt plotsol} is 
controlled by the subfields of {\tt p.plot} (and possible auxiliary arguments), 
and if {\tt p.plot.pstyle=-1}, then {\tt plotsol} immediately calls 
a function {\tt userplot}, to be user--provided.  Such user functions 
naturally must be in the \mlab--path, 
typically in the current problem directory, which \mlab\ scans first 
when looking for a file. We sometimes also exploit this to overload 
\pdep\ library functions that need modifications for a given problem. 
\eex}\erem

\subsection{\pdep\ setup for discrete differential geometry}
\label{ddgsec}
\subsubsection{Discrete differential geometry FEM operators}\label{ddgssec}
\begin{figure}[ht] %H]
\centering
\btab{lll}{{\sm (a)}\\
\ig[width=0.32\tew]{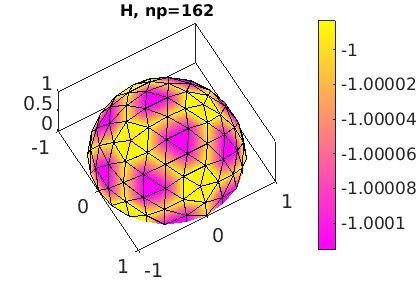}
&\ig[width=0.32\tew]{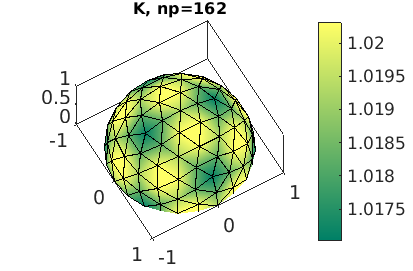}
&\quad \ig[width=0.32\tew]{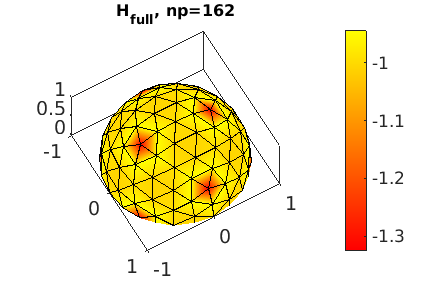}\\
{\sm (b)}\\[0mm]
\ig[width=0.32\tew]{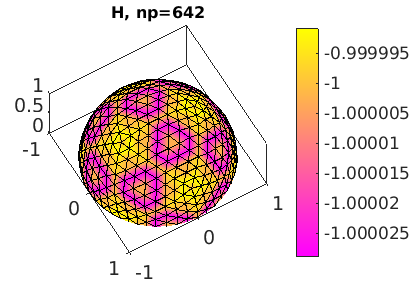}
&\ig[width=0.32\tew]{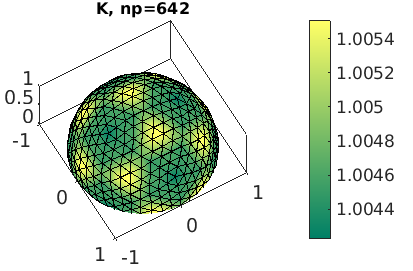}
&\quad \ig[width=0.32\tew]{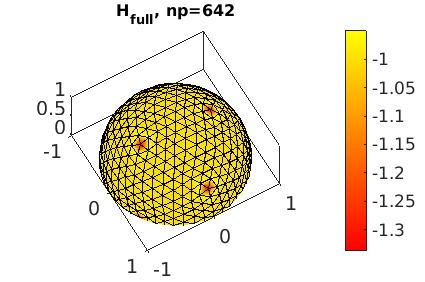}}
\vs{-3mm}
\caption{{\small Discrete $H$ (and $K$) on (coarse) 
meshes of the unit sphere (plots cropped). Two 
left columns: Convergence for $H=-\frac 1 2 M^{-1}\spr{L X,N}$ 
and $K=M^{-1}{\tt K}$ with Voronoi $M$. Right 
column: No convergence for $H$ (and similar for $K$) 
at valence 5 nodes when using $\Mf$. }\label{nf1}}
\end{figure}

We recall a few discrete differential geometry operators from 
\cite{MDSB03, Jdiss}, and shall use implementations of them 
from the \gptool\ \cite{gpgit}. 
 Given a triangulation 
${\tt X}\in\R^{n_p\times 3}$ (point coordinates) 
and ${\tt tri}\in \R^{n_t\times 3}$ (triangle corner indices) of $X$, 
and the piecewise linear element ``hat'' functions $\phi_i:X\to\R$, 
$\phi_i(X_j)=\del_{ij}$, we have 
\huga{\label{cotanL}
\int \nabla\phi_i\nabla\phi_j\dd S=-\frac 1 2 (\cot\al_{ij}+\cot\beta_{ij})=:L_{ij}, 
}
where $\al_{ij}$ and $\beta_{ij}$ are the angles opposite the edge 
$e_{ij}$ from point $X_i$ to point $X_j$. For $u:X\to\R$, 
$u=\sum_{i=1}^{n_p}u_i\phi_i$, this yields 
the FEM stiffness matrix $Lu$ corresponding to the Laplace--Beltrami operator 
$-\Delta u$ weighted 
by the mass matrix $M$. In \cite{MDSB03} it is explained that for 
geometric problems, with possibly rather distorted triangles, 
instead of the full mass matrix $\Mf$ with 
\huga{\Mf_{ij}=\int \phi_i\phi_j\dd S,} 
the Voronoi mass matrix 
\huga{\label{MV1} 
\MV=\diag(A_1,\ldots,A_{n_p}),
} 
should be expected to give better approximations, see also Fig.\,\ref{nf1}. 
Here, 
$A_i=\sum_{j=1}^{n_i}A_m(T_j)$ is the area of the Voronoi region at node $i$, 
where $T_j$, $j=1,\ldots,n_i$ are the adjacent triangles, 
and $A_m(T)$ is a ``mixed'' area: 
For non--obtuse $T$,  $A_m(T)$ is the area of the rhomb with corners 
in $X_i$, in the midpoints of the edges adjacent to $X_i$, 
and in the circumcenter of 
$T$, while for obtuse $T$ we let $A_m(T):=|T|/2$ if the angle at $X_i$ is 
obtuse, and $A_m(T):=|T|/4$ else. Alltogether, this yields the 
approximation 
\huga{\label{du1} 
-\Delta u=\MV^{-1}Lu, 
} 
where $\MV$ from \reff{MV1} is diagonal, and $L$ and $\MV$ are 
evaluated very efficiently via {\tt cotmatrix} and {\tt massmatrix} 
from the \gptool, see Table \ref{gptab}. 

However, as we always consider our problems such as \reff{feq} 
in weak form, we let ${\tt H}=-\frac 1 2\spr{L X,N}$, 
where for the vertex normals $N$ we can use 
{\tt per\_vertex\_normals}, and the weak form of, e.g., $H-H_0{=}0$ 
then is 
\huga{\label{cmcw}
-\spr{L X,N}-2MH_0=0, 
}
again with Voronoi $M$. 
Alternatively, we use {\tt [k,H,K,M]=discrete\_curvatures(X,tri)}, 
where {\tt K} and ${\tt k}=(k_1,k_2)$ 
are the (weighted, i.e., weak) discrete Gaussian and principal 
curvatures per vertex; 
these are computed from a discrete 
version of the Gau\ss--Bonnet theorem.%
\footnote{\label{GBfoot}On a manifold $X$ with boundary $\pa X$ we have 
$\int_X K\dd S+\int_{\pa X}\kap_g \dd s=2\pi\chi(X)$ 
where $\chi(X)$ is the Euler characteristic of $X$, 
and $\kap_g$ is the geodesic curvature of $\pa X$. This will play 
an important role for the biomembranes in \S\ref{wsec}. The discrete 
formula \reff{dKG1} is used at interior points of $X$, while at  
boundary points $X_i$ it is modified to $K(X_i)-\pi$. } 
Namely 
\huga{\label{dKG1}
{\tt K}(X_i)=2\pi-\sum_{j=1}^{n_i}\theta_j, 
\text{\quad (and $k_1=H+\sqrt{D}$ and $k_2=H-\sqrt{D}$),}
} 
where the $\theta_j$ are the 
angles at $X_i$, and where 
the discriminant $D=H^2-K$ (which is non-negative in the continuous case)  
in the discrete case is set to $0$ if negative. 
An approximations of 
$K$ is then obtained (cheaply, since $M$ is diagonal) from 
\huga{\label{GK}
K=\MV^{-1}{\tt K}. 
}

There are further schemes for $H$ and $K$, with different 
convergence behaviors, see \cite{xuxu09} and the references therein. 
Numerical experiments in 
\cite{xu04} show that a variety of natural schemes for $-\Delta$ 
in general do not converge, but that 
$M^{-1}L=-\Delta+\CO(h^2)$ with Voronoi $M$ at valence 6 nodes (6 neighbors) 
\cite[Theorem 2.1]{xu04}, 
where $h$ is a suitable triangle diameter.  
%(e.g., max edge length). 
In Fig.\,\ref{nf1} we give an illustration of the error and convergence 
behavior of our discrete $H=-\frac 1 2 M^{-1}\spr{L X,N}$ based on \reff{du1}, 
and of $K$ from \reff{GK} on (coarse) 
discretizations of the unit sphere, obtained via {\tt subdivided\_sphere} 
with 2 (a) resp.~3 (b) subdivisions, and with $N{=}$outer normal. 
Therefore, $H=-1$, $K=1$ are the exact values, 
and the two left columns indicate the convergence for $H$, but also that 
the node valence plays a role  on these 
otherwise very regular meshes.% 
\footnote{Euler's polyhedron formula yields that triangulations with 
all nodes of valence 6 do not exist, see, e.g., \cite{BoFi67}.}
However, the last column 
shows that using $\Mf$ in this example, i.e., 
$H_{{\rm  full}}=-\frac 1 2 \Mf^{-1}\spr{L X,N}$
gives a significant error (and similarly in $K$), 
and in fact no convergence 
at the valence 5 nodes;  see {\tt spheres/convtests.m} 
for the source of these experiments.

\subsubsection{The \pdep\ library {\tt Xcont}}\label{Xcsec}
Table \ref{gptab} lists the main (for us) 
functions from the {\tt gptoolbox} \cite{gpgit}, which 
we interface by functions from the  \pdep\ library {\tt Xcont}. 
The most important new data for continuation of a surface $X$ are 
{\tt p.X} and {\tt p.tri}, which essentially replace the data in 
{\tt p.pdeo.grid}. The most important switch, which also modifies 
the behavior of some standard \pdep\ functions, is 
\huga{\label{Xsw}
{\tt p.sw.Xcont}=\left\{\barr{ll} 0&\text{legacy setting (no {\tt X}),}\\
1&\text{switch on {\tt X}--continuation (default),}\\
2&\text{refined setting for {\tt X}--continuation.}
\earr\right.
}
The difference between {\tt p.sw.Xcont}=1 and {\tt p.sw.Xcont}=2 is as follows: 
For {\tt p.sw.Xcont}=1 we update {\tt p.X} after convergence of the Newton 
loop, i.e., set {\tt p.X=p.X+u*N0, p.up=u} (for plotting, see Remark \ref{Xcontrem}(b)), and 
{\tt u(1:p.np)=0} (zeroing out $u$ for the next continuation step). 
For {\tt p.sw.Xcont}=2 we do 
all this after each successful Newton-step, such that we obtain slightly 
different (more accurate) Jacobians as we also have a new $N$. 
A side-effect is that 
the last update in the Newton loop and hence the {\tt p.up} 
may be very small and will in 
general not represent the ``direction'' of continuation. 
For {\tt p.sw.Xcont}=1, {\tt p.up} will contain the complete step, 
and therefore we choose this as default. % in {\tt stanparamX}.  
For the spectral computations (bifurcation detection and localization) 
{\tt p.sw.Xcont}=1 vs {\tt p.sw.Xcont}=2 makes no difference as all 
spectral computations are done after convergence of the Newton loops, 
i.e., after the final update of $\tt p.X$. 
\taskip
\begin{table}[ht]\caption{
Main functions used from \cite{gpgit} and \cite{afemgit,FS21}, 
and some small additions/mods. 
 Here ${\tt X}{\in}\R^{n_p\times 3}$ are the point coordinates of the triangulation, and ${\tt tri}\in\R^{3\times n_t}$ the triangulation as rows of point numbers of triangles.
\label{gptab}}
\centering\vs{-2mm}
{\small 
\begin{tabular}{p{0.42\tew}|p{0.55\tew}} 
function&remarks\\
\hline
{\tt N=per\_vertex\_normals(X,tri)}&normals; should be interfaced as {\tt N=getN(p,X)}, see Table \ref{Xconttab} to possibly flip sign at one place;\\ 
{\tt L=cotmatrix(X,tri)}&cotangent discrete Laplace--Beltrami\\
{\tt M=massmatrix(X,tri,type)}&mass matrix, with type 'full', 'barycentric' or 
'voronoi'. \\
{\tt K=discrete\_gaussian\_curvature(X,tri)}&Gaussian curvature $K$. \\
{\tt [k,H,K,M]=discrete\_curvatures(X,tri)}&principal curvatures $k=(k_1,k_2)$, 
and $H$, $K$, $M$.\\
\hline
{\tt [X,tri,DBC,NBC]}&refinement, with DBC/NBC the Dirichlet/Neumann\\
{\tt =TrefineRGB(X,tri,DBC,NBC,elist)}&  boundary nodes, and elist the triangles to refine. Interfaced via {\tt p=refineX(p,sig);} where {\tt sig} is the usual 
factor of triangles to refine, and where after refinement {\tt p.u} and 
{\tt p.tau} are interpolated to the new mesh. See also 
{\tt TcoarsenRGB}, with similar signature.\\
\end{tabular}
}
\end{table}
\teskip

\taskip
\begin{table}[ht]\caption{
Main functions from the library {\tt Xcont}, which collects 
interfaces to functions from {\tt gptoolbox} and \cite{afemgit}, 
and some additions/mods; see sources for argument lists and detailed 
comments, and Remark \ref{Xcontrem} for further comments..\\
\label{Xconttab}}
\centering\vs{-2mm}
{\small 
\begin{tabular}{p{0.17\tew}|p{0.83\tew}} 
function&remarks\\
\hline
{\tt p=stanparamX(p)}&convenience function to set \pdep\ parameters to ``standard'' 
values for $X$--continuation; to be called {\em after} {\tt p=stanparam()} 
during initialization. \\
\hline 
{\tt getN}&interface to {\tt per\_vertex\_normals}; to flip orientation, 
make local copy with\\
&{\tt function N=getN(p,X);  N=-per\_vertex\_normals(X,p.tri); end} \\
{\tt getA, getV}&get area $A$ and volume $V$ according to \reff{area} and 
\reff{vol}. \\
{\tt qV, qjacV}&$V$ as a constraint $q$ for continuation, and its derivative 
$\pa_u q$. \\
{\tt qA, qjacA}&$A$ as a constraint $q$ for continuation, and its derivative 
$\pa_u q$.\\
{\tt c2P}&triangle center to point (nodes) interpolation matrix. \\
{\tt cmcbra}&default branch output in our demos below. \\
\hline
{\tt pplot}&plot of {\tt p.X}, usually called from {\tt userplot}; 
see Remark \ref{Xcontrem}b).  \\
{\tt oosetfemops}&usually only needed as a dummy (since called by \pdep\ 
library functions).\\
{\tt updX}&update {\tt p.X} after successful Newton steps.\\
{\tt plotHK}&given {\tt p.X}, plot $H(X)$ and $K(X)$ (mostly for checking).\\ 
\hline
{\tt e2rsA, e2rsAi}&element--to--refine--selector choosing triangles with large 
areas;  {\tt e2rsAi} chooses triangles with small areas, 
intended for coarsening. \\
{\tt e2rsshape1}&Select triangles according to \reff{delmdef}. See also {\tt e2rsshape*}.\\
{\tt meshqdat}&mesh quality data, see \reff{brd}. \\
{\tt refineX}&mesh refinement of {\tt p.X}; selecting triangles 
via {\tt p.fuha.e2rs} and calling {\tt TrefineRGB} (or {\tt Trefinelong} 
if {\tt p.sw.rlong==1} for bisection of only longest edges of triangles).  \\
{\tt coarsenX}&coarsening of {\tt p.X}, same syntax as {\tt refineX}; 
i.e., selecting triangles via {\tt p.fuha.e2rs}. Inverse to {\tt refineX}, 
as only triangles from prior refinement can be coarsened (to their common 
ancestors).    \\
{\tt degcoarsenX}&coarsening of mesh by removing degenerate triangles 
(via \gptool).\\
{\tt retrigX}&retriangulation of $X$, based on \cite{PeSt04}; 
using the adjacency in {\tt p.trig} to generate a new (Delauney) 
triangulation of $X$ while keeping the surface structure. \\
{\tt moveX}&retriangulate and move points of $X$ based on \cite{PeSt04} and the 
associated code.\\
\hline
{\tt geomflow}&simple explicit time integrator for $\dot X{=}fN$, 
where $f{=}{\tt p.fuha.flowf}$, e.g., $f(X){=}-H$ for mean curvature flow, implemented in {\tt mcff} 
(with DBCs). 
\end{tabular}
}
\end{table}
\teskip
Some main functions from the \pdep\ library {\tt Xcont} are listed 
in Table \ref{Xconttab}, and can be grouped as follows: 
\bcen 
\item Functions directly interfacing the {\tt gptoolbox} such 
as {\tt N=getN(p)} which in the default setting just calls 
{\tt N=per\_vertex\_normals(p.X,p.tri)}.  
%and {\tt M=massmatrix(p.X,p.tri,'Voronoi')}, respectively.  
These are meant as easy--to--change interfaces, as for a given problem 
it may be desired to, e.g., change the orientation of $X$. 
For this, make a local copy of {\tt getN.m} %in the problem directory 
and there set {\tt N=-per\_vertex\_normals(p.X,p.tri)}. 
\item  ``Template'' functions to compute $V(u)$ (corresponding to 
$V(X)=V(X_0+uN_0)$ in \reff{vol}),  
and $A(u)$ from \reff{area}, and wrappers to use them in constraints 
such as ${\tt q=qfV(p,u)}=V(u)-V_0$, and their derivative, 
such as ${\tt qu=qjacV(p,u)}=\pa_u V$. 
(Since $V$ is a scalar function, here we use 
``{\tt jac}'' in a loose sense.)  
\item Template convenience functions such as {\tt pplot} for plotting 
{\tt p.X}, and {\tt updX} for updating $X=X_0+uN_0$, cf.~\reff{n1}, 
and overloads of \pdep\ library function adapted 
to $X$ continuation, e.g., {\tt getGupde}, and 
a prototype {\tt oosetfemops}, which for $X$--continuation 
usually should only serve as a dummy needed by other \pdep\ functions. 
\item Functions related to mesh handling such as 
{\tt refineX}, {\tt coarsenX} and {\tt degcoarsenX}, 
and associated elements--to--refine--selectors such as 
{\tt e2rsA} (based on triangle areas) and {\tt e2rshape} (based on triangle 
shapes). 
Additionally there are {\tt retrigX} and {\tt moveX}, based on \cite{PeSt04}, 
see Remark \ref{Xcontrem2}. 
\item The template {\tt geomflow} for geometric flows such as 
MCF, interfaced in {\tt p.fuha.flowf}. 
\ecen 
All of these functions, in particular from 2.~and 3.~are ``templates''; 
for a given problem it may be necessary or useful 
to make local copies of some of these functions in the problem directory and 
change them there. Additionally, there is convenience function 
{\tt p=stanparamX(p)}, which (re)sets some \pdep\ parameters to 
typical values for $X$--continuation. Another important interface 
to the {\tt gptoolbox} is 
\bci
\item[6.] {\tt getM(p)}, which for scalar problems should call 
{\tt M=massmatrix(p.X,p.tri,'Voronoi')}, and for vector valued 
problems should build the system mass matrix 
from such ``scalar'' {\tt M},  cf.~\S\ref{wsec}. 
However, for compatibility with the non--$X$--setting, {\tt getM} is {\em not} 
part of {\tt Xcont}, but needs a local copy in each ($X$--continuation) 
problem directory. 
\eci

\brem{\rm\label{Xcontrem} 
a) As is, {\tt out=cmcbra(p,u)} puts the data 
\huga{\label{cmco} 
\text{out=[pars;V;A;E;meshq]}
} 
into {\tt out}, where {\tt pars} of length $m={\tt length(u)-p.nu}$ is the 
(user defined) parameter vector of the problem, $V$ and $A$ are 
volume and area of $X$, and $E=A+{\tt par(1)}V$, cf.\,\reff{cmcen}, 
which assumes that $H_0$ is at {\tt par(1)}, and is an {\em active} 
continuation parameter. 
Next, ${\tt meshq}=(\delm,a_{\max},a_{\min},h_{\max},h_{\min})$ 
gives measures of the mesh, computed in {\tt meshq=meshqdat(p)}. Here, 
 \huga{\label{delmdef}
\delm=\max(h/r)=\text{mesh distortion}, 
}
where, for each triangle, $h$ is the maximum edge--length, 
and $r$ the in--radius; $a_{\max}$ and $a_{\min}$ are 
the max and min of the triangle areas, and $h_{\max,\min}$ are the largest 
and smallest edge lengths in the triangulation. 
Thus, $\delm$ is our default mesh--distortion measure 
(one of several possible), see, 
e.g., Fig.\,\ref{spcf1c} and Fig.\,\ref{nodf2b} for plots of 
$\delm$ along branches.
For an equilateral triangle
$\delm=6/\sqrt{3}\approx 3.46$, and our experiments 
suggest that, as a rule of thump, triangles with 
$\delm<50$ are (still) of ``reasonable shape''. 

{\tt out=p.fuha.outfu} is appended to 
{\tt bradat(p,u)=[count;type;ineg;lam;err;L2]}%
\footnote{six values, where 
{\tt count=}step-counter, ${\tt type}\in\{0,1,2,3\}$ (regular point, 
BP, FP, Hopf point), {\tt ineg}=number of unstable eigenvalues (if 
${\tt p.sw.spcalc}>0$), {\tt lam}=value of primary continuation parameter, 
and  {\tt err} and {\tt L2} are {\em not} meaningful in the {\tt Xcont} setting, 
see {\tt bradat.m} for details)}, and %for convenience 
we list the component {\tt c} for branch plotting for {\tt p.fuha.outfu=@cmcbra} 
 (with m the number of parameters): 
{%\small
\huga{\label{brd} 
\btab{c|c|ccc|ccccc}{
c&1$\ldots$m&m+1&m+2&m+3&m+4&m+5&m+6&m+7&m+8\\
meaning&{\tt pars}&$V$&$A$&$E$&$\delm$&$a_{\max}$&
$a_{\min}$&$h_{\max}$&$h_{\min}$
}.
}
}%
Thus, to plot, e.g., $V$ over the {\em active} bifurcation parameter 
from a computed branch {\tt b1} into figure {\tt fnr}, 
use {\tt plotbra(b1,fnr,m+1,varargin)}, where {\tt varargin} 
gives {\tt many} options for colors, labels, etc. 
Similarly, to plot $\delm$, use {\tt plotbra(b1,fnr,m+4,varargin)}. 
%, see {\tt plotbra.m} for detailed documentation.%
 Moreover, 
{\tt c} in {\tt plotbra(b1,fnr,c,varargin)} can be a vector, and to, e.g., 
plot $\delm$ over $V$ use {\tt c=[m+1, m+4]}. 
%To plot data from {\tt bradat(p,u)} use 
%negative $c$, resp.~$c=0$ for {\tt L2}.}
%See Fig.\,\ref{nodf2b} for 
%some plots of $q_1$ in an experiment dealing with triangle shapes. 

b) {\tt pplot(p,fignr)} (or {\tt pplot(dir,pt,fignr)}, where {\tt p} is loaded from {\tt dir/pt}) colors {\tt p.X} by 
{\tt p.up}, which contains $u$ from the last continuation step and hence 
codes (together with $N$) the ``continuation direction'' 
(if {\tt p.sw.Xcont=1}), or just the last Newton update (if {\tt p.sw.Xcont=2}), 
see \reff{Xsw}. The behavior of {\tt pplot} can further be 
controlled by settings in 
the {\em global} structure {\tt p2pglob}.%
\footnote{This turned out to be convenient, i.e.:  When plotting from 
file we very often want to change the look of plots ``globally'', i.e., 
without first loading the point and then adapting settings.}
For instance, {\tt p2pglob.tsw} 
controls the titles, e.g., $(V,H)$ for {\tt tsw=1}, and 
$(A,H)$ for {\tt tsw=2}, which assumes that the parameters are ordered as 
$(H,V,A)$. For flexibility, if ${\tt tsw}>4$, then {\tt pplot} searches 
the current directory for a function {\tt mytitle} to generate a title, 
see, e.g., demo {\tt geomtut/enneper}. 
Again, {\tt pplot} is a template, if necessary 
to be modified in the problem directory for customized plots, 
but all demos from {\tt geomtut/} only use the given library version. 

c) In all demos below we use indices {\tt p.idx} and {\tt p.idN} 
of boundary points to set boundary conditions. In this, {\tt p.idx} 
should be thought as points for Dirichlet BCs and associated 
to  {\tt p.DBC} (the corresponding {\em edges}, 
updated in {\tt refineX} and {\tt coarsenX}), 
and {\tt p.idN} as Neumann BCs with edges {\tt p.NBC}. All of these 
can also be empty (for instance if $X$ is closed, or has only one 
type of boundary), and again, the use of {\tt p.idx, p.idN, p.DBC} and 
{\tt p.NBC} should only be seen as a template, for instance 
to be modified if a given problem has several different boundaries.
}\eex\erem 

\brem{\rm\label{Xcontrem2} 
a) For surface meshes ${\tt (X,tri)}$, 
mesh adaptation, i.e., refinement and coarsening,  
seems even more important than for standard (non--parametric) 
problems, because well behaved initial triangulations 
(well shaped triangles of roughly equal size) 
may  deteriorate as $X$ changes. 
The case of growing spherical caps in Fig.\,\ref{f1}(a) is rather harmless as 
triangle sizes grow but shapes stay intact, 
and can easily be dealt with by refinement of the 
largest triangles. For this, in {\tt e2rsA} we simply order the $n_t$  
triangles of {\tt tri} 
by decreasing size, and from these choose the first $\lfloor \sig n_t\rfloor$ 
for refinement by {\tt refineX}, i.e., we generally 
use $\sig={\tt p.nc.sig}$ as the parameter for the fraction of triangles 
to refine.  The refinement can be either done as RGB if {\tt p.sw.rlong=0}, 
or by refining only the longest edges of the 
selected triangles  if {\tt p.sw.rlong=1}. RGB is generally better if triangle shapes 
are crucial, but may result in rather long cascades to avoid hanging nodes 
(such that $\sig$ is only a lower bound for the fraction of triangles 
actually refined). 
Refine-long 
gives more control as {\em only} the selected triangles are bisected 
(plus at most one more triangle for each one selected), but 
may lead to obtuse triangles, and it seems that 
as for standard FEM very obtuse triangles are more dangerous than 
very acute triangles. A short computation shows that, e.g., 
for a right--angled 
triangle refine--long increases $\del=h/r_{\text{in}}$ by 45\%; however, 
this can often be repaired by combining refine--long with {\tt retrigX}, 
see b). See also \cite{Shew02} for a very useful discussion of mesh quality 
(in the planar setting, and in 3D).
%, and, e.g., \cite{...} for a thorough review of surface (re)meshing. 

Conversely, {\tt coarsenX} can be used to coarsen previously refined 
triangles, again from a list generated by {\tt p.fuha.e2rs}, 
which should be reset from the one chosen for refinement. For instance, 
{\tt e2rsAi} selects the $\lfloor \sig n_t\rfloor$ 
triangles of {\em smallest} area, but 
these have to be from the list of previously refined triangles. 

{\tt degcoarsenX} works differently: It calls the modification {\tt rmdegfaces} 
of the  \gptool\ function 
{\tt rm\_degenerate\_faces}, and essentially 
aims to remove obtuse and acute triangles by collapsing (short) edges. 
This works in many cases but may result in hanging nodes such that 
the FEM no longer works. 

Both, {\tt refineX} and {\tt degcoarsenX} can be told to {\em not} 
refine/coarsen boundary triangles, which is crucial for the case 
of periodic BCs. 

b) We also provide two small modifications of (actually interfaces to) 
code from \cite{PeSt04}. In {\tt retrigX.m} we generate a new 
(Delauney) triangulation of $X$, keeping intact the 
surface structure of $X$. This is in particular 
useful if $X$ has been obtained from {\em long} refinement, which 
typically results in nodes having $8$ adjacent triangles (valence $8$), while 
``standard'' triangulations (and the output of {\tt retrigX}) 
have valence $5$ and $6$, which generally seems to result in more 
robust continuations. In {\tt moveX} we combine {\tt retrigX} with 
motion of the points in $X$ due to ``truss forces'' of the 
triangulation, aimed at more uniform edge lengths. 
Due to the similarity of the triangulation truss forces and surface 
tension, this works best for minimal surfaces ($H{=}0$), or otherwise 
for surfaces with small $|H|$. 
}
\eex\erem

\section{Example implementations and results}\label{exsec}
Our demos are meant to show how to set up different geometric bifurcation 
problems, in particular with  different BCs. 
They mostly deal with classical minimal or more 
generally CMC surfaces, for instance the  Enneper and Schwarz--P 
surfaces, and so called nodoids (including physically relevant 
liquid bridges). Many demos start with CMC surfaces of revolution, 
%which were classified by Delauney and 
%consist of five categories: %(rotationally symmetric) planar domains, spheres, 
%cylinders, catenoids, undoloids, and nodoids. Our 
and our main interest then 
are bifurcations  breaking the rotational symmetry. The 
minimal surfaces in \S\ref{obcsec} are motivated by Plateau's problem 
of soap films spanning a given wire, 
and in \S\ref{wsec} we consider 4th order problems obtained 
from the Helfrich functional. 

All demos come with a number of function files namely (at least, 
with * a placeholder, usually to be replaced by a short problem name, 
cf.~\S\ref{ddgsec}): 
{\tt sG*.m} describing the rhs of the problem; {\tt *init.m} for 
initialization; {\tt userplot.m} and {\tt getM.m} for technical 
reasons (downward compatibility). Additionally, in some demos 
we overload functions from {\tt libs/Xcont}, e.g., {\tt cmcbra.m} 
for branch output. Finally, there are 
script files {\tt cmds*.m} with * a number if there is more than one 
script. In our descriptions of the first demos, we give tables listing 
the used files and their purpose (starting with the scripts), and 
 we give a few listings of (parts of) pertinent files to discuss 
important points. This becomes less for the later more advanced 
demos, for which we rather put more comments into the m--files themselves. 

\subsection{Spherical caps}\label{scsec}
\def\dhome{geomtut/spcap1} 
%\subsubsection{Dirichlet BCs at the unit circle}\label{dbcsec}
We start with the continuation in volume $V$ 
%(or area $A$) 
of spherical caps over the unit circle $\ga$ 
in the $x$--$y$ plane, 
as previewed in Fig.\,\ref{f1}(a). 
It is known  \cite{ALP99},\,\cite[\S2.6]{KPP14} that no bifurcations 
occur, and hence this only serves as an introductory toy model. 
Table \ref{sctab2} gives an overview of the used files, and Listings 
\ref{spl0}--\ref{spl3} 
show the initialization {\tt scinit.m}, the rhs {\tt sGsc.m}, 
and  the first script {\tt cmds1.m}. 
The BCs are $\pa X=\ga=\{(x,y,0)\in\R^3: x^2+y^2=1\}$, which since 
they hold for the initial unit disk translate into 
\hual{\label{bc1}
&\text{$u|_{\pa X}=0$.}
}

\brem\label{screm2}{\rm a) We can as well continue 
directly in $H$, without any constraints, and starting from 
the disk again obtain the same branch (see lines 14,15 of {\tt cmds1.m}). 
Our setup in {\tt cmds1.m} is motivated 
by applications, where typically the volume is the 
external parameter. I.e., the setup is a template how to use 
the volume ({\tt qfV}) or area ({\tt qfA}) constraints, together 
with the derivatives ({\tt qjacV} and {\tt qjacA}). 
Note that {\tt p.nc.ilam=[2,1]} for using 
$V$ as the {\em primary active} parameter in {\tt cmds1}, 
and {\tt p.nc.ilam=[3,1]} for using $A$, while $H={\tt par(1)}$ is 
a {\em secondary} active parameter in both cases. 

b) Only the {\em active} continuation parameters are updated in {\tt p.u}; 
thus, when continuing only in $H$, say, then to plot, e.g., $A$ over $H$ 
we cannot choose {\tt p.plot.bpcmp}=3 (the parameter index of $A$), 
but must take {\tt p.plot.bpcmp}=3+2=5. This is because the computed $A$ is 
put second after the parameters in the output function {\tt cmcbra}, and here we have 
three parameters $(H,V,A)$. But again, $A={\tt par(3)}$ is only updated  
if $3\in {\tt p.nc.ilam}$, i.e., if $A$ is an active parameter. 
} 
\eex\erem 

\taskip
\begin{table}[ht]\caption{
Files in {\tt pde2path/demos/geomtut/spcap1}; the last two are typical examples 
of (small) local mods of library functions. 
\label{sctab2}}
\centering\vs{-2mm}
{\small 
\begin{tabular}{p{0.2\tew}|p{0.79\tew}} 
%file&remarks\\ \hline
{\tt cmds1.m}&continuation in $(V,H)$ and in $(A,H)$, respectively.\\
{\tt cmds2.m}, {\tt cmds3.m}&tests of different mesh refinement options, 
and MCF tests.\\
%{\tt cmds3.m}&MCF tests.\\ 
{\tt getM.m}&standard (Voronoi) mass matrix.\\
{\tt scinit.m}&Init, data stored in {\tt p.u} (including 
computed $H, A$ and $V$), and in {\tt p.X} and {\tt p.tri}.\\
{\tt sGsc.m, scjac.m}&rhs based on \reff{cotanL}, and Jacobian based on \reff{f-mean}.\\
{\tt mcff.m}&mean curvature flow rhs $f$; problem dependent via choice of {\tt getN}. \\
\hline
{\tt cmcbra.m}&local copy and mod of library function {\tt cmcbra.m} to put 
error $e(X)$ \reff{sced} on branch.\\
{\tt refufu.m}&local copy and mod (and renaming) of {\tt stanufu.m} to do 
adaptive mesh refinement based on $e(X)$; ``switched on'' by setting 
{\tt p.fuha.ufu=@refufu}. \\
{\tt coarsufu.m}&similar to {\tt refufu.m}, used for mesh coarsening 
of decreasing caps; % ``switched on'' by setting {\tt p.fuha.ufu=@coarsufu}. 
\end{tabular}
}
\end{table}
\teskip

\hulst{caption={{\small  {\tt spcap1/scinit.m}; the 
pde-object {\tt pde} in line 4 is generated in a legacy \pdep\ setup 
but only used to generate the initial {\tt p.X}, with initial 
triangulation stored in {\tt p.tri} (line 10).}},label=spl0}{\dhome/scinit.m}

During init, we call {\tt pde=diskpdeo2} to generate a 
{\rm temporary} FEM object from which we extract the initial mesh 
to generate the initial {\tt p.X} and store {\tt p.tri} 
as the triangulation. Additionally we extract 
{\tt p.DBC} as the (Dirichlet) boundary {\em edge} index vectors, 
and {\tt p.idx} as the boundary {\em point} indices. This is 
in principle redundant, but it makes the setup of the DBCs in 
{\tt sGsc} shorter.  

\hulst{stepnumber=10}{\dhome/sGsc.m}
\hulst{caption={{\small {\tt spcap1/sGsc.m}, and 
start of {\tt cmds1.m} (omitting plotting).}},
linerange=3-17,label=spl3,firstline=3}{\dhome/cmds1.m} 

In {\tt cmds1.m} we then continue the initial disk (with $V=0$ and $A=\pi$) 
in $V$. For this we switch on the constraint $V(u)=V$ via 
{\tt p.nc.nq=1}, {\tt p.fuha.qf=@qfV; p.fuha.qfder=@qjacV}, with the {\tt Xcont} 
library functions {\tt qfV} and {\tt qjacV}, and set {\tt p.nc.ilam=[2,1]}, 
cf.~Remark \ref{screm2}.  
For mesh adaptation we use the triangle areas on $X$ 
as selector, and {\tt refineX} also updates {\tt p.DBC} and {\tt p.idx}, 
leading to Fig.\,\ref{f1}(a), where we use repeated mesh refinement 
every 5th step. This way we can accurately continue to arbitrary 
large $V$, i.e., arbitrary large ``cap radius'' $R$,  
where $H=1/R$ asymptotes to $H=0$. 
In the second part of {\tt cmds1.m} we continue in {\tt $A$} and hence set 
{\tt p.nc.ilam=[3,1]; p.fuha.qf=@qfA;} and {\tt p.fuha.qfder=@qjacA}. 
This yields exactly the same branch as the continuation in $V$, 
and all this works very robustly and fast. 

\brem\label{jacrem}{\rm 
The  numerical Jacobians of $G$ (for {\tt p.sw.jac=0} in line 5 of 
Listing \ref{spl3}) are sufficiently fast to not play a role for 
the speed of the continuation, at least for $n_p<2000$, say, 
because \mlab's {\tt numjac} can efficiently exploit the 
known sparsity (structure) of $\pa_u G$, given by the sparsity 
structure of the Laplacian $K$, or equivalently, by the 
sparsity structure of the (full, not Voronoi) mass matrix $M$. 
On the other hand, if  $q$ implements 
some integral constraints, e.g., area or volume, 
then $\pa_u q(u)\in\R^{n_q\times n_p}$ is dense, 
and numerical derivatives for $\pa_u q$ are a serious bottleneck. 
For illustration, in {\tt cmds1.m} we use the commands 
{\tt jaccheck} and {\tt qjaccheck}, which are rather important 
for ``debugging'' when numerical Jacobians become too slow. 
Both return relative errors between functional and numerical 
Jacobians, and as a rule of thump,  in $\pa_u G$ 
relative errors $\le 10^{-3}$ should be achieved, 
and do not affect the continuation or bifurcation results,
and for $\pa_u q$ even somewhat larger relative errors are usually no problem. 
}\eex\erem 

\begin{figure}[ht] 
\centering
\btab{lll}{{\sm (a)}&{\sm (b)}&{\sm (c)}\\
\hs{-2mm}\rb{2mm}{\ig[height=0.29\tew]{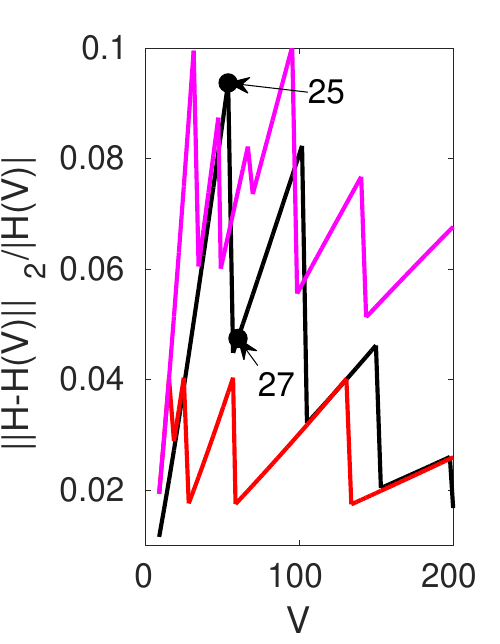}}
&\hs{-2mm}\rb{2mm}{\ig[height=0.29\tew]{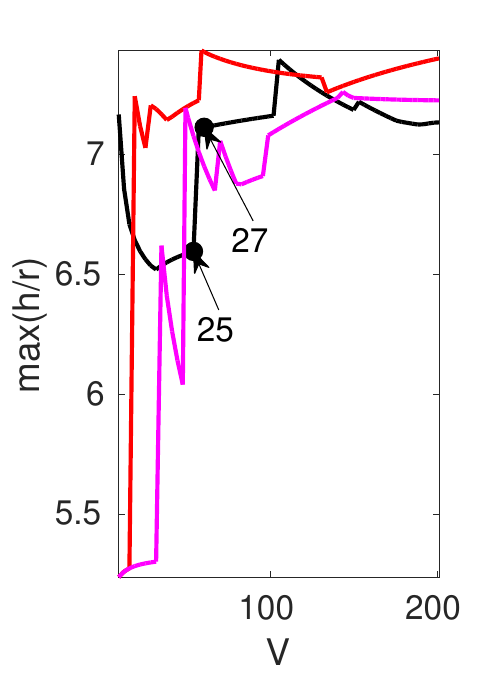}}
&\hs{-5mm}\ig[height=0.29\tew]{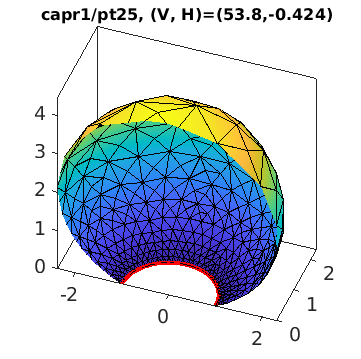}
\hs{0mm}\ig[height=0.29\tew]{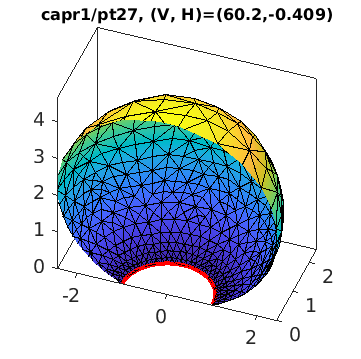}
%\hs{-3mm}\ig[height=0.33\tew]{p4/132}
}
\vs{-4mm}
\caption{{\small Results from {\tt spcap1/cmds2.m}. (a) 
Error $e(X):=\|H-H(V)\|_2/|H(V)|$ for refinement each 15th step 
({\tt capr1}, black) (starting at step 10), 
when $e(X)>{\tt p.nc.errtol}=0.05$, 
using {\tt p.fuha.ufu=@refufu} ({\tt capr3}, red), and when 
$\max(A)>0.3$ using  {\tt p.fuha.ufu=@refufumaxA} with $\sig=0.3$ 
({\tt capr4}, magenta). At $V=200$, $n_p=1452$ on {\tt capr1}, 
$n_p=1486$ on {\tt capr3}, and $n_p=636$ on {\tt capr4}. 
(b) Mesh distortion $\delm=\max(h/r)$ (edge-length over in--radius). 
(c) Illustration of meshes before/after refinement at pt25;  plots cropped 
at $y=0$ for better visibility of the meshes, and the boundary 
at $z=0$ marked in red. 
\label{spcf1c}}}
\end{figure}

In {\tt cmds2.m} we test different options for mesh adaptation, see 
 Listing \ref{spcl40} and Fig.\,\ref{spcf1c}.%
\footnote{Fig.\,\ref{spcf1c}(a,b) shows essentially verbatim output 
from {\tt plotbra} in {\tt cmds2.m}, where the dots and numbers 
indicate the continuation step, subsequently used also in the sample plots 
as in (c). This also holds for all subsequent plots, and the only 
``manual adjustments'' are the occasional repositioning of the 
numbers at the arrows by drag and drop, as this is not automatically 
optimized.} 
The black line  {\tt capr1} in (a), starting from {\tt cap1/pt10}, 
corresponds to adaptation each $15$th step, with 
``refinement factor'' $\sig=0.3$ (fraction of triangles marked for refinement). 
As we choose {\tt p.sw.rlong=1} we only bisect the 
longest edge of a selected triangle, and the actual 
fraction of refined triangles is between $\sig$ and $2\sig$.%
\footnote{For {\tt p.sw.rlong=0} (RGB refinement with possibly longer 
cascades to avoid hanging nodes) $\sig$ is only a lower bound.} 
For {\tt capr3} (red) 
we refine when the ``error'' $e(X)$ exceeds 
${\tt p.nc.errbound}=0.04$, where $e(X)$ is also used for plotting 
and defined as follows: For given $V$ 
we compute the (exact) $H(V)$ of the 
associated (exact) 
spherical cap $C(V)$ as %\huc{formel kurz begründen? aus Brubaker?} 
\huga{
H(V)=-\frac{\pi^{1/3}(3V+s-\pi^{2/3})(s-3V)^{1/3}}{s(3V+s)^{1/3}}, \quad 
\text{where}\quad s=\sqrt{9V^2+\pi^2}. 
}
We then define the ``relative $L^2$ error'' 
\huga{\label{sced}
e(X)=\|H(X)-H(V)\|_{L^2(X)}/|H(V)|, 
}
and put $e(X)$ on the branch in the modified local copy {\tt cmcbra.m} of the 
standard (library) {\tt cmcbra.m}.%
\footnote{We also put $\|H(X)-H(V)\|_{\infty}/|H(V)|$ and 
$z_{\max}(V)-z_{\max}$ on the branch, 
where $z_{\max}(V)$ is the height of $C(V)$ 
and $z_{\max}$ the numerical height; 
these can then also be plotted via {\tt plotbra}, and/or chosen as 
error indicators, but 
the $L^2$ error seems most natural. Also note that $e(X)$ 
is normalized by $|H(V)|$ (which decays in $V$), but not by 
$A$ (which increases with $V$).} 
$e(X)$ can then be plotted like any other output variable, and, moreover, 
can be used (without recomputing) in {\tt p.fuha.ufu} (user function), 
which is called after each successful continuation step. The default 
(library) setting {\tt p.fuha.ufu=@stanufu} essentially only 
gives printout, and to switch on the adaptive meshing we rename and 
modify a local copy as {\tt refufu.m}, and set {\tt p.fuha.ufu=@refufu}. 
Since $e(X)$ is at position 13 in (our modified) {\tt out=cmcbra(p,u)}, and  
since {\tt out} is appended to the six values from {\tt bradat}, 
cf.~Remark \ref{Xcontrem}a), in {\tt refufu.m} we then simply add the commands
\hugast{{\sm 
\text{{\tt if brout(6+13)$>$p.nc.errbound; p=refineX(p,p.nc.sig); end}. }}
}

Another ``natural'' alternative is to refine when 
\huga{\label{maxAe}
a_{\max}=\max(a_1,\ldots,a_{nt})>{\tt p.maxA}, 
}
i.e., when the maximum area of the {\tt nt} triangles 
exceeds a chosen bound. This 
is {\em not} an error estimator in any sense (as a plane can 
be discretized by arbitrary large triangles), but an ad hoc 
criterion, with typically an ad hoc choice of {\tt p.maxA}, which 
could be correlated to $H$. It is implemented in 
{\tt refufumaxA.m} which, if $\max A>{\tt p.maxA}$, calls 
{\tt refineX} with 
{\tt e2rsmaxA} to select all triangles with $A>(1-\sig){\tt p.maxA}$. 
With ${\tt p.maxA}=0.3$ and $\sig=0.2$ this yields the magenta line in 
Fig.\ref{spcf1c}(a). 
 
% The difference between {\tt capr1} ($k=30$, black), and {\tt capr2} 
% ($k=20$, blue) in (a) is obvious, i.e., the more frequent 
% refinement yields more frequent drops of $e(X)$ and, at $V=200$, 
% $n_p=1288$ mesh points on {\tt capr1}, but $n_p=1923$ on {\tt capr2}. 
% Visually, however, the solutions on the branches are the same.  
The samples in Fig.\ref{spcf1c}(c) illustrate a refinement step 
on the black branch, yielding a ``reasonable'' mesh 
also at large $V$. However, this naturally depends on 
the choice of steps between refinements 
(and on the refinement fraction {\tt sig} and 
continuation stepsize {\tt ds}).  
For the red line in Fig.\,\ref{spcf1c}(a), the refinement when 
the error $e(X)$ exceeds the chosen bound {\tt p.nc.errbound} 
is more genuinely adaptive, and this similarly holds for {\tt capr4} 
based on \reff{maxAe}, see also {\tt cmds2.m} for various further plots. 
(b) shows that the long--refinement generally yields a (mild) increase 
of the mesh distortion $\delm$, but overall the mesh--quality stays very good.  

\hulst{caption={{\small Selection from {\tt spcap1/cmds2.m}, 
refinement each 15th step, $e(X)$--dependent refinement 
via setting {\tt p.fuha.ufu=@refufu} and {\tt p.nc.errbound=0.04}, 
and refinement based on \reff{maxAe}.}},
linerange=4-16,label=spcl40,firstnumber=4}{\dhome/cmds2.m}

\begin{figure}[ht] 
\centering
\btab{l}{
\btab{lll}
{{\sm (a)}&{\sm (b)}&{\sm (c)}\\
\ig[height=0.25\tew]{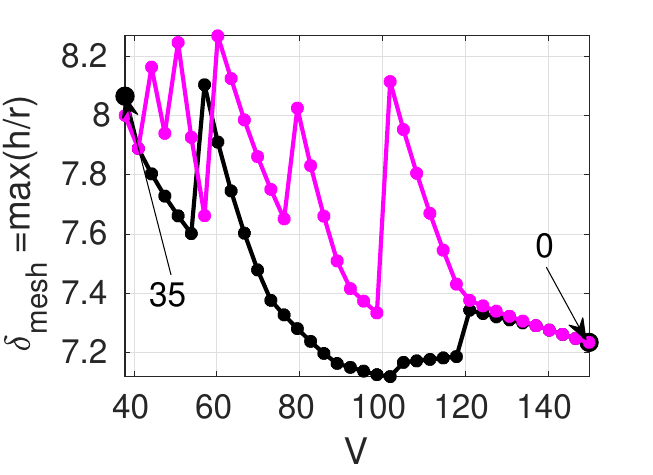}&\ig[height=0.25\tew]{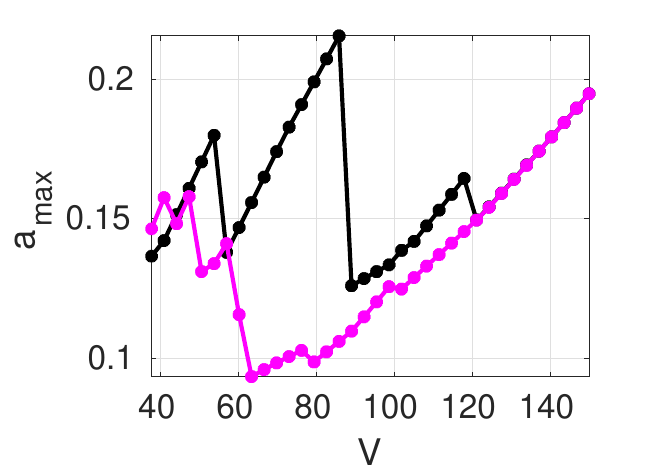}&
\ig[height=0.25\tew]{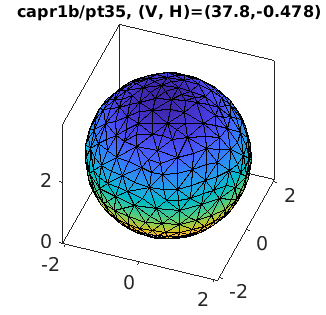}
}\\
\btab{l}{{\sm (d)}\\[-7mm]
\hs{0mm}\ig[height=0.26\tew]{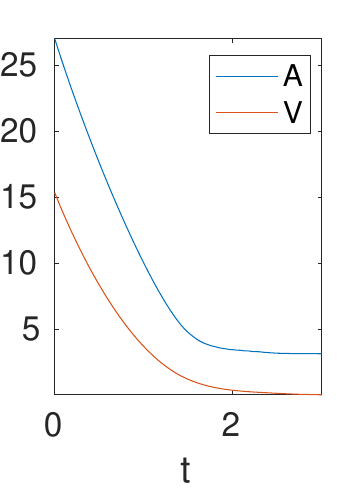}
\hs{0mm}\ig[height=0.26\tew]{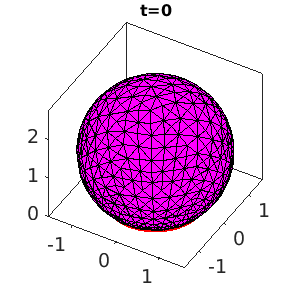}
\hs{-2mm}\rb{2mm}{\ig[height=0.25\tew]{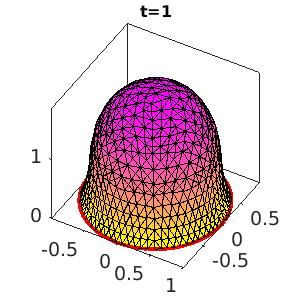}}
\hs{3mm}\rb{26mm}{\btab{l}{\ig[height=0.22\tew]{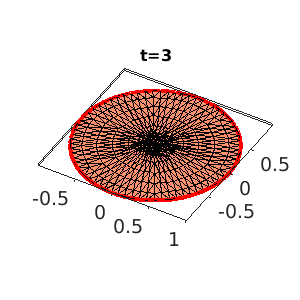}\\[-5mm]
\ig[height=0.13\tew]{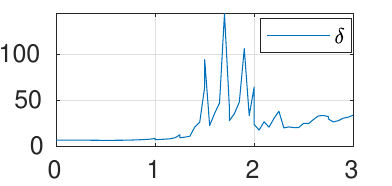}
}}}
}
\vs{-2mm}
\caption{{\small Results from {\tt spcap1/cmds3.m}. (a)-(c) 
continuation backwards in $V$ from $V{\approx}150$ ($n_p{=}1452$); coarsening 
each 5th step ({\tt capr1b}, black, $n_p{=}644$ at $V{=}40$) 
vs coarsening when $\delm>8$ (magenta, $n_p{=}650$ at $V{=}40$)). 
(d) MCF from the spherical cap at $V{\approx}15$.  
time series of $A$ and $V$, sample plots, and time series of $\delm$ 
(last plot). Coarsening at times $t=0.25 j$, 
altogether from $n_p=773$ at $t=0$ to $n_p=450$  at $t=3$.  
\label{spcf1b}}}
\end{figure}

In {\tt cmds3.m} and Fig.\,\ref{spcf1b} we {\em decrease} $V$ from $V\approx 150$ 
(running the branch {\tt capr1} from Fig.\,\ref{spcf1c} backwards), 
and test the MCF from a spherical 
cap at $V\approx 15$. For both, 
because the shrinking of the caps gives mesh distortions, 
the main issue is 
that we now need to alternate continuation/flow and mesh--{\em coarsening}. 
For the continuation we give two options: similar to the refinement 
for increasing $V$ in Fig.\,\ref{spcf1c}, we either coarsen after a fixed 
number of steps (black branch), or when $\delm>8$ (magenta branch). 
Both here work efficiently only 
until $V\approx 35$, after which new parameters for the coarsening 
should be chosen. 
 For the MCF in (d) we similarly coarsen after a 
given number of time steps. With this we can flow back to the disk, 
more or less reached at $t=3$, but the last plot in (d) shows 
that along the way we have strongly distorted meshes, which are 
somewhat repaired in the coarsening steps, and the final distortion 
with $\delm\approx 30$ is not small but OK. 

\hulst{linerange=1-5,firstnumber=1}{\dhome/cmds3.m} 
\hulst{caption={{\small Selection from {\tt spcap1/cmds3.m}; 
 decreasing 
$V$ by going backwards, and MCF; both need to be combined with coarsening. 
Omission between lines 5 and 14 deal with plotting, and further experiments 
are at the end of {\tt cmds3.m}.}},
linerange=14-22,label=spcl4,firstnumber=14}{\dhome/cmds3.m}

\brem\label{mcfrem1}{\rm The performance of the 
MCF as in Fig.\,\ref{spcf1b}, based on our simple explicit 
Euler stepping, depends on the choice of flow parameters, i.e., 
step size {\tt dt}, number  $n_{\text{f}}$ of steps before coarsening, 
and coarsening factor $\sig$. With too weak coarsening (large 
$n_{\text{f}}$, or small $\sig$), triangles may degenerate. 
%(``move through'' each other). 
Too aggressive coarsening (large $\sig$) 
may lead to wrong identification of boundary edges. Altogether, at this 
point we must recommend trial and error. 
We also tested the use of {\tt degcoarsenX} instead of 
{\tt coarsenX} (see {\tt cmds3.m}) for backward continuation and for MCF, 
but this here does not give 
good results; see Fig.\ref{nodf2b} for a successful use of {\tt degcoarsenX} 
in a related problem. %, where it works well. 
}\eex\erem

\subsection{Some minimal surfaces}
\label{obcsec}
\def\dhome{geomtut/bdcurve} 
Plateau's problem consists in finding soap films $X$ spanning a 
(Jordan) curve (a wire) $\ga$ in $\R^3$, and minimizing area $A$. 
Mathematically, we seek a {\em minimal} surface $X$, i.e., $H(X)\equiv 0$, 
with $\pa X=\ga$. 
Such problems have a long history, and already Plateau discussed 
non--uniqueness and bifurcation issues, called ``limits of stability'' 
in \cite{Pl73}. 

A classical example for which a bifurcation is known is Enneper's surface, 
see \S\ref{ensec}. However, in the demo {\tt bdcurve} we first start 
with other BCs, meant to illustrate options (and failures) for prescribing 
boundary values in our numerical setup. We introduce 
parameters $\al\in\R$ and 
$k\in\N$ (angular wave number)  and a switch {\tt p.bcsw}, and 
consider BCs of the form 
\hual{\label{bc1b}
&u|_{\pa X}=0, \quad
\text{(for {\tt p.bcsw=0})},\\
&\label{bc2} 
X_3|_{\pa X}=\al\sin(k\phi), \quad
\text{$\phi=\arctan(y/x)$\quad (for {\tt p.bcsw=1})},\\
&\pa X=\ga(\cdot;\al,k)\quad \text{(for {\tt p.bcsw=2})},\label{bc3}
}
where $\ga$ in \reff{bc3} is a prescribed boundary curve, depending 
on parameters $\al,k$.  Specifically, in \S\ref{bc2sec} we choose 
\huga{\label{gadef}
\ga(\phi;\al,k)=\bpm \beta\cos(\phi)\\\beta\sin\phi\\
\al\cos(k\phi)\epm, \quad \phi\in[0,2\pi), \quad 
\beta=\sqrt{1-\al^2\cos^2(k\phi)}. 
}

For \reff{bc2}, $\pa X$ is not uniquely determined by the parameter $\al$ 
(and fixed $k$), and this illustrates how our scheme 
\reff{n1} can fail, and that the condition $Y\in \CN_C$ 
in general cannot be dropped in Lemma \ref{l-normvar}. Relatedly, 
for \reff{bc2} the continuation can genuinely depend on the 
continuation stepsize {\tt ds}, 
as different predictors give different BCs \reff{bc2}. In other words, 
the problem is under--determined and 
the continuation algorithm itself ``chooses'' the BCs. 
Thus, \reff{bc2} is a cautionary example, though it produces interesting 
minimal surfaces. 
On the other hand, \reff{bc3} is a genuine DBC with unique continuation, 
which however requires a modification of the 
``standard'' {\tt updX.m}, and careful mesh handling. 
Together, \reff{bc2} and \reff{bc3} are meant to illustrate 
options. 

%\brem\label{bdrem1}{\rm  
The condition on $\beta$ in \reff{gadef} 
yields that 
$\|\ga\|_2=1$, i.e., that $\ga$ lies on the unit sphere, for $\al\in[0,1]$. 
Moreover, the projection of $\ga$ into 
the $x$--$y$ plane is injective, and this is useful since we then can 
extract $\phi$ from $\pa X$. In \S\ref{ensec} we treat a 
variant of \reff{gadef}, associated to the Enneper surface, 
where the discretization of $\ga$ requires a further trick. 
On the other hand, $\ga$ from \reff{gadef}  becomes singular 
at $\phi=j\pi/k$ as $\al\to 1$, which is useful to test mesh--handling. 
Thus, \S\ref{bc2sec} and \S\ref{ensec} are quite related, but 
illustrate different effects. 
%}\eex\erem 

Table \ref{bdctab} shows the used files, 
Listing \ref{bdcl1} shows {\tt sGbdcurve}, and 
Listing \ref{bdcl2} the ``new'' (compared to {\tt spcap1/}) files needed 
to run the BCs \reff{bc3}. The other files are essentially as in 
{\tt spcap1/}, except that we now have altogether five parameters 
$(H,V,A,\al,k)$, and that we use the additional parameter {\tt p.bcsw}. 

\taskip
\begin{table}[ht]\caption{
Files in {\tt pde2path/demos/geomtut/bdcurve}. 
\label{bdctab}}
\centering\vs{-2mm}
{\small 
\begin{tabular}{l|l}%p{0.15\tew}|p{0.8\tew}} 
{\tt cmds1a.m}&continuation in $\al$ (and $H$) for \reff{bc2}, 
see Fig.\ref{spcf3}; MCF tests in {\tt cmds1b}. \\
{\tt cmds2.m}&continuation in $\al$ for \reff{bc3}, see Fig.\,\ref{bdcf1}.\\
{\tt bdcurveinit.m}&Initialization, very similar to {\tt scinit}.\\
{\tt sGbdcurve.m}&very similar so {\tt sGsc}, except for the BCs.\\
{\tt updX.m}&mod of standard {\tt updX};  
for {\tt p.bcsw=2} setting the boundary curve. \\
{\tt bcX.m}&user function to give $\ga$, here implementing \reff{gadef}. 
\end{tabular}
}
%\ece\vs{-4mm}
\end{table}
\teskip

\hulst{caption={{\small {\tt bdcurve/sGbdcurve.m}, with BCs 
depending on {\tt p.bcsw}.}},label=bdcl1}{\dhome/sGbdcurve.m}

\hulst{firstnumber=1}{\dhome/bcX.m}
\hulst{caption={{\small {\tt bcX.m} and 
{\tt updX.m} from {\tt bdcurve/}}, needed to run with the 
BCs \reff{bc3}.}, label=bdcl2, firstnumber=1}{\dhome/updX.m}

\subsubsection{Prescribing one component of $X$ at the boundary} 
\label{bc1sec}
\begin{figure}[H]
\centering 
\btab{l}{{\sm (a)}\hs{45mm}{\sm (b)}\\[-1mm]
\hs{0mm}\ig[height=0.26\tew]{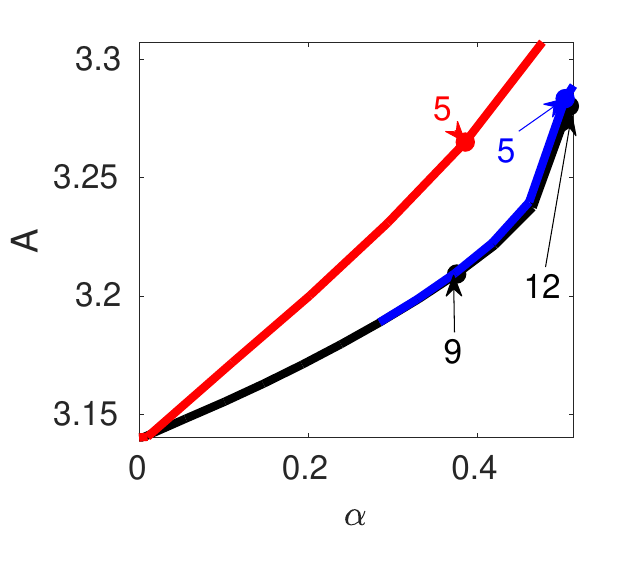}
\hs{0mm}\rb{2mm}{\ig[height=0.25\tew]{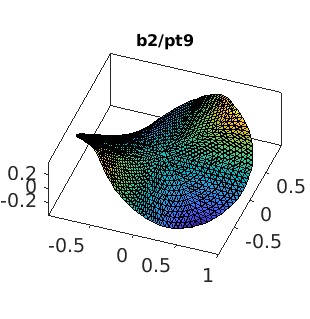}}
\hs{-1mm}\rb{6mm}{\ig[height=0.23\tew]{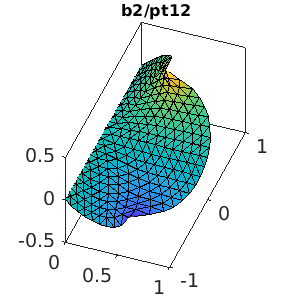}}
\hs{-0mm}\rb{6mm}{\ig[height=0.22\tew]{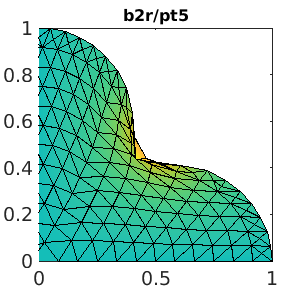}}\\[-3mm]
{\sm (c)}\hs{40mm}{\sm (d)}\\[-1mm]
\hs{-1mm}\ig[height=0.25\tew]{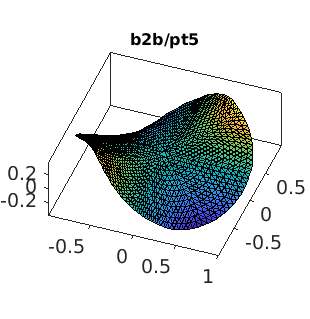}
\hs{0mm}\ig[height=0.23\tew]{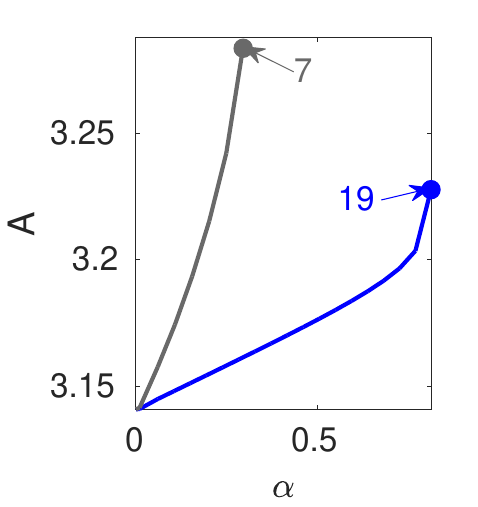}\hs{0mm}
\hs{0mm}\ig[height=0.2\tew]{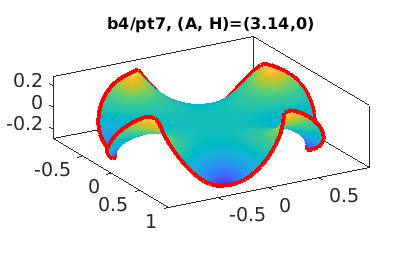}
\hs{-3mm}\ig[height=0.22\tew]{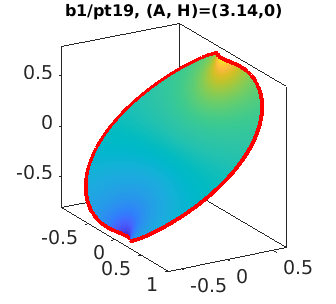}\\[-3mm]
{\sm (e)}\hs{40mm}{\sm (f)}\\[-1mm]
\hs{0mm}\ig[height=0.24\tew]{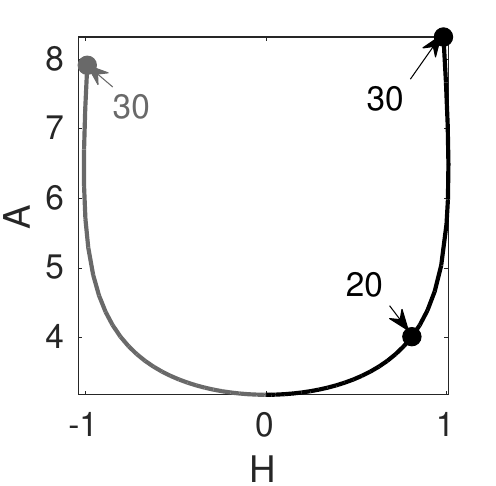}\hs{3mm}\ig[height=0.26\tew]{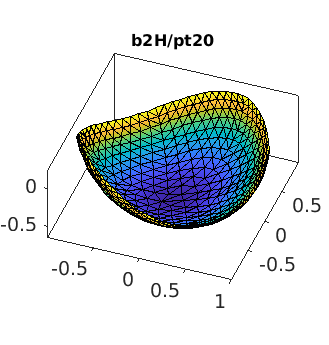}
\hs{-0mm}\ig[height=0.26\tew]{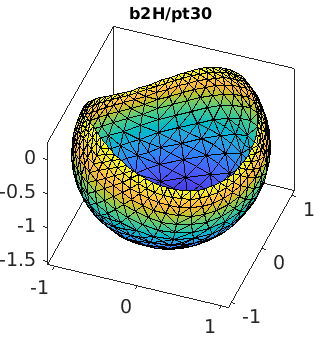}
\hs{-0mm}\ig[height=0.26\tew]{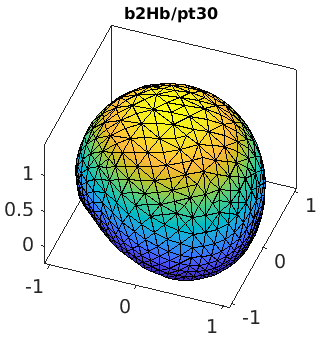}
}
\vs{-4mm}
\caption{{\small Results from {\tt bdcurve/cmds1a.m}. 
(a) Continuation in $\al$ with BCs \reff{bc1}, $k=2$:  
black and red branches with $n_p{=}945$ mesh points and 
fixed {\tt ds=0.05} (black) and {\tt ds=0.1} (red), 
illustrating the step--length dependence of the continuation;   
blue branch {\tt b2r} via refinement at boundary.  
Samples in (b) and (c), partly with cropping. %(not showing the full $X$), 
%and {\tt b2/pt12} and {\tt b2r/pt5} 
%shortly before continuation failure due to ``vertical $X$'' at the boundary. 
(d) Like (a,b) but on finer meshes, with $k{=}1$ and $k{=}4$, and with 
marking $\pa X$ in red. % (by setting {\tt p2pglob.showbd=2}). 
All of (a--d) are {\em minimal} surfaces, 
i.e., $H\equiv 0$. 
(e,f) Switching back to continuation in 
$H$ at {\tt b2/pt6}.  \label{spcf3}}}
\end{figure}

In {\tt cmds1a.m} (Listing \ref{bdcl2b}) we continue \reff{bc2} in $\al$, 
starting 
with $\al=0$ at the flat disk, and first with angular wave number $k=2$. 
Some results are shown in Fig.\,\ref{spcf3}.

\hulst{caption={{\small Start of {\tt bdcurve/cmds1a.m}, running 
BCs \reff{bc2}.}},label=bdcl2b,linerange=1-4}{\dhome/cmds1a.m} 

As we increase $\al$, the surface lifts at $\phi=\pi/4$ and $\pi=5\pi/4$ 
according to $X_3=\al\sin(2\phi)$, and sinks at $\phi=3\pi/4$ and $\phi=7\pi/4$. 
Near $\al=0.5$ ({\tt b2/pt12}), $X$ becomes vertical at these angles, 
and hence our scheme \reff{n1} can no longer continue to fulfill the BCs. 
To better resolve  the boundary we use some mesh--refinement 
{\em only at the boundary}. 
For this we choose {\tt p.fuha.e2rs=@e2rsbdry} at {\tt b2/pt6} and 
obtain the blue branch (with a sample top view as last plot in (b)), 
which however naturally runs into the same 
continuation failure at $\al\approx 0.5$. Although this was 
on quite coarse meshes, none of this changes on finer meshes, 
and hence this mainly serves as an example of necessary failure of 
the algorithm \reff{n1}, and as an example of mesh refinement with 
{\tt e2rsbdry}.

The red branch in (a) together with sample (c) 
shows that here the branches are continuation stepsize {\tt ds} dependent. 
In (d) we choose finer meshes and wave numbers $k=1$ (blue branch {\tt b1}) and 
$k=4$ ({\tt b4}, grey), and get analogous results up to continuation 
failure. In (e,f) we switch back 
to continuation in $(H,A)$ from {\tt b1/pt6}, in both directions 
of positive (black branch) and negative (grey branch) $H$.  
As $\al$ is now fixed again, $\pa X$ stays fixed even with the BCs \reff{bc2}. 
The branches are {\tt ds}--independent again, 
and $H$ asymptotes to nonzero $\pm H_\infty$ as $A\to\infty$. 
In {\tt cmds1b.m} we run MCF (not shown) from selected solutions from (f), 
where again we need to  undo the refinement which happened, e.g., 
during the continuation $H$ from {\tt b2Hb/pt0} to {\tt pt30}, 
and  thus we alternate between {\tt geomflow} and {\tt coarsenX} as 
in {\tt spcap1/cmds2.m} and Fig.\ref{spcf1b}.

\subsubsection{A Plateau problem}\label{bc2sec}
In {\tt cmds2.m} we choose the BCs \reff{bc3} with $\ga$ from \reff{gadef}. 
We again continue in $\al$, 
for $k=2,3$, starting at $\al=0$ with the unit disk. 
The basic idea (Listing \ref{bdcl2}) 
to implement \reff{bc3}, \reff{gadef} is to 
\huga{
\text{set $\pa X=\ga$ in {\tt updX}, and $u|_{\pa X}=0$ in {\tt sGbdcurv}.}
}

\hulst{caption={{\small First 15 lines from 
{\tt bdcurve/cmds2.m}}, using the BCs \reff{bc3}.},label=bdcl4, 
linerange=1-15}{\dhome/cmds2.m}

\begin{figure}[ht] 
\bce 
\btab{llll}{{\sm (a)}&{\sm (b)}&{\sm (c)}&{\sm (d)}\\ 
\hs{-2mm}\ig[height=0.27\tew]{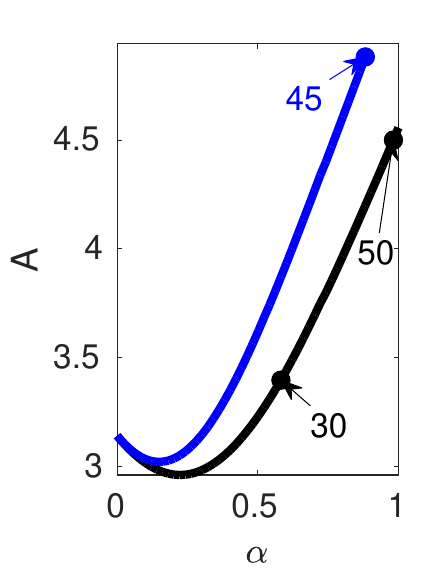}&
\hs{0mm}\ig[height=0.28\tew]{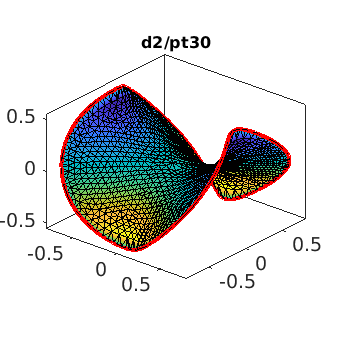}&
\hs{-4mm}\ig[height=0.28\tew]{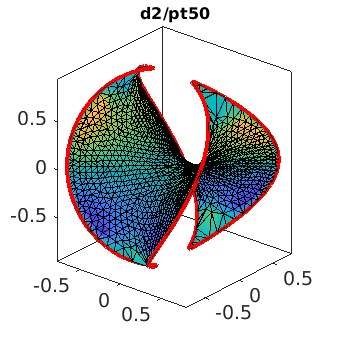}&
\hs{-4mm}\ig[height=0.28\tew]{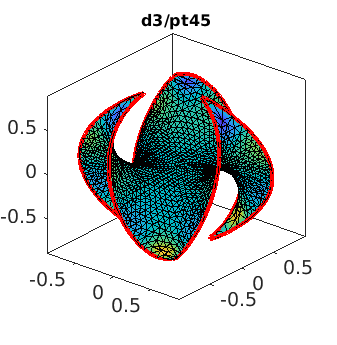}
%\hs{-4mm}\ig[height=0.24\tew]{p3/al2r-7}
}
\ece 

\vs{-5mm}
\caption{{\small Results from {\tt cmds2.m} for BCs \reff{bc3} 
with $k{=}2$ (black branch {\tt d2}) and $k{=}3$ (blue branch {\tt d3}). 
Along the way in, e.g., {\tt d2} we do 3 refinements and coarsenings, 
and the total number of mesh points only increases mildly from $n_p{=}945$ 
to $n_p{=}1177$. These are really two different ($k{=}1$ vs $k{=}2$) continuation 
problems, out of inifinitely many ($k{\in}\N$), and hence the two 
branches in (a) are from different problems and are both stable. 
 \label{bdcf1}}}
\end{figure}

Figure \ref{bdcf1} shows some results from {\tt cmds2.m}. 
The crucial points are that as we increase 
$\al$ (in particular beyond $\al=0.2$, say) we 
\bci 
\item move mesh points via {\tt moveX(p,dt,it)} after {\tt ncs} continuation steps (here {\tt ncs=1});   
\item after {\tt nis=15} ``inner'' steps refine $X$ (introduce new points), here near the boundary, {\em and} coarsen $X$, here via {\tt degcoarsenX(p,sigc)}, 
to remove ``bad'' triangles.
\eci 
The parameter {\tt dt} in {\tt moveX} is the Euler step size to balance 
the ``truss forces'' \cite{PeSt04} 
({\tt nit} gives the number of iterations), while {\tt sigc} in 
{\tt degcoarsenX} has a similar meaning as in {\tt refineX} and 
{\tt coarsenX}, i.e., giving the fraction of triangles to coarsen.%
\footnote{\label{dcfoot}In more detail, {\tt degcoarsenX} can also be called as 
{\tt p=degcoarsenX(p,sigc,iter)}, where {\tt iter} (default 5) gives 
the number of internal iterations, or as 
{\tt p=degcoarsenX(p,sigc,iter,keepbd)} where {\tt keepbd=1} (default 0) 
means that boundary triangles are kept, which is mainly needed for 
periodic BCs, see \S\ref{nodnumpBC}.}
Again, the parameters {\tt ncs, nis, sigr} and {\tt sigc} 
are generally highly problem 
dependent and it may require (educated) trial and error to find 
good values. 
In summary, Fig.\,\ref{bdcf1} shows that with 
a good combination of 
{\tt moveX, refineX} and {\tt degcoarsenX} we can continue 
rather complicated minimal surfaces $X$ 
($X$ with complicated boundary curve $\ga$) 
with reasonable meshes.%
\footnote{As already said in Rem.~\ref{Xcontrem2}b), 
due to the analogy between the truss forces and surface tension 
(constant in minimal surfaces) {\tt moveX} 
works particularly well for minimal $X$.}  

\subsubsection{Bifurcation from the Enneper surface}\label{ensec}
\begin{figure}[ht] 
\centering
\btab{lll}{{\sm (a)}&{\sm (b)}&{\sm (c)}\\[0mm]
\hs{-2mm}\ig[height=0.28\tew]{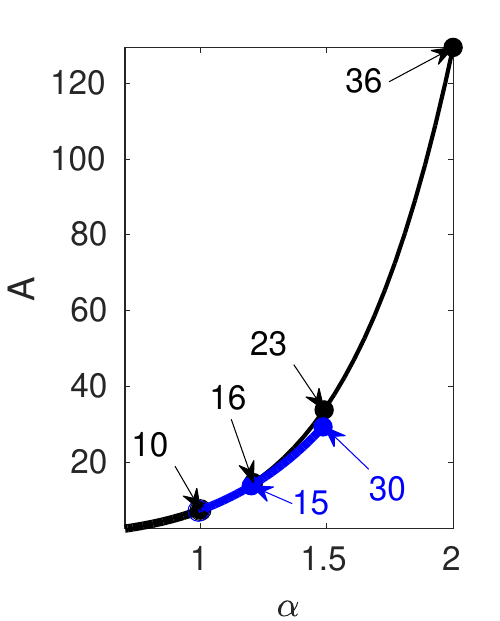}
&\hs{-4mm}\ig[height=0.28\tew]{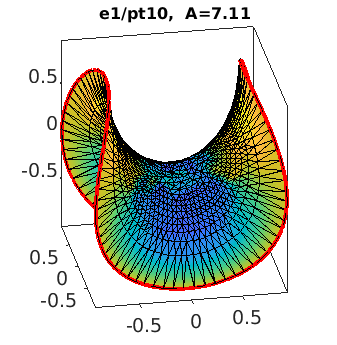}
&\hs{-3mm}\ig[height=0.28\tew]{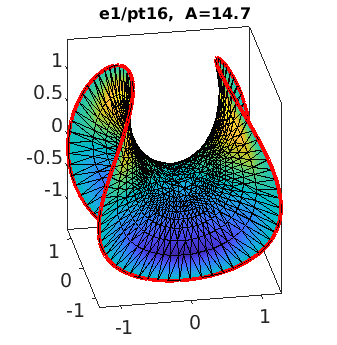}
\hs{-3mm}\ig[height=0.28\tew]{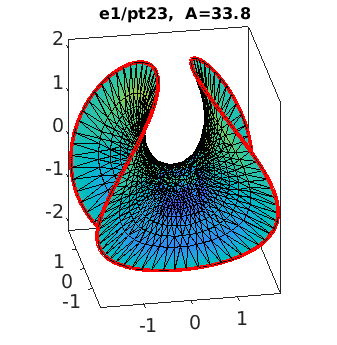}\\[-2mm]
{\sm (d)}&{\sm (e)}&{\sm (f)}\\ [-0mm]
\hs{-2mm}\ig[height=0.28\tew]{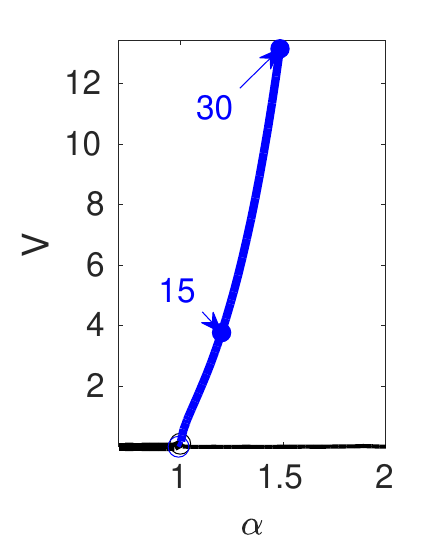}
&\hs{-3mm}\ig[height=0.28\tew]{p5/e1-36}
&\hs{-3mm}\ig[height=0.28\tew]{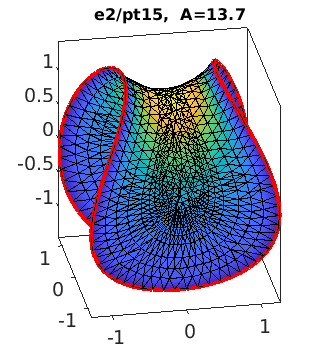}
\hs{-6mm}\ig[height=0.28\tew]{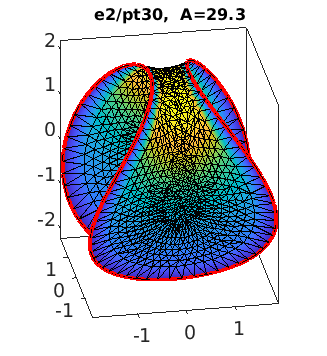}\\[-3mm]
{\sm (g)}&{\sm (h)}\\[0mm]
\hs{-2mm}\ig[height=0.24\tew,width=0.23\tew]{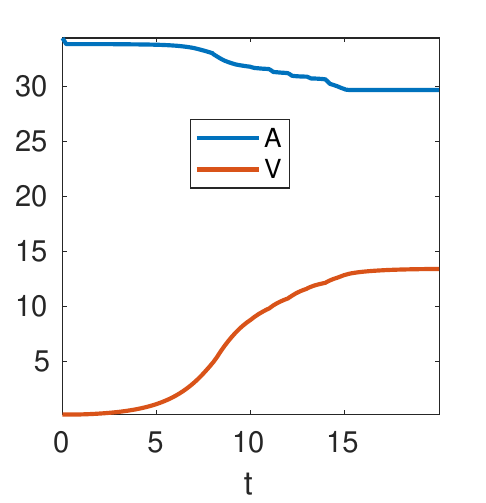}
&\hs{-1mm}\ig[height=0.25\tew]{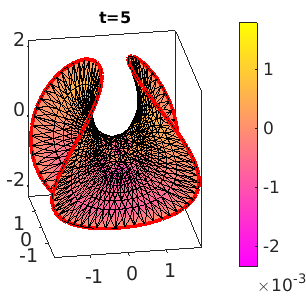}
&\hs{-1mm}\ig[height=0.25\tew]{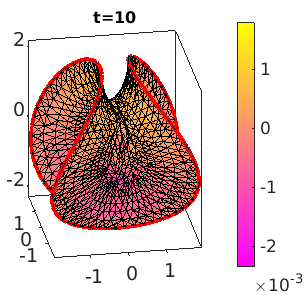}
\hs{-1mm}\ig[height=0.25\tew]{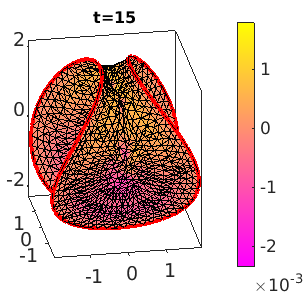}
}
\vs{-3mm}
\caption{{\small Bifurcation from the Enneper surface $X_E$, 
$A$ over $\al$ (a), and $V$ over $\al$ (d). At $\al=1$ ({\tt e1/pt10} in (b)), 
the branch {\tt e1b} (blue) with smaller $A$ bifurcates from {\tt e1} (black), 
samples in (b,c) and (e,f). (g,h) MCF from perturbation of {\tt e1/pt23} 
to {\tt e2/pt30}, samples showing $H$. %\huc{producing nice mesh, subseq.cont also nice.}
 \label{ef1}}}
\end{figure}

\def\dhome{geomtut/enneper} 
The Enneper surface is a classical minimal surface. 
Bounded parts of it can be parameterized by%
\footnote{see also Remark \ref{Prem2} for the Enneper--Weierstrass 
representation}  
\huga{\label{en1}
X_E=X_E(r,\vt)=\bpm r\cos(\vt)-\frac{r^3}3\cos(3\vt)\\
-r\sin(\vt)-\frac{r^3}3\sin(3\vt)\\
r^2\cos(2\vt)\epm,\quad (r,\vt)\in D_\al=[0,\al)\times [0,2\pi), 
}
see Fig.\ref{ef1}. 
We start by reviewing some basic facts, see \cite{BT84} and the references 
therein. 
For $\al\le 1/\sqrt{3}$, the boundary curve 
\huga{\label{egd} 
\ga(\vt;\al)=\bpm \al\cos(\vt)-\frac{\al^3}3\cos(3\vt), -\al\sin(\vt)-\frac{\al^3}3\sin(3\vt), \al^2\cos(2\vt)\epm,\quad \vt\in[0,2\pi)
} 
has a convex projection to the $x$--$y$--plane, and for 
$1/\sqrt{3}< \al\le 1$ the projection  is still injective. 
This yields uniqueness (of the minimal surface spanning $\ga$) for 
$0<\al\le 1$ (see \cite{ruch81} for $\al\in(1/\sqrt{3},1]$). 
For $\al>1$ uniqueness of $X_E$ fails, i.e., at 
$\al=1$ we have a (pitchfork, by symmetry) bifurcation of different minimal surfaces spanning 
$\ga_\al$  \cite{Nit76}. This has been analyzed in detail in \cite{BT84} 
as a two--parameter bifurcation problem, showing a so called 
cusp catastrophe.%
\footnote{See, e.g., \cite[Example 1.30]{p2pbook} and the 
references therein for comments on cusps (and other catastrophes).}

In the demo {\tt enneper} 
we simply choose $\al$ as a continuation/bifurcation parameter for 
\huga{\label{ennpde}
H(X)=0,\quad \pa X=\ga_\al, 
} 
and get the pitchfork bifurcation at $\al=1$. The 
used files {\tt bcX.m, cmds1.m, cmds2.m, enninit.m, sGenn.m, updX.m} are 
very similar to those from the demo {\tt bdcurve}, 
but we also include a Jacobian {\tt sGjacenn.m}, a function {\tt mytitle.m} 
for customized titles, and {\tt thinterpol.m}, discussed next. 

The problem \reff{ennpde} is ``easy'' in the sense that we have the explicit 
parametrization \reff{en1} which we can use at any $\al$, 
but like in \S\ref{bc2sec} it does require care with the meshing, 
and compared to \reff{gadef} it requires an additional trick 
to update $\vt$ on $\pa X$ after mesh adaption (at the boundary): 
Since we cannot in general extract $\vt$ from \reff{egd} 
from the projection to the $x$--$y$ plane (which is not 
injective for $\al>1$), we keep a field $\vt={\tt p.th}$ associated 
to {\tt p.idx} (the indices of $\pa X$) in the given discretization. 
Then, if {\tt p.X1} is obtained from refining {\tt p.X} 
with new mesh--points {\tt p.nX} on $\pa X$, then we need to update 
the $\vt={\tt p.th}$ values of {\tt p.nX}. This is done in 
{\tt p=thinterpol(p,idxold,thold)} by 
\bci
\item finding the (old--point) neighbors of the new points on $\pa X$; 
\item linear interpolation of the neighbors' $\vt$ values to the new points. 
\eci 
This is a question of indexing, and we refer to the source of {\tt thinterpol.m} for comments. A refinement step thus takes the 
form
\bce 
{\tt idold=p.idx;\ \ thold=p.th;\ \ p=refineX(p,sigr);\ \ p=thinterpol(p,idold,thold); }
\ece 
see {\tt cmds1.m} which produces Fig.\,\ref{ef1}. 
The other files in {\tt enneper/} are very much like in {\tt bdcurve/}. 

At $\al=1$ we find a supercritical pitchfork bifurcation from $X_E$, 
branch {\tt e1} (black), to a branch {\tt e2} (blue) 
which breaks 
the $(x,y,z)\mapsto (-y,x,-z)$ symmetry of $X_E$ 
(rotation by $\pi/2$ around the $z$ axis and mirroring at the $z=0$ plane). 
The solutions ``move up'' (or down) in the middle, 
which decreases $A$ compared to $X_E$, cf.~(c) vs (f). (d) illustrates that 
the (algebraic) volume $V$ of $X_E$ is always zero. The numerical continuation 
of {\tt e1} to large $\al$ is no problem, using suitable mesh--adaption, 
even as $\ga(\cdot;\al)$ self--intersects for $\al>\sqrt{3}$, because 
the associated parts of $X_E$ do not ``see'' each other, cf.~(e) for an 
example. The continuation of {\tt e2} 
to larger $\al$ is more difficult, and fails for $\al>1.5$, 
as for instance shortly after {\tt e1b/pt30} 
we can no longer automatically adapt the mesh near the top. 

However, physically the change of 
stability at the symmetry breaking pitchfork at $\al=1$ is most 
interesting. Using suitable combinations of {\tt geomflow}, 
{\tt refineX}, {\tt degcoarsenX} and {\tt moveX} we can use MCF to 
converge for $\al>1$ and $t\to\infty$ to {\tt e2}, from a 
variety of ICs, for instance from perturbations of {\tt e1}, 
see Fig.\,\ref{ef1}(g,h), and {\tt enneperflow.avi} in \cite{sig}. 
After convergence we can then 
again continue the steady state, see {\tt cmds2.m}.

\subsection{Liquid bridges and nodoids}
\label{lbsec}%\label{nodsec}
\def\dhome{geomtut/libri} 
Weightless liquid bridges are CMC surfaces with prescribed boundary 
usually consisting of two parallel circles wlog centered 
on the $z$-axis at a fixed distance $l$ and parallel to the 
$x$--$y$ plane. Additionally there is a volume constraint, which makes the 
problem different from Plateau's problem.  See for instance \cite{SAR97} 
and the references therein for physics background and results (experimental, 
numerical, and semi-analytical). 
%They for instance model interfaces between two liquids, neglecting gravity. 

We consider liquid bridges between two fixed circles $C_1$ and $C_2$ of 
\huga{\label{libribc}
\text{radius $r=r^*=1$, parallel to the $x$--$y$ axis and 
centered at $z=\pm l=\pm 1/2$. } 
}
A trivial solution $X_0$ is the cylinder, with $H=1/2$, volume 
$V=2\pi l$ and area $A=4\pi rl$ (without the top and bottom disks). 
Further explicit solutions are known in 
the class of surfaces of revolution, for instance nodoids.
We first review some theory for nodoids with DBCs, 
and then continue  basic liquid 
bridges (embedded nodoids), with 
bifurcations to non axial branches, see 
Figures \ref{nodf0} and \ref{nodf00}. In Figure \ref{nodf2}
we then start directly with nodoids with one ``inner loop''. 
Nodoids with ``periodic'' BCs are studied in \cite{MaPa02}, 
and numerically in \S\ref{nodnumpBC}, where we also comment on the theory for these. 

\subsubsection{Nodoid theory} \label{nodsec}
In \cite{KPP17}, a family of nodoids $\mathcal N(r,R)$ is 
parameterized by the neck (smallest) radius $r$ and the buckle 
(largest) radius $R$. 
Let $l\in\R$ and $C_1,C_2\subset\R^3$ be two circles of radius $r^*$ centered 
at heights $z=\pm l$ and parallel to the $x$--$y$ plane. With the
two parameters $a,H\in\R$ the nodoids are parameterized by the
nodary curve 
\hual{\label{nodpar} (x,z):[-t_0,t_0]&\rightarrow\R^2, \quad 
  t\mapsto\bpm x(t),z(t)\epm=\bpm \frac{\cos t+\sqrt{\cos^2 t+a}}{2
    \abs{H}},\ \ \frac{1}{2\abs{H}} \int_0^t \frac{\cos \tau+\sqrt{\cos^2 \tau
      +a}}{\sqrt{\cos^2 \tau+a}}\cos \tau \dd \tau\epm, 
} 
which is then rotated around the $z$ axis, i.e., 
\hual{\label{paranod}\CN_{t_0}:
  M&\rightarrow \R^3,\qquad (t,\theta)\mapsto\bpm x(t)\cos
  \theta,x(t)\sin \theta,z(t)\epm,
} 
where $M=[-t_0,t_0]\times[0,2\pi)$. 
Thus, in terms of \S\ref{geom} 
these nodoids are immersions of cylinders. While \reff{nodpar} 
only gives nodoids with an even number of self
intersections (or none), shifting $t_0$ also gives odd numbers of self 
intersections. From the immersion $\CN_{t_0}$,
we can determine geometric quantities by evaluating the
parametrization at the endpoints. For example the height and the
radius are given by 
\hueq{\label{lenrad} 2l=\frac{1}{
    \abs{H}}\int_0^{t_0} \frac{\cos t+\sqrt{\cos^2 t +a}}{\sqrt{\cos^2
      t+a}}\cos t \dd t,\quad 
r^*=\frac{\cos t_0+\sqrt{\cos^2
      t_0+a}}{2 \abs{H}}, 
} 
and the buckle radius (at $t=0$) is 
$\ds R=\frac{1+\sqrt{1+a}}{2\abs{H}}$. 
Implicitly, the equations in \eqref{lenrad} define $a(t_0)$, hence 
also the mean curvature $H$, and thus $t_0$ parameterizes a family of 
nodoids $t_0\mapsto \CN_{t_0}$.  Conversely, given $r,l$ in \reff{libribc}, 
the implicit equation 
\hueq{\label{impnod}
{l\over 2r}\left(\cos t_0+\sqrt{\cos^2t_0+a}\right)-
\left(\sin t_0+\int_0^{t_0}{\cos^2\tau\over \sqrt{\cos^2\tau +a}}
\dd \tau \right)=0 
}
defines all possible combinations of $a$ and $t_0$ satisfying the 
boundary condition, which we exploit to relate our numerics 
to results from \cite{KPP17}, see Remark \ref{koirem}. 

In order to detect bifurcations from the family \reff{paranod}, we 
search for
Jacobi fields vanishing on the boundary, cf.~\reff{jac}. 
The unit normal vector (field) of $\CN_{t_0}$ is 
\hugast{
N=\bpm \cos t\cos\theta,\cos t \sin\theta,\sin t\epm, 
\quad t\in[-t_0,t_0),\ \vt\in[0,2\pi), 
}
and for every fixed vector $\vx\in\R^3$, the function 
$\spr{\vx,N}$ is a solution to \eqref{f-mean}. 
So the task is to find  $\vx$ and $t_0$ such that the 
Dirichlet BCs are fulfilled. The components of 
$N$ have zeros if the nodoid meets the boundary 
horizontally (parallel to the $x$--$y$ plane), which happens 
at $t_0=\frac \pi 2+n\pi$, or vertically, which happens at 
$t_0=n\pi$ for $n\in\N$. 
Choosing the unit basis $(e_i)_{i=1,2,3}$, we have in the horizontal 
case that $\spr{e_i,N}\lvert_{\pa \CN_{t_0}}=0$ for $i=1,2$, 
and in the vertical case $\spr{e_3,N}\lvert_{\pa \CN_{t_0}}=0$. 

\blem\label{l-deins} \cite[Lemma 3.4 and Proposition 3.6]{KPP17} 
Consider the one parameter family $\CN_{t_0}$. If for
some $t_0\in\R_+$ the normal vector at $\pa \CN_{t_0}$ is
\begin{enumerate}
\item $N=\bpm 0,0,\nu(x)\epm$, then $L=\pa_u H(u)$ has
  a double zero eigenvalue. 
\item $N=\bpm\nu_1(x),\nu_2(x),0\epm$ then $L=\pa_u H(u)$ has a
  simple zero eigenvalue. 
\end{enumerate}
The immersions are isolated degenerate, i.e., there exists
an $\eps>0$ such that $(\CN_{t})_{t\in[t_0-\eps,t_0+\eps]}$ has a jump
in the Morse index.  In $\textit 1.$ this occurs for
$t_0=\frac \pi 2 +k\pi$, and in $\textit 2.$ for $t_0=k\pi$, for
every $k\in\N$.  
\elem 

Now general bifurcation results (see the discussion after Lemma \ref{derlem}) 
yield the existence of bifurcation points at the horizontal and
vertical cases presented in Lemma \ref{l-deins}.
\bthm \label{t-nodbif} \cite[Propositions 3.5 and 3.6]{KPP17} 
In cases 1.~and 2.~in Lemma \ref{l-deins} we have bifurcation points 
for the continuation in $H$. Moreover, 
\begin{enumerate}
\item if $\psi=\spr{e_i,N}\in\ker L$ for $i=1,2$, then the
  bifurcating branch breaks the axial symmetry; 
\item if $\psi=\spr{e_3,N}\in\ker L$, then the bifurcating 
branch breaks the $z\mapsto -z$ symmetry. 
\end{enumerate}\ethm

\subsubsection{Nodoid continuation with fixed boundaries}\label{nodnumDBC}
\def\dhome{geomtut/nodDBC} 
Nodoids with DBCs at the (fixed) top and bottom circles are treated in 
the demo {\tt nodDBC}. Table \ref{nodtab} lists the pertinent files. 
We treat two cases: 
\bci 
\item Short embedded nodoids (liquid bridges)  in {\tt cmds1.m}, 
starting from the cylinder (eventually continued to self-intersecting nodoids). 
\item Long nodoids (with self--intersections from the start) in {\tt cmds2.m}. 
\eci 

\taskip
\begin{table}[ht]\caption{Files in {\tt pde2path/demos/geomtut/nodDBC}; 
\label{nodtab}} 
\centering\vs{-2mm} {\small
\begin{tabular}{p{0.25\tew}|p{0.76\tew}}
{\tt cmds1.m}&continuation in $(A,H)$ of ``short'' nodoids, starting from a cylinder. First yielding classical liquid bridges, but eventually turning into self intersecting nodoids, requiring restarts. 
See Figs.~\ref{nodf0},\ref{nodf00}, and {\tt cmds1plot.m} for the plotting. \\
{\tt cmds1A2t.m}&Relating the numerical BPs from {\tt cmds1.m} to Theorem \ref{t-nodbif}.\\
{\tt cmds2.m}&Continuation in $(A,H)$ of ``long'' nodoids with one inner 
loop, see Figs.~\ref{nodf2},\ref{nodf2b}.\\ 
{\tt bdmov1.m, bdmov2.m}&scripts to make movies of Fig.\,\ref{nodf0} and 
Fig.\,\ref{nodf2}, see also \cite{sig}. \\ 
  \hline
{\tt nodinit.m, nodinitl.m}&Initialization of ``short'' and ``long'' nodoids\\
{\tt sGnodD.m, sGnodDjac.m}&rhs with DBCs, and Jacobian.\\
 {\tt qfArot.m}&area and rotational constraints, see {\tt qjacArot} for 
derivative.\\
 {\tt getN.m}&overload of {\tt getN} to flip $N$.\\
{\tt coarsufu.m}&mod of {\tt stanufu.m} for adaptive coarsening, see 
Fig.\,\ref{nodf2b}. 
\end{tabular}
}
\end{table}
\teskip 

For solutions without axial symmetry 
we additionally need to set a rotational phase condition (PC): 
If $X$ is a solution to \reff{feq}, so is $R_\phi X$,  
where $\phi$ is the angle in the
$x$--$y$ plane, and 
\hueq{
R_\phi\vec{x}=\bpm \cos\phi & \sin\phi& 0\\
  -\sin\phi&\cos\phi &0\\0&0&1 \epm\vec{x}. 
} 
Thus, if $\pa_\phi(R_\phi X)|_{\phi=0}=
{1\over x^2+y^2}\left(-y\pa_x X+x\pa_y X\right)\in\R^3$ is non--zero, 
then it gives a non--trivial kernel of $L$, 
which makes continuation unreliable and 
bifurcation detection impossible. See, e.g., \cite[\S3.5]{p2pbook} 
for further discussion of such continuous symmetries. 
Here, to remove the kernel we use the PC 
\huga{\label{Hq1}
q(u):=%\spr{\spr{\pa_\phi X_0,X}}:=
\int_X \spr{\pa_\phi X_0,X_0+uN_0}\dd S 
=\int_X \spr{\pa_\phi X_0,N_0}u\dd S=:\int_X {\rm d}\phi\ u\dd S\stackrel!=0, 
}
where $X_0$ is from the last step, with normal $N_0$, where 
$\phi$ is the angle in the $x$--$y$ plane, and hence 
$\pa_\phi X=-X_2\nabla_{X_1} X+X_1\nabla_{X_2}X$, where $\nabla_{X_j}$ 
are the components of the surface gradient, cf.~\reff{coZ}. On the discrete 
level we thus obtain the linear function 
\huga{\label{qfrot1}
\text{$q(u)=({\rm d}\phi)^T u$, with derivative $\pa_u q= ({\rm d}\phi)^T$,}
} 
$\dd \phi=\spr{-X_2\nabla_{X_1} X+X_1\nabla_{X_2}X,N}$, node--wise, i.e., 
$\nabla_{X_j} X$ is interpolated to the nodes via {\tt c2P}, with 
Voronoi weights. 
We then add $\srot q(u)$ to $E$ from \reff{cmcen} with Lagrange
 multiplier $\srot$, 
and thus modify the PDE to 
$G(u):=H(u)-H_0+\srot {\rm d}\phi \stackrel!=0$. 
This removes the $\phi$--rotations of non-axisymmetric $X$ 
from the kernel of $\pa_u G(u)$, and, moreover, $|\srot|<10^{-8}$ 
for all the continuations below. 

Since the (algebraic) volume $V$ of self--intersecting nodoids is not intuitive, 
here we use continuation in area $A$ and $H$. Thus, we start with 
the constraint {\tt qfA} with derivative {\tt qjacA}. 
For non--axisymmetric branches we 
set up {\tt qfArot} and its derivative,   
where we put \reff{Hq1} as a second component of {\tt qfA}, 
and similarly for the derivatives, 
and when we bifurcate to a non-axisymmetric branch, we set {\tt p.nc.nq=2} 
(2 constraints, area and rotational phase) and {\tt p.fuha.qf=@qfArot}.

\subsubsection{Short nodoids}\label{snsec}  
Listing \ref{nodl0} shows how we initialize by either a cylinder ({\tt icsw=0}) 
or parametrization \reff{paranod}. Here, 
{\tt pde.grid} is a 2D rectangular FEM mesh, of which we use the second 
component as $\phi\in[-\pi,\pi]$. Lines 16--20  implement \reff{nodpar} 
and \reff{paranod} with {\em twice} the line $\phi=\pm\pi$ where 
the rotation around the $z$--axis closes. To obtain a mesh 
without duplicate points from this, in the 
last line of Listing \ref{nodl0} we use {\tt clean\_mesh} from the \gptool. 

\hulst{caption={{\small From {\tt nodDBC/nodinit.m}; setting initial {\tt p.X} 
as a cylinder ({\tt icsw==0}) or via the parametrization 
\reff{nodpar} and \reff{paranod} with (x=t) 
and subsequent removal of duplicate points.}},
label=nodl0,linerange=12-23,firstnumber=12,stepnumber=5}{\dhome/nodinit.m}

\hulst{caption={{\small From {nodDBC/\tt cmds1.m}; 
branch switching at double BP, and continuation 
with rotational PC.}},label=nodl1,linerange=12-17,
firstnumber=10, stepnumber=5}{\dhome/cmds1.m}

\begin{figure}[ht]
\centering
\btab{l}{
\btab{ll}{{\sm (a)}&{\sm (b)}\\[-1mm]
\hs{-5mm}\btab{l}{\ig[width=0.42\tew]{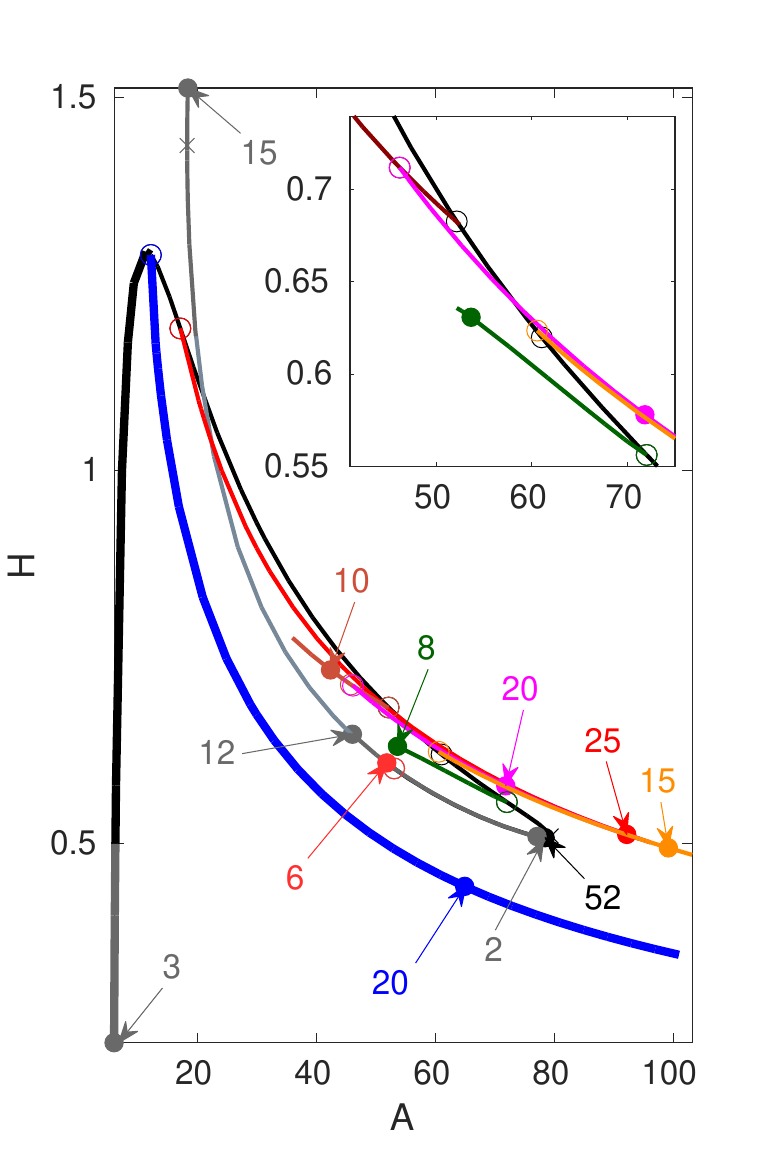}
%\\[-4mm]{\sm (c)}\\\hs{2mm}\ig[width=0.31\tew]{nopi1/at2}
}
&
\hs{-3mm}\rb{0mm}{\btab{l}{
\hs{0mm}\ig[width=0.28\tew]{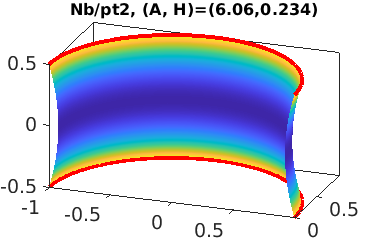}
\hs{0mm}\ig[width=0.3\tew]{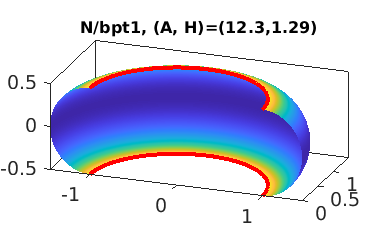}\\
\hs{0mm}\ig[width=0.28\tew]{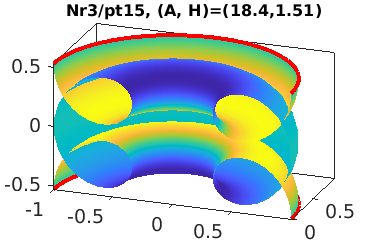}
\hs{0mm}\ig[width=0.3\tew]{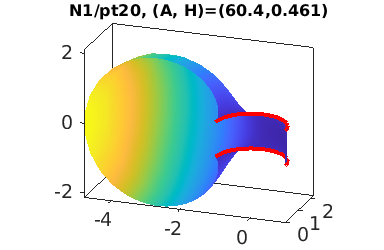}\\
\hs{0mm}\ig[width=0.29\tew]{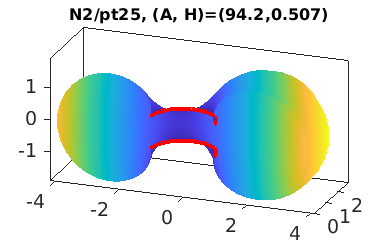} 
\hs{0mm}\ig[width=0.29\tew]{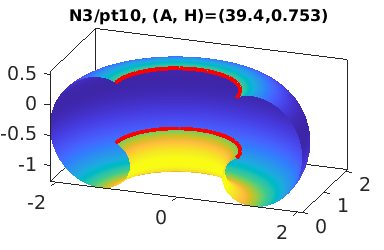}}}}\\[-3mm]
{\sm (c)}\\
\hs{-4mm}\ig[width=0.29\tew]{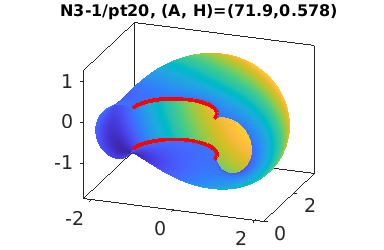}
\hs{-3mm}\ig[width=0.29\tew]{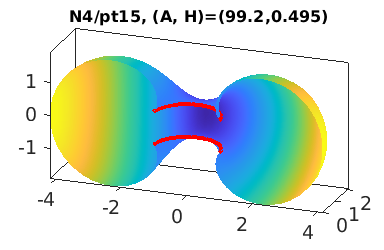}
\hs{-3mm}\ig[width=0.24\tew]{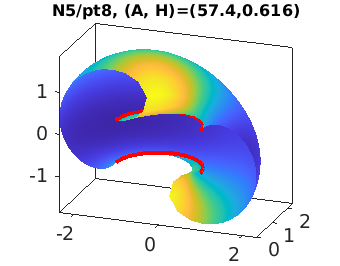}
\hs{-3mm}\ig[width=0.29\tew]{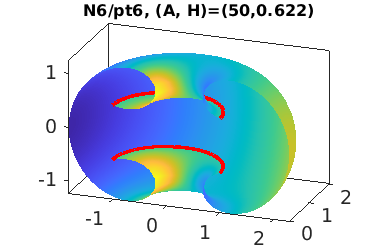}
}
%\ece 
\vs{-2mm}
\caption{\small Bifurcation diagram of (mostly) 
embedded nodoids (a), with samples in (b,c) cut open at the $x$--$z$ plane 
($y=0$). %, and (c) relating the numerical results to Theorem \ref{t-nodbif}, see Remark \ref{koirem}. 
Branches {\tt N} (black),  {\tt Nb} (grey), 
{\tt N1} (blue), {\tt N2} (red), 
{\tt N4} (orange), {\tt N5} (green), {\tt N6} (light blue), 
{\tt N3-1}  (magenta), and {\tt Nr1}, {\tt Nr2} and {\tt Nr3} (``restarts'' of {\tt N}, grey). 
See text for details, and Fig.~\ref{nodf00} for plots of 
{\tt N/pt52}, {\tt Nr1/pt2}, and {\tt Nr2/pt12}. 
\iffalse 
Branch {\tt N} (black) starts with the cylinder at $(A,H)=(2\pi,0.5)$; 
increasing $A$, at {\tt N/bpt1} there bifurcates 
to {\tt N1} (blue) with $m=1$ angular dependence and change of stability. 
Up to the 
fold at $(A,H)\approx (79.46,0.51)$  there are 4 more BPs to 
{\tt N2} (red, $m=2$), {\tt N3} (dark red, broken $z\mapsto -z$ symmetry), 
{\tt N4} (orange, $m=2$ again), {\tt N5} (green). 
On {\tt N3} there are secondary bifurcations, and we follow the 
first {\tt N3-1}  (magenta). After the fold we do a restart 
(see Fig.~\ref{nodf00}) and find several further BPs, e.g., 
{\tt N6} (light blue). 
\fi 
\label{nodf0}}
\end{figure} 

Figure \ref{nodf0} shows results from {\tt cmds1.m} (see also the movie 
{\tt nodDBCs.avi} from \cite{sig} to go step by step through the 
bifurcation diagram). We start at the cylinder 
and first continue to larger $A$ (black branch {\tt N}). 
The first BP at $(A,H)\approx (12.24,1.29) $ is double with 
angular wave number $m=1$. 
We simply select one of the kernel vectors to bifurcate, 
and do two steps without PC 
(blue branch {\tt N1}, lines 11-13 of Listing \ref{nodl1}). Then we 
switch on the rotational PC in line 14 and continue further.  
As predicted, BP1 occurs when $X$ meets the 
lower and upper boundary circles horizontally, and the stability 
changes from {\tt N} to {\tt N1}.%
\footnote{\label{stfoot}
{\tt N} up to {\tt BP1}, {\tt Nb}, and {\tt N1} are the only 
stable (in the sense of VPMCF) 
branches in Fig.\,\ref{nodf0}, and hence 
physically most relevant; the further branches we compute are 
all unstable, and hence of rather 
mathematical than physical interest.}
 The second BP yields the $m=2$ branch {\tt N2} (red). 
These results fully agree with those from \cite{bru18}. The branch 
{\tt Nb} (grey,  with {\tt pt3}) is the continuation of {\tt N} 
to smaller $A$ (and $V$), where the cylinder curves inward.

\begin{figure}[ht]
\centering
\btab{ll}{{\sm (a)}&{\sm (b)}\\[-1mm]
\hs{-5mm}\ig[width=0.5\tew]{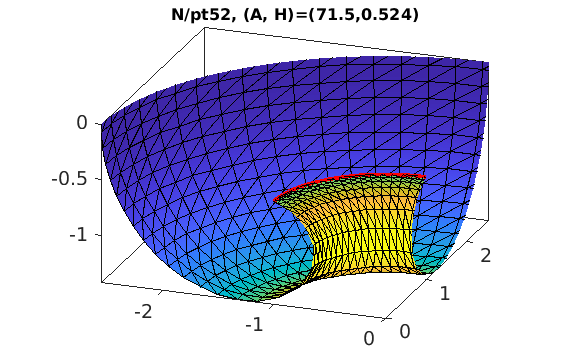}
&\hs{-7mm}\ig[width=0.5\tew]{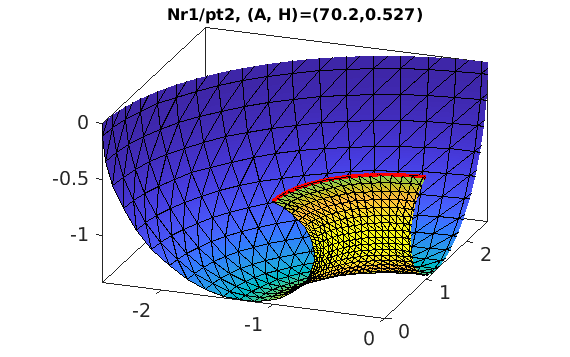}\\
{\sm (c)}&{\sm (d)}\\[-1mm]
\hs{-5mm}\ig[width=0.5\tew]{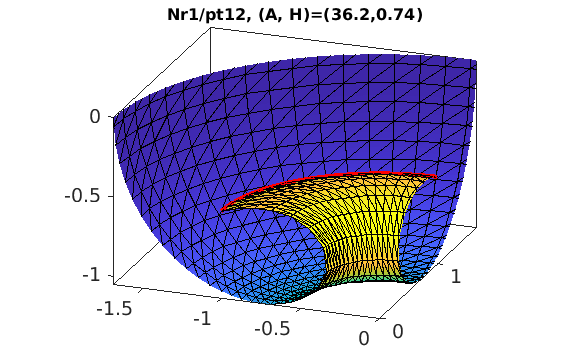}
&\hs{0mm}\ig[width=0.43\tew]{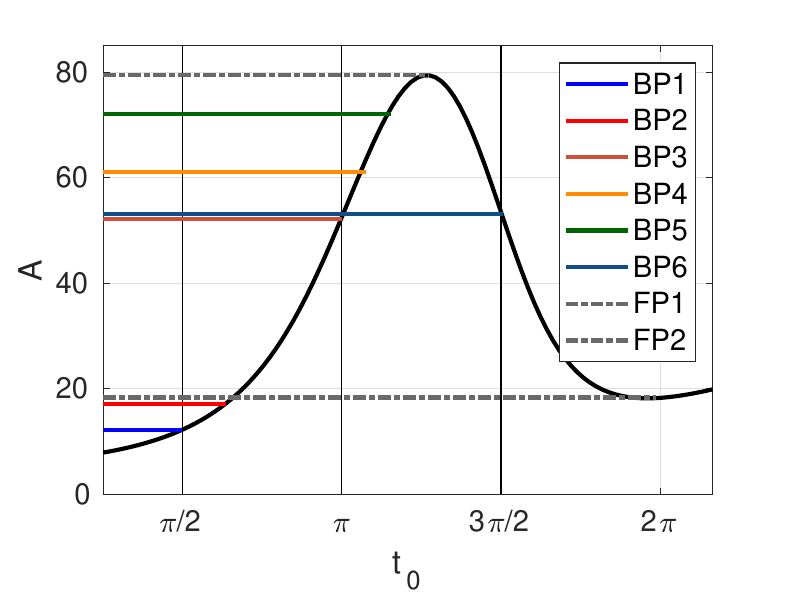}
}
\vs{-3mm}
\caption{\small (a--c) Continuation of Fig.~\ref{nodf0}; plots 
of ($1/8$th of) solutions on 
{\tt N} before and after remeshing; 
{\tt Nr2} from Fig.~\ref{nodf0} is from remeshing {\tt Nr1/pt12}, 
and   {\tt Nr3} from remeshing {\tt Nr2/pt10}. 
(d) Results from {\tt cmds1A2t.m}, see Rem.~\ref{koirem}(b). 
\label{nodf00}}
\end{figure}

The third BP on {\tt N} is simple with $z\mapsto -z$ symmetry breaking, 
yielding branch {\tt N3} (brown).  
On {\tt N3} there are secondary bifurcations, and following 
the first we obtain {\tt N3-1} (magenta). 
The 4th BP on {\tt N} again has $m=2$ but is different from the 2nd BP on {\tt N} as 
the nodoid has already ``curved in'' at the boundary circles, 
which is inherited by the bifurcating branch {\tt N4} (orange). 
The 5th BP on {\tt N} yields a skewed $m=2$ nodoid {\tt N5} (green).%
\footnote{BP5 is an example of a BP qualitatively predicted in 
\cite[Prop.3.9]{KPP17} at large $t_0$.\label{koif3}}
After the fold, the mesh in {\tt N} becomes bad at the necks, 
see {\tt N/pt52} in Fig.\,\ref{nodf00}. Thus, for accurate 
continuation we use \reff{paranod} to remesh, see {\tt Nr1/pt2} and 
Remark \ref{koirem}(a) and Fig.\,\ref{nodf00}(a--c), 
yielding the branch {\tt Nr1} (grey) in Fig.\,\ref{nodf0}(a). 
{\tt Nr1/pt12} in Fig.\,\ref{nodf00} shows that as after a number of 
steps the nodoid bulges further in, the mesh at the neck deteriorates again, and 
so we remesh again to {\tt Nr2} (light grey). 
The nodoid then self--intersects at $(A,H)\approx (22.9,1.05)$, 
and at {\tt Nr2/pt10} we do the next restart to {\tt Nr3}. Using such 
remeshing we can continue the branch {\tt N} (as {\tt Nr1, Nr2, Nr3, \ldots}) 
to many loops and self--intersections, with many further BPs 
as predicted in Lemma \ref{l-deins}.  In any case, 
although by branch switching from 
{\tt Nr1/bpt1} instead of from {\tt N/bpt6} 
we use a somewhat adapted mesh to compute 
branch {\tt N6} (red), we only compute a rather short segment of {\tt N6} 
because on {\tt N6} we quickly run into bad meshes again. See 
also \S\ref{lnsec} for further comments/experiments on the meshing of nodoids. 
In Fig.\,\ref{nodf00}(d) we illustrate 
the correspondence of our numerical results for the continuation in 
$A$ to Theorem \ref{t-nodbif}, see Remark \ref{koirem}(b). 

\brem\label{koirem}{\rm a) For axi-- and $Z_2$ symmetric  
nodoids, we can easily extract $a{=}(2HR{-}1)^2-1$ from our numerical data, 
with $R$ the radius on the $z{=}0$ plane. We can then numerically solve 
the second equation in \reff{lenrad}, i.e., 
$\ds 1=r^*=\frac{\cos t_0+\sqrt{\cos^2 t_0+a}}{2 \abs{H}}$ for $t_0$, 
and use this for restarts with a new mesh, for instance from {\tt N/pt52} 
to {\tt Nr1/pt1} in Fig.\,\ref{nodf00}. 

b) Similarly, given $r^*=1$ and $l=0.5$, we can solve 
\reff{impnod} for $a$ and $t_0$ in a continuation process, see 
{\tt cmds1A2t.m}. Then computing $A=A(a,t_0)$ 
gives the black curve in Fig.\,\ref{nodf00}(d), and 
intersecting the $A$ values of our numerical BPs 
gives the $t_0$ values for BP1, BP3 and BP6 as predicted, 
and explains the folds FP1 and FP2. In summary, the BPs on {\tt N}, 
their multiplicities, 
and their relation to Theorem \ref{t-nodbif} (if applicable) are
\huga{\label{koitab}
\text{\btab{r|cccccc}{
BP number&BP1&BP2&BP3&BP4&BP5&BP6\\
multiplicity&2&2&1&2&2&2\\
Theorem \ref{t-nodbif}&1.&NA&2.&NA&NA&1.\\
$t_0$&$\pi/2$&1.995&$\pi$&3.377&3.622&$3\pi/2$
}}
} 
where NA means not applicable, and where for BP1, BP3 and BP6 we give the exact values, with as indicated in 
Fig.\,\ref{nodf0}(c) very good agreement of the numerics.%
\footnote{This also holds for further BPs and folds, but we 
refrain from plotting these in the already cluttered BD in Fig.\,\ref{nodf0}.} 
}\eex\erem

\subsubsection{Long nodoids}\label{lnsec} 
\begin{figure}[ht]
\btab{ll}{{\sm (a)}&{\sm (b)}\\
\hs{-2mm}\ig[width=0.26\tew]{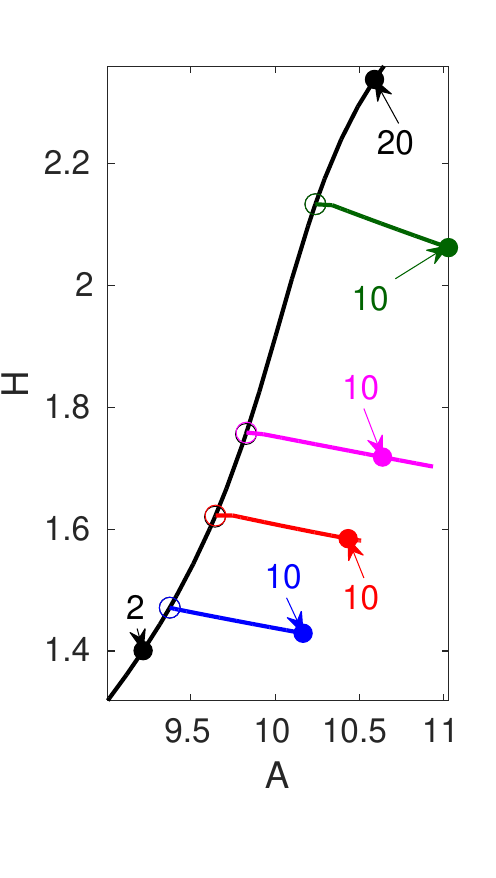}
&
\hs{-3mm}\rb{40mm}{\btab{ll}{
\hs{0mm}\ig[width=0.25\tew]{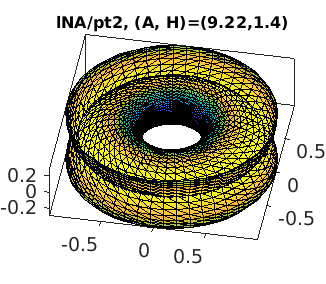}
\hs{0mm}\ig[width=0.25\tew]{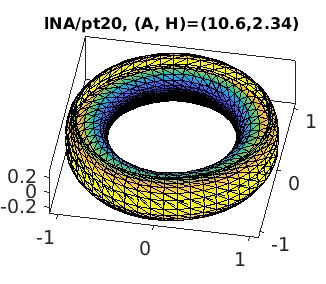}
\hs{0mm}\ig[width=0.25\tew]{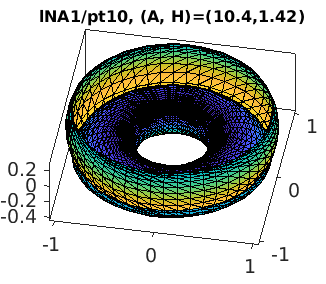} \\ 
\hs{0mm}\ig[width=0.24\tew]{nopi1/lNA2-10b}
\hs{0mm}\ig[width=0.25\tew]{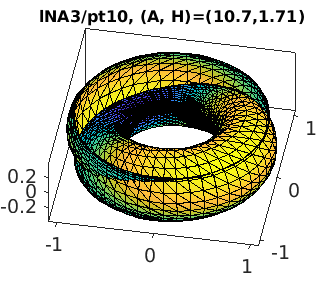}\hs{0mm}
\ig[width=0.23\tew]{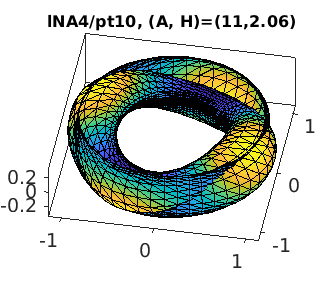}
}}}
\vs{-6mm}
\caption{\small (a) Bifurcation diagram of self--intersecting 
nodoids; 
branch {\tt lNA} (black) starts near $(A,H)=(9,1.3)$  
and shows four BPs to {\tt lNA1} (blue, broken $z$--symmetry), 
{\tt lNA2} (red, $m=2$), 
{\tt lN3A} (magenta, $m=1$), and {\tt lNA4} (green, $m=4$). Samples in (b). 
\label{nodf2}}
\end{figure}

In {\tt nodDBC/cmds2a.m} and Fig.\,\ref{nodf2}  
(see also {\tt nodDBCs.avi} from \cite{sig})  
we consider ``long'' nodoids with self--intersections.%
\footnote{\label{stfoot2}All of the branches from Fig.\,\ref{nodf2} 
are unstable, i.e., Footnote \ref{stfoot} applies more strongly.} 
As a slightly more explicit alternative to 
\reff{paranod}, in {\tt nodinitl.m} we now parameterize an initial 
axisymmetric nodoid $\tilde\CN_{r,R} $ following \cite{Mla02} by
\hual{\label{inodpara} 
X:[-\pi/2,\pi/2]\times[0,2\pi]&\rightarrow\R^3, \qquad 
  (x,\bphi)\mapsto \bpm\frac{r}{\del(x,k)}\cos(\bphi)\\
  \frac{r}{\del(x,k)}\sin(\bphi)\\ R
  E(x,k)-rF(x,k)-Rk^2\frac{\sin(x)\cos(x)}{\del(x,k)}\epm,
} 
where $r$ and $R$ are the neck and the buckle radius, $F$ and $E$ are the 
elliptic integrals of the first and second kind, $k=\sqrt{(R^2-r^2)/R^2}$, and 
$\del(x,k)=\sqrt{1-k^2\sin(x)}$.  
 It turns out that 
here we again need to be careful with the meshes, and besides adaptive 
mesh refinement we also use suitable initial meshes. We discretize the 
box $[-\pi/2,\pi/2]\times[0,2\pi]$ (pre-image in \reff{inodpara}) 
by Chebychev nodes in $x$ and equidistant nodes in $y$. This is 
implemented in  a slight modification 
of {\tt stanpdeo2D.m} in the current directory, and adapted 
to the parametrization \reff{inodpara}, 
which ``contracts'' the mesh for the loop in the middle. 

\begin{figure}[ht]
\btab{ll}{%{\sm (a)}&{\sm (b)}\\
\rb{12mm}{\hs{8mm}\ig[width=0.41\tew]{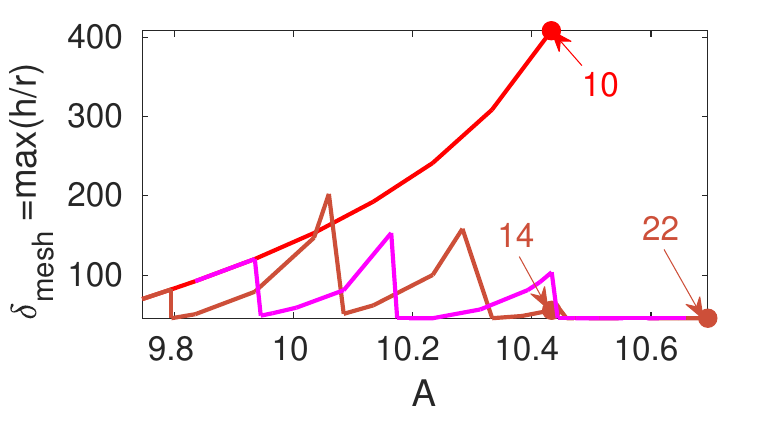}}&
\hs{-4mm}\ig[width=0.51\tew]{nopi1/lNA2-10z}\\[-11mm]
{\sm (c)}\\[-6mm]
\hs{-4mm}\ig[width=0.51\tew]{nopi1/lNA2c14}&\hs{-4mm}
\ig[width=0.51\tew]{nopi1/lNA2c22}
}
%\vs{-12mm}

\vs{-110mm}{\sm (a)}\hs{90mm}{\sm (b)}\\[40mm]
{\sm (c)}\hs{90mm}{\sm (d)}\\[45mm]

\caption{\small Results from the end of {\tt cmds2.m}, example of a {\tt degcoarsenX}--{\tt cont} loop (lNA2c, brown), 
and adaptive coarsening (lNA2cc, magenta) for 
{\tt lNA2}. 
(a) mesh quality $\delm=\max(h/r)$ over $A$, original {\tt lNA2} 
in red. 
(b) original {\tt lNA2/pt10} (cut open), $n_p{=}3430$; 
(c,d)  samples from {\tt lNA2c} with $n_p=2744$ and $n_p=2347$. 
\label{nodf2b}}
\end{figure}

As we continue the axisymmetric branch {\tt lNA} (black) to larger $A$, 
the inner loop ``contracts and moves out'', 
cf.~{\tt lNA/pt2} vs {\tt lNA/pt20} in Fig.\,\ref{nodf2}(b). 
Along the way we find several BPs, the first 
yielding a branch {\tt lNA1} (blue) with broken $z\mapsto -z$ symmetry. 
The next three BPs yield branches with angular wave numbers $m{=}2, m{=}1$, and 
$m{=}3$. As we continue these branches to larger $A$, the mesh quality 
deteriorates due to very acute triangles where the inner loop 
strongly contracts. 
This suggests coarsening by removing degenerate triangles, 
which we exemplarily discuss in  Fig.\,\ref{nodf2b}, see also 
the end of {\tt cmds2.m}. The red line in (a) shows the 
mesh--distortion along {\tt lNA2}, and (b) shows a zoom 
of {\tt lNA2/pt10}; the very acute triangles on the inner loop 
($\delm\approx 400$ at {\tt pt10}) lead to stepsize reduction 
and eventual continuation failure. The brown line in (a) and the 
samples in (c,d) show results from the {\tt degcoarsenX}--{\tt cont} loop
\bce
{\tt for i=1:6; p=degcoarsenX(p,sigc,nit,keepbd); p=cont(p,4); end;} 
\ece
\vs{-0mm}
with   ${\tt sigc}{=}0.5, {\tt it}{=}6,
{\tt keepbd}{=}1$ (cf.~footnote \ref{dcfoot})
starting at {\tt lNA2/pt3}. The distortion stays smaller 
(with $\delm\approx 50$ actually at the boundary), and the continuation 
runs faster and more robustly (larger stepsizes 
feasible) than on the original mesh. The magenta line is 
from setting {\tt p.fuha.ufu=@coarsufu} which adaptively coarsens 
(cf.~{\tt refufu.m} for refinement in Fig.\,\ref{spcf1c}) when 
$\delm$ exceeds $100$.  
Both here yield quite similar results, and naturally, similar use 
of {\tt degcoarsenX} is also useful for the other nodoids 
from Fig.\,\ref{nodf2}, and for those from Fig.\,\ref{nodf0}, 
in addition to the very specific remeshing used 
there, which is only possible because of the explicit formulas. 
Nevertheless, 
we remark again that the parameters for {\tt degcoarsenX} 
need trial and error for robustness and efficiency.

\subsection{Nodoids with pBCs in $z$} % (demo {\tt nodpBC})}
\label{nodnumpBC}
\def\dhome{geomtut/nodpBC} 
In \cite{MaPa02}, 
bifurcations of axisymmetric to non--axisymmetric nodoids are studied with 
the period (the ``height'') along the axis of revolution  
(wlog the $z$--axis) 
as the continuation/bifurcation parameter. 
This uses a different parametrization 
of the nodoids than \reff{paranod} or \reff{inodpara}, which we do not 
review here, as we shall again use \reff{paranod} for the 
initialization. For fixed $H=1$, 
\cite{MaPa02} proves that there is a $r_0>0$ such that 
for neck radii $r>r_0$ ($r<r_0$) there are (are not) bifurcations from nodoids, 
and gives detailed asymptotics 
of bifurcation points in a regime ($\tau\to-\infty$ in \cite{MaPa02}) 
which corresponds 
to $(R-r)/R\to 0$ with outer radius $R$, 
see below. In particular, the 2nd variation of the area functional 
around a given nodoid $\CN_\tau$ is analyzed with $z\in\R$, i.e., 
for the full non--compact nodoid, not just for one period cell. 
This proceeds by Bloch wave analysis, and first establishes the
 band structure of the spectrum. 
Using a parametrization similar to 
\reff{inodpara}, a detailed analysis of the 
second variation of the area functional, and ultimately two 
different numerical methods, \cite{ross05} shows that $r_0=1/2$, 
and the first bifurcation (i.e., at $r_0$) leads to 
non--axisymmetric nodoids with angular wave number $m=2$ and same periodicity 
in $z$, i.e., Bloch wave number $\al=0$ in \cite{MaPa02}. 

\begin{figure}[ht] 
\centering
\btab{lll}{{\sm (a)}&{\sm (b)}&{\sm (c)}\\
\hs{-2mm}\rb{4mm}{ \ig[width=0.28\tew]{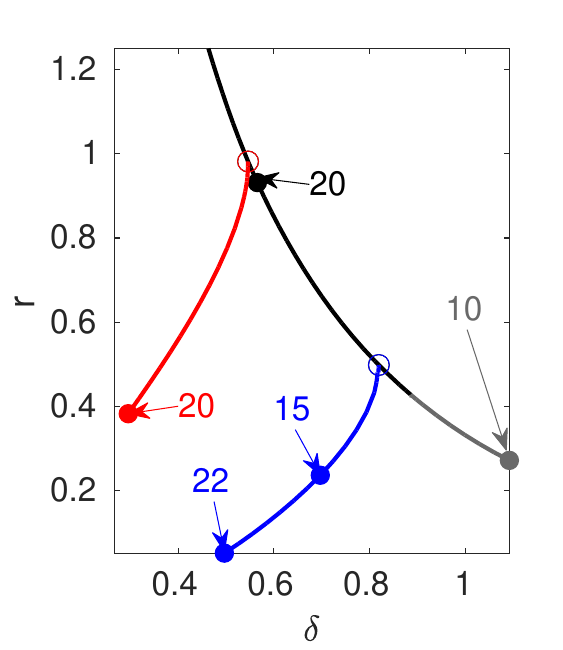}}
%\hs{-3mm}\rb{30mm}{\btab{ll}{
&\hs{-12mm}\ig[width=0.36\tew]{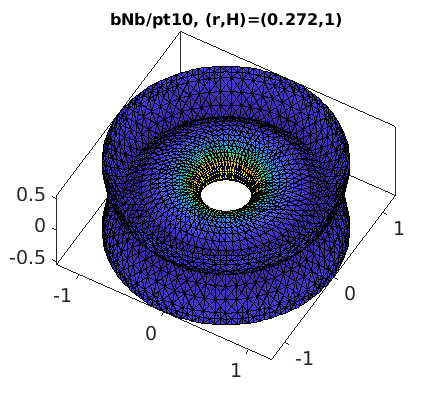}
&\hs{-5mm}\ig[width=0.36\tew]{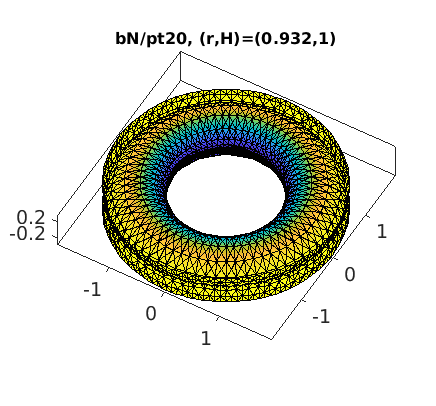}\\[-5mm]
{\sm (d)}&{\sm (e)}&{\sm (f)}\\
\hs{-4mm}\ig[width=0.36\tew]{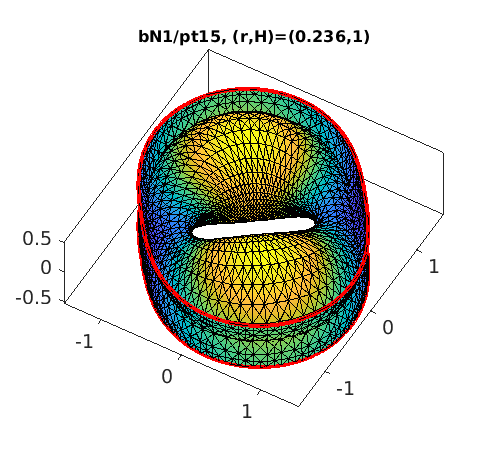}
&\hs{-5mm}\ig[width=0.36\tew]{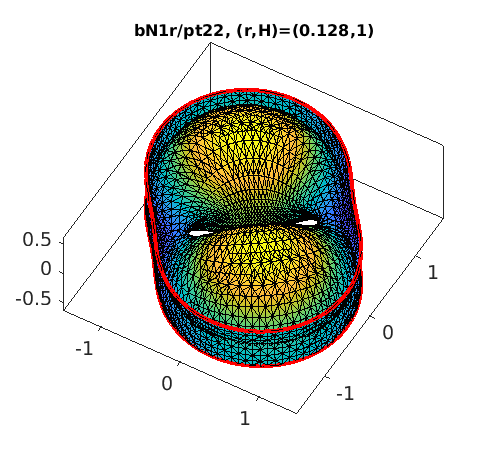} 
&\hs{-5mm}\ig[width=0.36\tew]{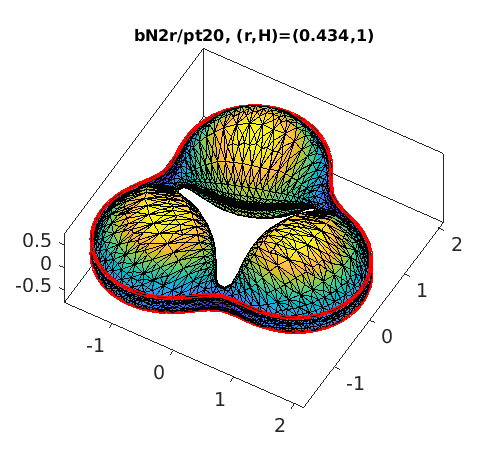}
}
\vs{-6mm}
\caption{\small (a) Bifurcation diagram of nodoids parametrized by 
height $\del$, fixed $H=1$. The axisymmetric 
branch {\tt bN} (black) starts near $\del=0.88$ via \reff{paranod}, 
and in direction of decreasing $\del$ shows a sequence of 
BPs to nodoids with broken $S^1$ symmetry, here {\tt bN1} (blue, $m=2$) and 
{\tt bN2} (red, $m=3$). Samples in (b--f), with {\tt bN1r} and {\tt bN2r} 
after some refinement. 
\label{nodf3}}
\end{figure}

\begin{figure}[H]
\centering
\btab{lll}{{\sm (a)}&{\sm (b)}&{\sm (c)}\\
\hs{-2mm}\rb{4mm}{ \ig[width=0.33\tew,height=50mm]{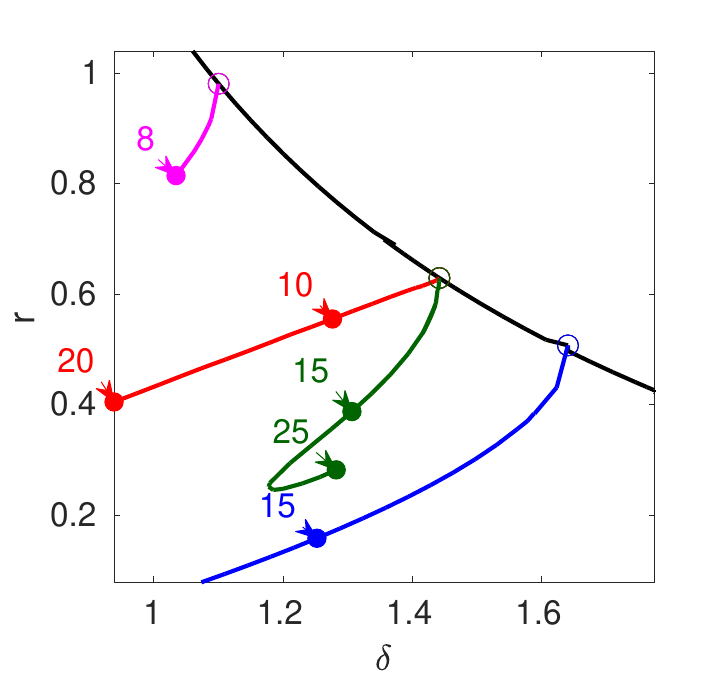}}
%\hs{-3mm}\rb{30mm}{\btab{ll}{
&\hs{-12mm}\ig[width=0.36\tew]{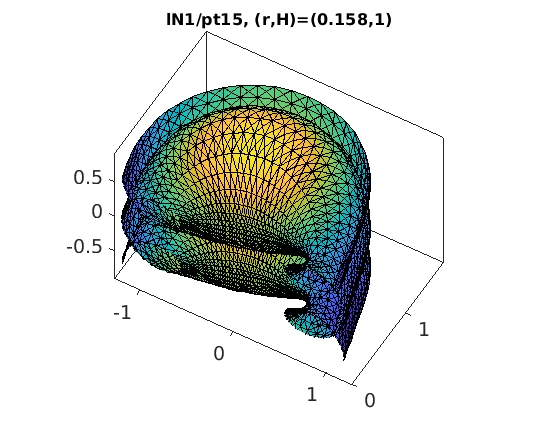}
&\hs{-8mm}\ig[width=0.36\tew]{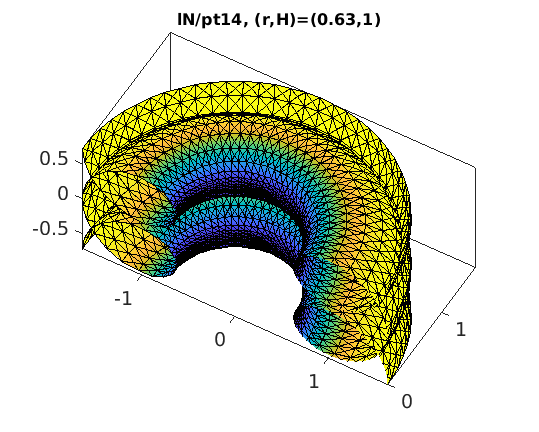}\\[-5mm]
{\sm (d)}&{\sm (e)}&{\sm (f)}\\[-1mm]
\hs{-4mm}\ig[width=0.37\tew]{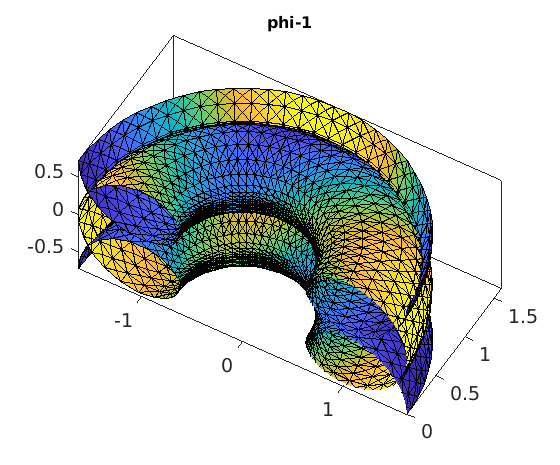}
&\hs{-5mm}\ig[width=0.37\tew]{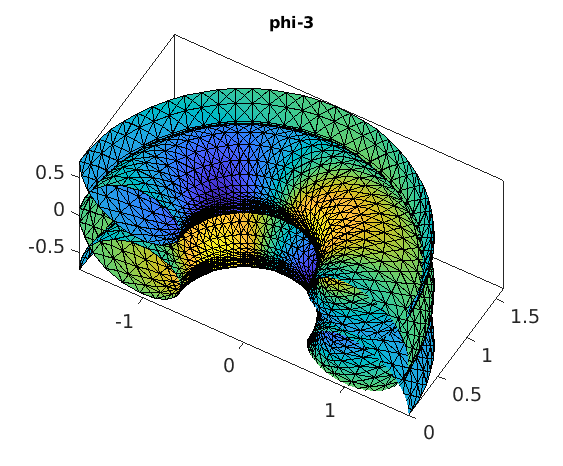} 
&\hs{-4mm}\rb{-3mm}{\ig[width=0.34\tew]{mpi3/ln2a-10b}}\\[-5mm]
{\sm (g)}&{\sm (h)}&{\sm (i)}\\[-1mm]
\hs{-5mm}\ig[width=0.35\tew]{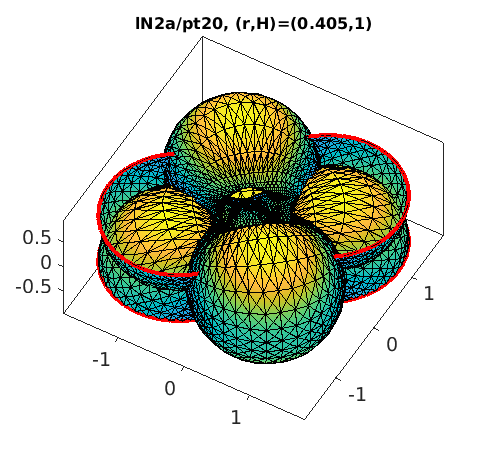} 
&\hs{-5mm}\ig[width=0.35\tew]{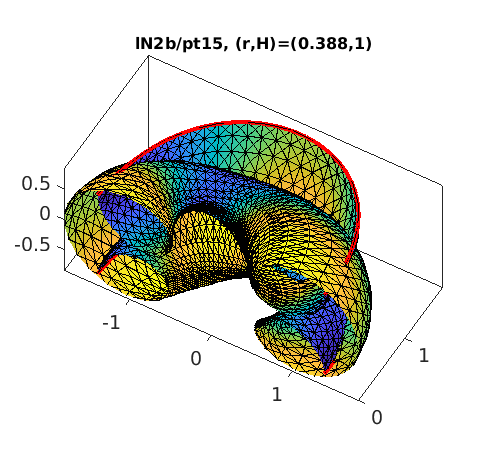}
&\hs{-5mm}\ig[width=0.35\tew]{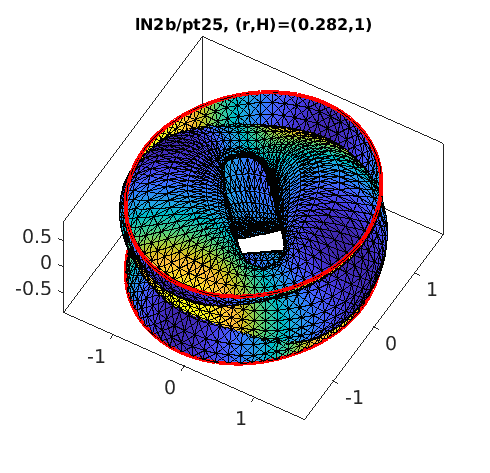} 
}
\vs{-6mm}
\caption{\small (a) Branching from the $S^1$ nodoid from Fig.~\ref{nodf3} 
on twice the minimal cell, see text for details.
\iffalse giving a new (fourfold) BP at $r\approx 0.63$. 
Solutions on lN1 (blue) and lN3 (magenta) consist of two copies 
of the respective 
solutions from Fig.~\ref{nodf3}, like in (b) (plot cut open). 
Two kernel vectors at BP2 (see (c)) 
are in (d),(e);  $\phi_1$ yields {\tt ln2a} (f,g), with broken $D_4$ 
symmetry, and the nodoid bulging out around 
$\vt=\pi/4$ and $\vt=5\pi/4$ in the $z=\pm\del$ plane, and at 
$\vt=3\pi/4$ and $\vt=7\pi/4$ in the $z=0$ plane. The behavior 
in $\phi_3$ and accordingly on {\tt ln2b} 
is similar but rather on the inner buckle, leading to 
a more twisted shape, and ultimately to continuation failure (with 
the given mesh) shortly after {\tt ln2b/pt25}. The other 
eigenvectors $\phi_2,\phi_4$ are related to $\phi_1,\phi_3$, respectively, 
by rotation around $z$. 
\fi 
\label{nodf4}}
\end{figure}

Here we also consider periodic (in $z$) nodoids with fixed $H=1$ using 
the height $\del$ as continuation/bifurcation parameter. We recover 
the primary bifurcation at $r=r_0=1/2$ from \cite{ross05}, and 
further bifurcations, see Figs.~\ref{nodf3} and \ref{nodf4}. 

\brem\label{nprem}{\rm Similar to \S\ref{nodnumDBC} we distinguish between 
``short'' and ``long'' nodoids. Here, this merely corresponds to 
computing on one respectively two period cells in $z$, and the 
main distinction is as follows: All 1--periodic solutions are 
naturally $n$--periodic for any $n\in\N$. With respect to bifurcations, 
the 1--cell computations then correspond to Bloch wave numbers $\al=0$ 
in \cite{MaPa02}. For $n\ge 2$ periods cells we obtain 
further discrete Bloch wave numbers, e.g., additionally $\al=\pi$ 
for $n=2$. This then allows bifurcations which simultaneously 
break the $S^1$ and the $Z_2$ symmetry of the symmetric nodoid, and 
this is illustrated in Fig.\,\ref{nodf4}, which only gives a basic 
impression of the extremely rich bifurcation picture to be expected 
when the computational domain is expanded further in $z$. %In particular, 
To avoid clutter 
we refrain from putting the cases $n=1$ and $n=2$ in one figure.}
\eex\erem 

Numerically, to set up ``periodic boundary conditions in $z$'', 
we proceed similar to the \pdep\ setup for periodic boundary 
conditions on fixed (preimage) domains, see \cite[\S4.3]{p2pbook}. 
The basic idea is to identify 
points on $\pa X$ at $z=\pm \del$. Thus, before the main step 
$X_0\mapsto X_0+uN_0$ for all our computations, we 
transfer the values of $u$ from $\{X_3=-\del\}$ to $\{X_3=\del\}$ 
via a suitable ``fill'' matrix {\tt p.mat.fill}, which has to be generated at 
initialization and regenerated after mesh--adaptation. 
The essential command is {\tt box2per}, which calls {\tt getPerOp} 
to create {\tt p.mat.fill} (and {\tt p.mat.drop} which is used to 
drop redundant variables), and which rearranges $u$ by dropping 
the (redundant) nodal values at points which are filled by periodicity. 

Similar to \S\ref{nodnumDBC} we need a rotational PC for non--axisymmetric 
branches, but here for all computations 
we additionally need translational PCs in $x,y$ and $z$ directions, 
i.e. $S_i\vx=\vx+e_i$. These translations act infinitesimally 
in the tangent bundle 
as $S_i X_0=\nabla_i X_0$, and hence the pertinent PCs are 
\huga{\label{qf-trans}
q_i(u)=\spr{\nab_iX_0,X_0+uN_0}=\spr{\nab_i X_0,N_0}u, \quad i=1,2,3, 
}
with derivatives $\pa_u q_i(u)=\spr{\nab_iX_0,N_0}$. Like 
\reff{qfrot1}, they are implemented node--wise, and their derivatives 
are added to $G$ with Lagrange multipliers $s_x,s_y,s_z$. 
Table \ref{nodptab} comments on the files used, and Listings 
\ref{npbcl1}--\ref{npbcl4} show the main new issues from the 
otherwise typical function and script files.

\taskip
\begin{table}[ht]\caption{
Files in {\tt pde2path/demos/geomtut/nodpBC}; similar to {\tt nodDBC}, 
and we mainly comment on the differences.  
\label{nodptab}} 
\centering\vs{-2mm} {\small
\begin{tabular}{p{0.18\tew}|p{0.83\tew}}
  {\tt cmds1.m}&script for Fig.\,\ref{nodf3} (single period cell); experiments with mesh-adaptions 
in {\tt cmds1b.m}. \\
 {\tt cmds2.m}&script for Fig.\,\ref{nodf4} (double period cell).\\
{\tt lnodpBCmov.m}&script to make a movie of Fig.\,\ref{nodf4}, see also \cite{sig}. \\ 
{\tt cmds2b.m}& experiments with trying to solve 
quadratic(cubic) bifurcation equations at BP2, and checking mode crossing 
at BP2 by varying $H_0$. \\
  {\tt nodbuckinit.m}&initialization, using parametrization \reff{paranod}, 
extracts initial height $\del$, and sets the ``fill'' and ``drop'' matrices 
for the pBCs. \\
  {\tt sGnodpBC.m}&rhs with pBCs; here using numerical Jacobians.\\
  {\tt qf.m, qjac.m}&3 translational constraints (for $S^1$ symmetric branches), and derivative.\\
{\tt qfrot.m}&{\tt qf.m} extended by rotational phase condition, 
derivative in {\tt qjacrot.m}.\\
{\tt getN.m}&flips $N$, and additionally forces horizontal normals at bottom and top.\\
{\tt getM.m}&mod of standard {\tt getM} with subsequent {\tt filltrafo} 
for reduction to active nodes.\\
{\tt getMf.m}&renaming of standard {\tt getM} to get $M$ for the full $X$ 
(no dropping of per.~boundaries).\\
{\tt getA.m}&small mod of standard {\tt getA} to work on full $X$.\\
{\tt mytitle.m}&helper function for plots, called in {\tt pplot} since 
{\tt p2pglob.tsw=5}. 
\end{tabular}
} 
\end{table}
\teskip

\hulst{linerange=24-26,firstnumber=24, 
caption={\small Initializing pBCs in {\tt nodpBC/nodbuckinit.m} via 
{\tt box2per}, which automatically generates the pertinent {\tt drop} and {\tt 
fill} matrices to deal with the pBCs in direction 3 ($=z$--direction).},label=npbcl1}
{\dhome/nodbuckinit.m}
\hulst{linerange=1-4, 
caption={\small Start of {\tt nodpBC/sGnodpBC.m}, l4 fills u, and afterwards 
we proceed as before.},label=npbcl2}
{\dhome/sGnodpBC.m}
\hulst{linerange=1-9, 
caption={\small {\tt nodpBC/qf.m}, implementing the three translational PCs; 
the gradient matrices are computed on the full {\tt p.X}, and their  
periodic parts are dropped for their actions on {\tt u(1:p.nu)} (line 8). 
For non axi--symmetric nodoids we extend {\tt qf} in {\tt qfrot} 
by a rotational PC as an additional 4th line.},label=npbcl3}{\dhome/qf.m}
\hulst{linerange=5-15, firstnumber=5, 
caption={{\small Start of {\tt nodpBC/cmds1.m}. Initialization, first step, 
and initial mesh refinements. }},
label=npbcl4}{\dhome/cmds1.m}

Fig.\,\ref{nodf3} shows some results from {\tt cmds1.m}. For robustness 
(essentially due to the strong contractions at the inner loops later in 
the branches) it turns 
out to be useful to initialize with a rather coarse mesh and after 1 or 2 steps 
refine by area. As we then decrease $\del$ 
from the initial $\del\approx 0.88$, we find the first BP at 
$\del\approx 0.82$ and with $r=0.5$, corroborating \cite{ross05}, 
to the angular wave number $m=2$ branch {\tt bN1}. Using suitable 
mesh refinement along the way we can continue {\tt bN1} to small $\del$, 
where in particular we have multiple self--intersections; first, 
the inner loops extend the ``height'' $\del$ for $\del<\del_0\approx 0.78$, 
and second the inner loops intersect in the plane 
$z=0$ for $\del<\del_1\approx 0.43$ 
(not shown), making the inner radius $r=0$ (or rather undefined). 
The branch {\tt bN2} 
from the next BP at $\del\approx 0.54$ has $m=3$, and otherwise 
behaves like the $m=2$ branch.  All these branches are rather strongly unstable, with $\ind(X)>4$, and Footnotes \ref{stfoot} and \ref{stfoot2} 
again apply. 

As indicated in Remark \ref{nprem}, the branching behavior of the periodic 
nodoids very much depends on which period cell in $z$ we prescribe, 
with Fig.\,\ref{nodf3} corresponding to one cell.  
To illustrate the richness that can be expected for larger cells, 
in {\tt cmds2.m} and Fig.\,\ref{nodf4} 
we consider twice the minimal cell, see also {\tt nodpBCl.avi} from \cite{sig}. 
This yields the same primary nodoid branch, and as a subset of bifurcations the 
bifurcations from Fig.\,\ref{nodf3}, with two stacked copies of 
the solutions from Fig.\,\ref{nodf3} along all these branches. 
Additionally we have a new BP2 around $\del=1.44$, 
with small eigenvalues 
$\mu_{1,2}\approx 0.0003$ and $\mu_{3,4}\approx 0.004$, 
and the next eigenvalues are $\mu_{5}\approx -0.67$ (simple) and 
$\mu_{6,7}\approx 0.87$. 
The (approximate) kernel vectors associated to $\mu_{1,3}$ are 
$\phi_1, \phi_3$ given in Fig.\,\ref{nodf4}(d,e),  
and additionally we have $\phi_2=R_{\pi/2}\phi_1$ and $\phi_4=R_{\pi/2}\phi_1$ 
(rotation around the $z$ axis). From the 4 small eigenvalues separated 
from the rest of the spectrum we might guess that BP2 
is fourfold, which should have important consequences 
for the branching behavior at BP2, in the sense of ``mixed modes'', 
%similar to ``horizontal stripes + vertical stripes=spots'' \cite{E91}, 
see \cite[\S2.5.4]{p2pbook}. 
However, $\phi_1$ and $\phi_3$ do not seem related by any 
symmetry, and, moreover, using {\tt qswibra} and {\tt cswibra} 
(see {\tt cmds2b.m}) to search 
for solutions of the algebraic bifurcation equations at BP2 only 
yields the ``pure'' modes $\phi_1$ and $\phi_3$ (and their rotations).  
%and then the same bifurcating branches {\tt lN2a} and {\tt lN2b} as before. 
Thus we conclude that BP2 is {\em not fourfold}, but 
corresponds to two double BPs close together.%
\footnote{See also \S\ref{SPCsec}, where we discuss the opposite effect 
in more detail: 
There, an analytically double (due to obvious symmetry) 
zero eigenvalue is split up, 
with a larger split there than between $\mu_{1,2}$ and $\mu_{3,4}$ here.}

To further corroborate 
this, in {\tt cmds2b.m} we compute the ``{\tt lN}'' branches 
for $H_0=0.8$ and $H_0=1.1$. In both cases we find a similar 
spectral picture at the pertinent BPs as at BP2 for $H_0=1$ 
(4 small eigenvalues, well separated from the rest of the spectrum, 
with kernel vectors similar to Fig.\,\ref{nodf4}(d,e)), 
but the two pertinent pairs are themselves more clearly separated, 
and $\phi_1,\phi_3$ flip order between $H_0=0.8$ and $H_0=1.1$. 
If BP2 at $H_0=1$ was fourfold, then we would expect this 
to be due to symmetry, and hence to also hold for $H_0\ne 1$. That this 
is not the case suggests that near $H_0=1$ we rather have 
a ``mode crossing'' at BP2, and moreover do not expect 
bifurcating branches of ``mixed modes''.

Therefore, in {\tt cmds2.m} 
we simply do a branch switching in direction of the modes $\phi_1$ 
(see Listing \ref{npbcl5}) 
and $\phi_3$ separately, and obtain the branches {\tt lN2a} (red) and 
{\tt lN2b} (green).%
\footnote{\label{qtrickfn}
This trick can also be summarized as follows: We do not localize 
each of the close--together BPs near BP2, which would require very 
small (arclength) stepsizes {\tt ds}, and possibly many bisections. 
Instead, we just approximately localize some BP$^*$ near BP2, subsequently 
compute the (approximate) kernel at BP$^*$ using {\tt qswibra}, 
and then select a kernel 
vector and try branch--switching in that direction. This ``usually'' works 
(in particular it works here), and is in particular useful if many 
BPs (including BPs of higher multiplicity) are close together. See also 
\S\ref{SPCsec} for application of this trick in a different setting, and, 
e.g., \cite[\S10.1]{p2pbook} for further discussion.}
In detail, we do an initial step 
in direction $\phi_1$ (resp.~$\phi_3$), and then switch on the 
rotational PC as before. The branch {\tt lN2a} continues to small 
$\del$ without problems.  The branch {\tt lN2b} folds back near 
$\del=1.18$, and after {\tt pt25} continuation fails, due to 
mesh degeneration in the inner bends. This can be fixed to 
some extent  by careful mesh adaptation (yielding later continuation failure), 
but we do not elaborate on this here as Fig.\,\ref{nodf4} is mostly 
intended as an illustration of the  rich bifurcation behavior over larger 
period cells.

\hulst{caption={{\small Selection from {\tt nodpBC/cmds2.m}. Branch switching 
from BP2 via {\tt qswibra} and {\tt gentau}.}},linerange=27-37, firstnumber=27, label=npbcl5}
{\dhome/cmds2.m}

%\input{lnf3b3}

%\clearpage
\subsection{Triply periodic surfaces}\label{tpssec}
\def\dhome{geomtut/nodpBC} 
Triply periodic surfaces (TPS) are CMC surfaces in $\R^3$ 
which are periodic wrt three independent (often but not always 
orthogonal) directions. Triply periodic {\em minimal} surfaces (TPMS) 
(this implicitly also means embedded, sometimes abbreviated as TPEMS) 
have been studied since H.A.~Schwarz in the 19th 
century, and have found renewed interest partly due to the discovery 
of new TPMS by A.~Schoen in the 1970ies, 
and due to important (partly speculative) 
applications of TPMS (and their non--zero $H$ TPS companions) 
in crystallography, mechanics and biology.  
See for instance  \cite{AHLL87} and \cite{STFH06}, and 
\cite{TPMSbra} for a long list of TPMS. 
%From the applications point of view it is also natural 
%to study non--zero TPS, i.e., to vary $H$, cf.~\cite{...}.

From the PDE point of view, 
TPS solve \reff{feq} with periodic BCs on a bounding box. 
Some families of TPMS were studied as bifurcation problems in \cite{KPS18},  
using a cell length (period) in one direction as continuation/bifurcation 
parameter, and combined with numerical results from \cite{EjSh18}. 
Much of the theory of TPMS is based on Enneper--Weierstra\ss\ 
representations. 
See Remark \ref{Prem2}, where we relate some of 
our numerical results for the Schwarz P surface family to 
results from \cite{KPS18} obtained via Enneper--Weierstra\ss\ representations. 
A way to {\em approximate} TPS 
is as zeros of Fourier expansions of the form
\hugast{F(\vec{r})=\sum_{k\in\Z^3, |k|\le N}F(k)\cos(2\pi k
  \vec{r}-\alpha(\vec{r})). 
}  
A simple first order approximation of the Schwarz P surface 
(cf.~Fig.\ref{f1}(b)) is 
\hueq{\label{SPa}
\text{ Schwarz P surface }\approx
 \{(x,y,z)\in\R^3: \cos(x)+\cos(y)+\cos(z)=0\},
} 
Better approximations with some higher order terms are known, also 
for many other ``standard'' TPS, 
see, e.g., \cite{GBMK01} for a quantitative evaluation of such approximations. 
In the demo {\tt TPS} we focus on the Schwarz P family, and some 
CMC companions.%
\footnote{The approximation \reff{SPa}, 
and higher order corrections, also arise from solving the amplitude 
equations for a Turing bifurcation on a simple cubic (SC) lattice, 
where hence the Schwarz P surface, or, depending on volume fractions a 
CMC compagnion of Schwarz P, occurs as the phase separator 
between ``hot'' and ``cold'' phases. See, e.g., \cite{CKnob97} and \cite[\S8.1,8.2]{p2pbook}, and similarly \cite{WBD97} for the occurence of Scherk's surface in 3D 
Turing patterns.}

\taskip
\begin{table}[ht]\caption{
{\sm Selected files in {\tt TPMS/}, others 
({\tt getA.m,\,qf.m,\,qjac.m,\,updX.m}) are 
rather standard, with minor mods to account for the ``filling'' of 
the periodic boundaries. 
\label{tpstab}}}
\centering\vs{-2mm}
{\small 
\begin{tabular}{p{0.142\tew}|p{0.86\tew}} 
%file&remarks\\ \hline
{\tt cmds1.m}&Schwarz P family, relation to 
Weierstrass representation in {\tt cmdsaux.m}, Figs.~\ref{tpsf1} and \ref{tpsf2}.\\
{\tt cmds2.m}&CMC companions of Schwarz P, Fig.\,\ref{tpsf5}\\
{\tt Pinit.m}&Initialization, based on \reff{SPa} and {\tt distmesh}. \\
%{\tt getA, getV}&mods of standard to first fill $u$\\
{\tt getN.m}&mod of standard {\tt getN}, applying corrections at the boundaries of $X$, see 
Remark \ref{Prem1}.\\
%{\tt qf, qjac}&3 translational constraint, and derivative.\\
%{\tt sG}&scaling, fill, otherwise as usual\\
%{\tt updX}&mod of standard {\tt updX}; scaling!  \\
{\tt getperOpX.m}&here we also fill X; see also {\tt Xfillmat.m}. 
\end{tabular}
}
\end{table}
\teskip

\subsubsection{The Schwarz P minimal surface (family)}\label{SPsec}
In {\tt TPS/cmds1.m} we study continuation (and bifurcation) 
of the Schwarz P surface in the period $\del$ in $z$--direction, 
focusing on one period cell, i.e., the box 
\huga{\label{Qdeldef} 
B_\del:=[-\pi,\pi)^2\times [-\del/2,\del/2).
}
To get an initial (approximate) $X$ on $B_{2\pi}$, we use 
\reff{SPa} and the mesh 
generator {\tt distmesh} \cite{PeSt04}, on one eighth of $B_{2\pi}$, which 
we then mirror to $B_{2\pi}$. The continuation in $\del$ proceeds similar 
to \S\ref{nodnumpBC}, by first scaling $X=S_\del{\tt p.X}$ to period $\del$ 
in $z$ and then setting $X=X+uN$ and solving for $(u,\del)$. Subsequently, 
the same scaling is applied in {\tt updX} to set the new {\tt p.X}. 
As in \S\ref{nodnumpBC} we have translational invariance in $x,y$ and $z$, 
and hence exactly the same PCs, implemented in {\tt qf.m}, with derivatives 
in {\tt qjac.m}. 

Somewhat differently from \S\ref{nodnumpBC} we now also ``fill'' $X$ by 
taking the $\pa X$ values from the left/bottom/front of the box 
to the right/top/back of the box. 
While $u$ is stilled filled via $u={\tt  p.mat.fill*} u$, 
for filling $X$ we compute matrices {\tt p.Xfillx, p.Xfilly, p.Xfillz} 
(in {\tt Pinit.m}, via {\tt Xfillmat}, which calls {\tt getPerOpX}) 
similar to {\tt p.mat.fill}, but with 
$-1$ (instead of $1$) where we want to transfer $X$ values from one side 
of the box to the opposite side (assuming symmetry wrt the origin).  
%The functions {\tt getA.m} and {\tt getV.m} are local copies of their standard versions, with first filling $u$. 
Finally, it turns out that the continuation is slightly 
more robust if in {\tt getN} we correct $N$ at the boundaries to 
lie {\em in} the boundaries of $B_\del$, see Remark \ref{Prem1}. 

\begin{figure}[ht] %H]
\centering
\btab{llll}{{\sm (a)}&{\sm (b)}&{\sm (c)}&{\sm (d)}\\
\hs{-2mm}\ig[width=0.27\tew]{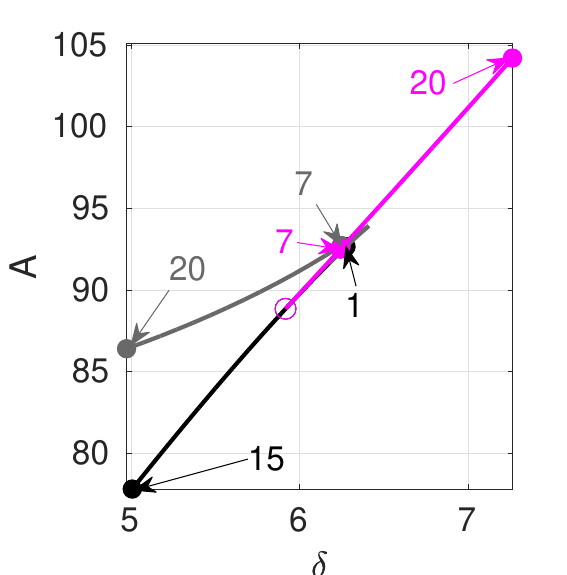}
&\hs{-5mm}\ig[width=0.26\tew]{tps/p-1}
&\hs{-5mm}\ig[width=0.26\tew]{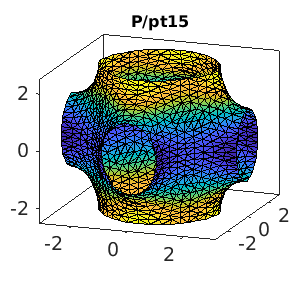}
&\hs{-5mm}\ig[width=0.26\tew]{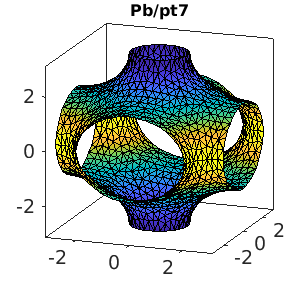}\\
{\sm (e)}&{\sm (f)}&{\sm (g)}&{\sm (h)}\\
\hs{-2mm}\ig[width=0.26\tew]{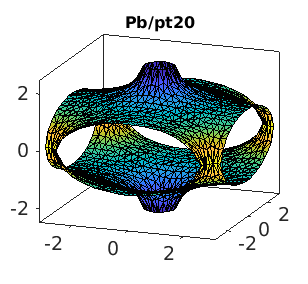}
&\hs{-4mm}\ig[width=0.26\tew]{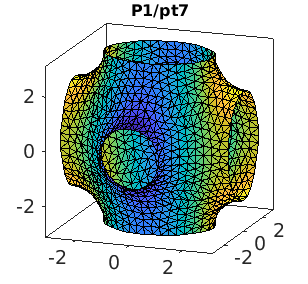}
&\hs{-3mm}\ig[width=0.26\tew]{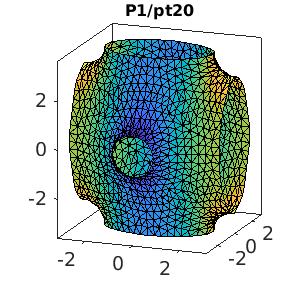}
&\hs{-4mm}\ig[width=0.24\tew]{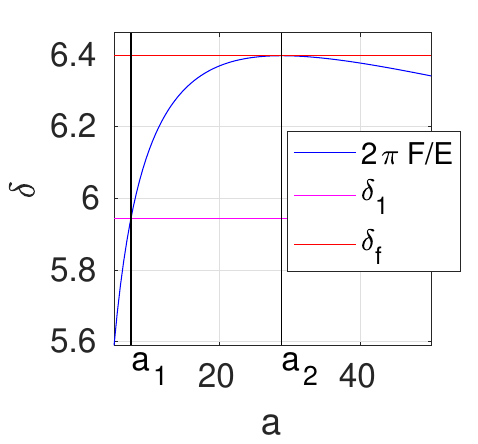}
}
\vs{-3mm}
\caption{\small (a) Bifurcations in the Schwarz P family, black (P) and 
grey (Pb) branch; bifurcating magenta branch (P1) breaks $D_4$ symmetry. 
Samples in (b--g). Comparison with \cite{KPS18} in (h), cf.~Remark 
\ref{Prem2}. 
\label{tpsf1}}
\end{figure}

Figure \ref{tpsf1} shows some results from {\tt cmdsP.m}. 
Decreasing $\del$ from $2\pi$ ({\tt P/pt1} in (b) at $\del=6.2732$), 
$X$ gets squashed in $z$ direction, and at $\del=\del_1\approx 5.9146$ 
we find a $D_4$ symmetry breaking pitchfork bifurcation (with the 
two directions corresponding to interchanging the $x$ and $y$ axis 
wrt shrinking and expansion) to a branch {\tt P1}, 
which then extends to large $\del$. On the other hand, increasing 
$\del$ from $2\pi$ (branch {\tt pB}, grey), we find a fold on the P branch at 
$\del=\del_f\approx 6.408$. Both $\del$ values agree well with 
results from \cite{KPS18} based on the Enneper--Weierstrass representation, 
summarized in Fig.\,\ref{tpsf1}(h), see Remark \ref{Prem2}. 

\brem\label{Prem1}{\rm a) The results from Fig.\,\ref{tpsf1} can also 
be obtained by choosing ``Neumann'' BCs on $\pa B_\del$. However, 
for other TPMS we need the pBCs. For instance, 
we can also continue the H surface family on a suitable (almost minimal) 
rectangular box, where however solutions fulfill pBCs but not Neumann BCs. 
 Due to the necessary larger period cell, and due to branch points of 
higher multiplicity, the numerics for the H family are more elaborate, 
and these results will be presented elsewhere. 

b) In fact, in the local copy {\tt TPS/getN.m} we apply a trick and 
zero out $N_1$ at $x{=}\pm \pi$, $N_2$ at $y{=}\pm\pi$ and $N_3$ at $z{=}\pm\del/2$. 
Thus, $N$ is forced to always lie in the cube's faces, yielding a 
``combination of NBCs and pBCs'' in the sense that the trick forces $X$ to meet 
the cube's faces orthogonally, while the pBCs keep $X$ on opposite faces together. 
However, the trick is for convenience as without it we get the same branches but in a less 
robust way, i.e., requiring finer discretizations and smaller continuation stepsizes.
 }\eex\erem 

%\subsubsection{Enneper--Weierstrass representation}\label{ewsec} 
\brem\label{Prem2}{\rm 
The Enneper--Weierstrass representation of a minimal 
surface is 
\huga{\label{ew1}
\bpm x,y,z\epm=
\mathrm{Re}\biggl[\er^{\ri \vt}
\int_{p_0}^p (1-z^2,\,\ri(1+z^2),\,2z)R(z)\dd z\biggr], 
}
$p_0,p\in \CM$ with $\CM$ a Riemannian surface, where $\vt$ is called  
Bonnet angle, and $R:\CM\to\C$ is called Weierstra\ss\ function. 
The Enneper surface $E$ from \S\ref{ensec} 
is given by the data $\CM=D_{\al}$ (disk of radius $\al$) and 
$R(z)\equiv 1$. 
For TPMS, $R$ is a meromorphic function, and $\CM$ 
consists of sheets connected at branch points 
given by poles of $R$. 
See, e.g., \cite[\S8]{oss86} for a very readable introduction 
to Weierstrass data and the connection of minimal surfaces 
and holomorphic functions, 
\cite{Hoff90} for a basic discussion of the Weierstrass data of TPMS, 
\cite{Ross92} for identifying the Riemannian surface $\CM$ for the 
Schwarz P surface with $S^2\times S^2$ 
by stereographic projection, where $S^2$ is the unit sphere, and 
\cite{FW09} for further examples for construction of TPMS from 
Weierstra\ss\ data. 

Following \cite{KPS18}, we consider $\CM$ a double cover of $\C$, and, 
for $a\in(2,\infty)$, let 
\huga{\label{SR} 
R(z)=1/\sqrt{z^8+az^4+1}, 
}
where the Schwarz P surface with period cell $[-\pi,\pi)^3$ 
is obtained for $\vt=0$ and $a=14$.%
\footnote{For $\vt=\pi/2$ we obtain the Schwarz D family, and 
for $\vt\approx 0.9073$ Schoen's gyroid, as two further TPMS. 
Moreover, since these have the same Jacobians as Schwarz P, 
all bifurcation results from Schwarz P carry over to Schwarz D 
and the gyroid, but these appear to be much more difficult to treat 
in our numerical setting.}
See also \cite{Ga00} for the explicit computation of a fundamental 
patch of Schwarz P based on \reff{ew1} and \reff{SR} with $a=14$ and a small 
planar preimage $\subset\C$. 

In \cite{KPS18}, $a$ is taken as a bifurcation parameter 
along the Schwarz P  family % (and hence also Schwarz D and gyroid),  
with the {\em periods} for Schwarz P given by \cite[\S7.3]{KPS18} 
\hual{\label{pper}
&E=2\int_0^1 \frac{1-t^2}{\sqrt{t^8+at^4+1}}\dd t+
4\int_0^1\frac{\dd t}{\sqrt{16t^4-16t^2+2+a}}\quad \text{(periods in $x$ and $y$)}, \\
&F=8\int_0^1\frac t{\sqrt{t^8+at^4+1}}\dd t\quad \text{ (period in $z$)}, 
}
up to homotheties (uniform scaling in all directions). 
We have $\del=2\pi F/E$ for our $\del$, and evaluating 
$E,F$ numerically (or as elliptic integrals) and 
plotting $\del(a):=2\pi F/E$ as a function of $a$ we get the blue curve 
in Fig.\,\ref{tpsf1}(h), which corresponds to \cite[Fig.13]{KPS18}. 
In particular, $\del(a)$ has a maximum at $a=a_2\approx 28.778$, 
and $\del(a_2)=\del_f$ agrees with our fold position 
in Fig.\,\ref{tpsf1}(a). On the other hand, 
with suitable mesh adaptation 
the branch P1 continues to at least $\del=10$. 
Next, \cite{KPS18} based on \cite{EjSh18} 
gives a bifurcation from the P family at $a=a_1\approx 7.4028$, 
and again we find excellent agreement $\del(a_1)=\del_1$ with our 
BP at $\del_1$. 
}\eex\erem

\begin{figure}[ht] %H]
\centering 
\btab{lll}{{\sm (a)}&{\sm (b)}&{\sm (c)}\\ 
\hs{0mm}\ig[width=0.26\tew]{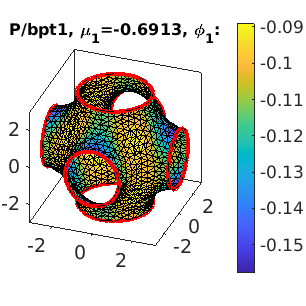}&
\hs{0mm}\ig[width=0.26\tew]{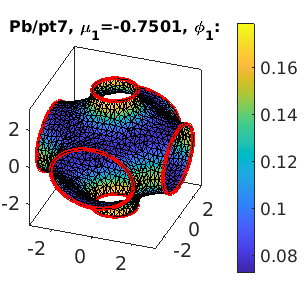}&
\hs{0mm}\ig[width=0.26\tew]{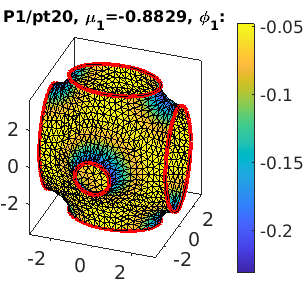}\\[-3mm]
\hs{0mm}\ig[width=0.26\tew]{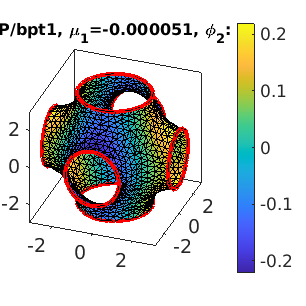}&
\hs{0mm}\ig[width=0.26\tew]{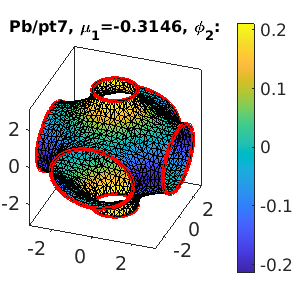}&
\hs{0mm}\ig[width=0.26\tew]{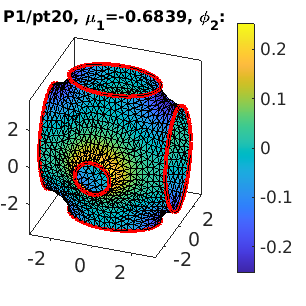}
}
\vs{-6mm}
\caption{\small Selected eigenvectors at points as indicated, 
cf.~Remark \ref{Prem2}b). 
Top: the most unstable direction, which does not change sign. 
Bottom: the second eigenvector; in (a) this approximately spans the kernel. 
\label{tpsf2}}
\end{figure}

\brem\label{TPrem2}{\rm 
a) The fact that P does not extend to ``large'' $\del$ (but folds back) has also 
been explained geometrically in \cite{Hoff90}, without 
computation of the fold position. 

b) The stability of Schwarz P (and hence also Schwarz D) on a 
minimal period cell and wrt {\em volume preserving} variations 
is shown in \cite{Ross92}. However, ``larger pieces'' of P, e.g., 
P on $[-\pi,\pi)^2\times [-2\pi,2\pi)$ are {\em always} unstable, even 
wrt volume preserving variations. See also \cite[\S8]{Bra96} 
for a useful discussion, and illustrations. 
Numerically, in Fig.\,\ref{tpsf1} 
we find: $\ind(X)=2$ except on the segment $\CS$ of P (and Pb) 
between the fold and the BP at $\del_1$, where $\ind(X)=1$. 
However, the most (and on $\CS$ only) unstable eigenvector has a sign, see Fig.\,\ref{tpsf2}, 
and hence the solutions on $\CS$ 
are stable wrt volume preserving variations. 
 }\eex\erem 

\subsubsection{CMC companions of Schwarz P}\label{SPCsec}
\begin{figure}[ht] %H]
\centering
\btab{ll}{\btab{l}{{\sm (a)}\\ 
\hs{0mm}\ig[width=0.27\tew]{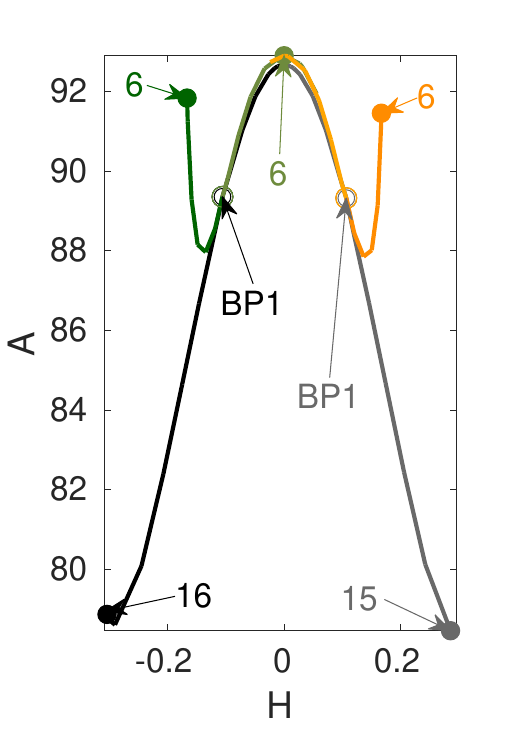}}&
\hs{-4mm}\rb{0mm}{\btab{l}{{\sm (b)}\hs{40mm}{\sm (c)}\\
\hs{0mm}\ig[width=0.28\tew]{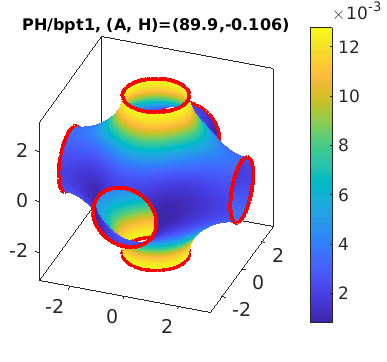}
\hs{0mm}\ig[width=0.24\tew]{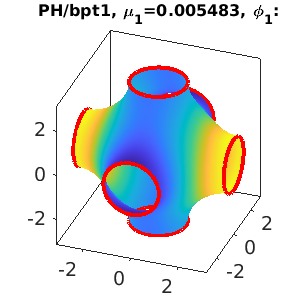}
\hs{-5mm}\ig[width=0.24\tew]{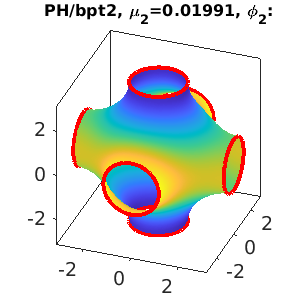}\\[-3mm]
{\sm (d)}\\%\hs{40mm}{\sm (d)}\\
\hs{0mm}\ig[width=0.24\tew]{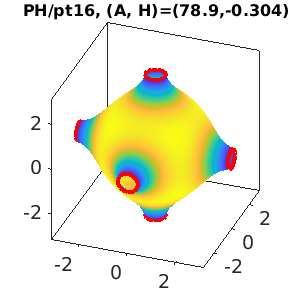}
\hs{-2mm}\ig[width=0.24\tew]{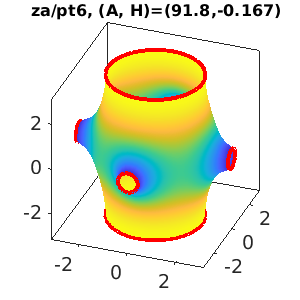}
\hs{-2mm}\ig[width=0.24\tew]{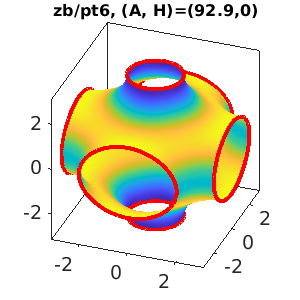}\\[-2mm]
%{\sm (e)}\\
\hs{0mm}\ig[width=0.24\tew]{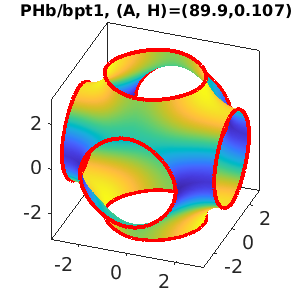}
\hs{-2mm}\ig[width=0.24\tew]{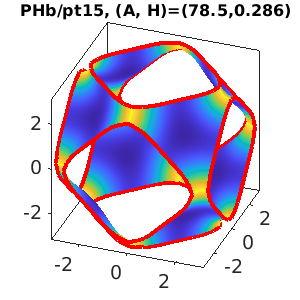}
\hs{-2mm}\ig[width=0.24\tew]{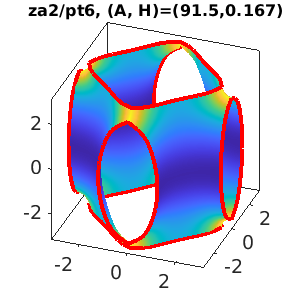}
}}
}
\vs{-4mm}
\caption{\small Some results from {\tt TPS/cmds2.m}. 
Continuation of Schwarz P in $H$ at fixed $\del=2\pi$. BD in (a): 
{\tt PH} (black), {\tt PHb} (grey), {\tt za} (dark green) and {\tt zb} 
(lighter green), and {\tt za2} and {\tt zb2} (orange), 
but altogether these are just two different branches. 
BP1 on {\tt PH} and approximate kernel vectors in (b,c), 
and further samples in (d). See text for details. 
\label{tpsf5}}
\end{figure}

In {\tt TPS/cmds2.m} we compute some CMC companions of the 
Schwarz P surface, see Fig.\,\ref{tpsf5}, where all we have to do is set 
{\tt ilam=[1,5,6,7]} as $H$ sits at position 1 in the  
parameter vector (and the translational Lagrange 
multipliers at {\tt [5,6,7]}). 
Continuing first to smaller $H$ 
(black branch {\tt PH}), 
$X$ (the volume enclosed by $X$ and the boundaries of the cube) 
``shrinks'' and we find a BP at $H\approx -0.1$. In the other direction 
(grey branch {\tt PHb}), 
$X$ (the volume enclosed by $X$) ``expands'', with a BP at $H\approx 0.1$. 
The continuation of both these branches fails at $H\approx -0.3$ and 
$H\approx 0.3$ (respectively), though they can be continued 
slightly further with careful mesh adaptation. 

Our main purpose here is 
to show how symmetry considerations and some tricks can 
help to avoid numerical pitfalls. 
By symmetry, the BP {\tt PH/bpt1} (and similarly {\tt PHb/bpt1}) 
must be double, although the smallest (in modulus) 
eigenvalues reported at {\tt PH/bpt1} 
are $\mu_1\approx 0.005$ and $\mu_2=0.02$.%
\footnote{Additionally, there is a simple negative eigenvalue $\mu_0\approx
-0.7$, 
and the next two eigenvalues are $\mu_{3,4}\approx 0.5$, i.e., $\mu_{1,2}$ are 
well separated from the rest of the spectrum.}  
See Fig.\,\ref{tpsf5}(b) for {\tt PH/bpt1}, and (c) 
for the (approximate) kernel vectors $\phi_1$, $\phi_2$. In fact, 
the plot in (b) (stronger correction along the $z$--axis) 
shows that at least the last step in the localization 
of  {\tt PH/bpt1} violated the $S_4$ symmetry of the (now fixed) cube, 
which explains the rather significant splitting of the in principle double 
eigenvalue $\mu_1=0$. Clearly, we expect $\phi_{1,2}$ to approximate 
two bifurcation directions, with $D_4$ symmetry along the $x$ axis ($\phi_1$) 
and $y$ axis ($\phi_2$). By symmetry we then must have at least one 
more bifurcating branch, with $D_4$ symmetry along the $z$ axis. To find this bifurcation 
direction, we can use {\tt qswibra} with 
numerical derivation and solution of the algebraic bifurcation 
equation (ABE) \cite[\S3.2.2]{p2pbook}. 
However, this is expensive and not always 
reliable. Here, the three bifurcation directions (oriented along $x$, 
along $y$, and along $z$) {\em are} returned, 
but we have to relax the tolerance {\tt isotol} for identifying 
solutions of the ABE as isolated. 
Alternatively, cf.~also Footnote \ref{qtrickfn}, 
we can use {\tt qswibra} with {\tt aux.besw=0} 
(bifurcation equation switch$=0$) to let {\tt qswibra} just compute 
and plot the (approximate) kernel $\phi_1, \phi_2$. This lets us 
guess to approximate the third direction as $\phi_3=0.2\phi_1+\phi_2$. 
This turns out to be sufficiently accurate and 
gives the transcritical branch(es) {\tt za} (dark green) and {\tt zb} 
(other direction, lighter green). 

On {\tt za}, the continuation fails after {\tt pt6}. {\tt zb/pt6} is 
at $H=0$ and corresponds to {\tt Pb/pt7} from Fig.\,\ref{tpsf1} 
%(which however in contrast to {\tt zb/pt6} does not sit exactly at $\del=2\pi$). 
Subsequently, {\tt zb} continues to {\tt PHb/bpt1}, and is indeed 
identical to the branch(es) {\tt za2} (and {\tt zb2}), transcritically bifurcating 
there. In particular, {\tt PHb/bpt1} is again double, and 
we can compute the three branches oriented along $x,y$ or $z$ as 
above (see {\tt cmds2.m}). {\tt zb2} (light orange) then continues 
back to {\tt PH/bpt1}, while {\tt zb2} fails after {\tt pt6} 
(last sample in (d)). The continuation failures of {\tt za} and {\tt za2} 
after {\tt pt6} are due to poor meshes as the different boundaries of $X$ come 
close to each other, like after {\tt PH/pt16} and {\tt PHb/pt15}, and it seems 
difficult to automatically adapt these meshes.

\subsection{Fourth order biomembranes, cylinders}\label{wsec}
The (dimensionless) Helfrich functional  \cite{H73} is 
\huga{\label{H1}
E=\int_X (H-c_0)^2+bK\dd S+\lam_1(A-A_0)+\lam_2(V-V_0), 
}
where $c_0\in\R$ is called spontaneous curvature parameter, 
$b\in\R$ is called 
 saddle--splay modulus,  $A$ and $V$ 
are the area and the volume of $X$, and $\lam_1$ and $\lam_2$  
are Lagrange multipliers for area and volume constraints 
$A(X)=A_0$ and $V(X)=V_0$.  For closed $X$, 
$\lam_1$ corresponds to 
a surface tension, and $\lam_2$ to a pressure difference between 
outside and inside, and the Euler-Lagrange equation is 
\huga{\label{H2}
\Delta H+2H(H^2-K)+2c_0K-2c_0^2H-2\lam_1H-\lam_2=0, 
}
together with the area and volume constraints. In this case, 
the term b$\int_X K\dd S$ in \reff{H1} can be dropped due to 
the Gauß--Bonnet theorem, cf.~Footnote \ref{GBfoot}, as 
$\int_X K\dd S=2\pi\chi(X)$ is a topological constant.

If $X$ is not closed, then usually the constraints $A=A_0$ and $V=V_0$ 
are dropped, and $\lam_1$ in \reff{H1} denotes an external surface 
tension parameter, and wlog $A_0=0$, and $\lam_2$ is fixed to $0$. 
If in 
\huga{\label{GB1} 
\int_X K\dd S=2\pi\chi(X)-\int_{\pa X}\kap_g \dd s 
}
we assume $\ga=\pa X$ to be parameterized by arclength, 
then the geodesic curvature $\kap_g$ 
is the projection of the curvature vector $\ga''(\vx)$ onto 
the tangent plane $T_{\vx}(X)$, see, e.g., \cite[\S4.3]{Tapp16}. 
%If $X$ has no boundary, then $\int K\dd S$ is constant. 
If as before we restrict to normal variations 
$\psi=uN$, which moreover fix the boundary, i.e.,
\huga{\label{hdbc}
\text{ $u|_{\pa X}=0$, }
}
then 
\hualst{
\pa_\psi E&= 
\int_X (\Delta H+2H(H^2-K)+2c_0K-2Hc_0^2-2\lam_1 H)u\dd S
+\int_{\pa X} (H-c_0+b\kap_n)\pa_n u \dd s,}
where $\kap_n=\spr{\ga'',N}$ is the normal curvature of $\ga=\pa X$, 
i.e., the projection of the curvature vector onto the normal plane, 
see, e.g.,  \cite{PP22} and the references therein. 
Thus we again obtain \reff{H2} (with $\lam_2=0$), and additionally 
to \reff{hdbc} we can consider either of 
\hual{
\pa_n u&=0 \text{ on $\pa X$ (clamped BCs, or Neumann BCs),}\label{cBCs}\\
H-c_0+b\kap_n&=0 \text{  on $\pa X$ (stress free BCs).}\label{natBCs}
}
In the case of \reff{cBCs} we have %$N(X)\perp X$ on $\pa X$ such that 
$\int_{\pa X}\kap_g \dd s=0$ in \reff{GB1}, and hence $\int bK\dd S$ again 
becomes constant and can be dropped from \reff{H1}. 

Following \cite{DDG21}, we first % in the demo {\tt biocyl} 
consider ``Helfrich cylinders'', i.e., \reff{H2} 
with cylindrical topology and BCs  \reff{hdbc} and \reff{cBCs}. 
In \S\ref{hcap} we then consider ``Helfrich caps'', i.e., 
disk type solutions with BCs  \reff{hdbc} and \reff{natBCs}. 

\brem\label{vesrem}{\rm a) The original 
motivation of \reff{H1} are the shapes of closed vesicles 
 with a lipid bilayer membrane, 
in particular red blood cells (RBCs).   This 
motivated much work, e.g., 
\cite{SBL90, sei97, NT03, VDM08, OYT14}, aiming to understand the various shapes 
of RBCs, 
mostly in the axisymmetric case. 
See also \cite{LWM08,KIPM20} for further biological and mechanical 
background, and \cite{JQJQY98} for non--axisymmetric shapes. 
Applying our algorithms to RBCs (without a priori enforcing any symmetry) 
we recover many of 
the results from the above references. However, the bifurcation 
diagrams quickly become very complicated, and therefore our results will 
be presented in detail elsewhere. See also 
\cite{FVGK22} for a 1D version with an extremely rich bifurcation structure. 

b) For $X$ not closed it is an open problem for what parameters, and 
boundaries and BCs, the minimization of \reff{H1} is a well--posed problem. 
In \cite{Ni93}, the following conditions on $c_0,b$ and $\lam_1$ are posed 
for $E$ with $\lam_2=0$ 
to be definite in the sense that $E\ge C_0$ for some $C_0>-\infty$ 
for all connected orientable surfaces $X$ of regularity $C^2$ with 
or without boundary: 
\huga{\label{ndefc} 
\text{(i) }\lam_1\ge 0,\qquad \text{(ii) }-1\le b\le 0, 
\qquad\text{and (iii) }-bc_0^2\le \lam_1(1+b).
}
This proceeds as in \cite{Ni91} by scaling properties of $E$ for various 
surfaces 
composed of planes (of area $A$), cylinders (of lengths $l$ and radius $r_c$), 
and (hemi)spheres (of radius $r_S$), and considering the asymptotics of $E$ as 
$A,l\to\infty$ and/or $r_c,r_s\to 0$. For instance, 
the condition \reff{ndefc}(i) arises most naturally by considering 
$X$ to contain a plane with $A\to\infty$, which for $\lam_1<0$ 
gives $E\to -\infty$. 

On the other hand, in the physics literature no restrictions 
on $c_0, b\in\R$ are 
given, and in a given problem a fixed $\pa X$ and the BCs 
\reff{cBCs} or \reff{natBCs} may make $E$ definite for much larger 
ranges than given in \reff{ndefc}. In our experiments below we do take 
parameters to rather extreme values, e.g., $\lam_1<0$ in Fig.\,\ref{hf2} 
in the sense of a ``detour'', and $b=-4$ in Fig.\,\ref{hcaf2}, where we find 
interesting solutions, which can then again 
be continued to moderate parameter regimes. 
For closed $X$, we additionally 
in general have $\lam_2\ne 0$, and, moreover, what then matters is 
the {\em reduced} volume defined by dividing $V$ by the volume of 
the sphere with the same area. 
%Thus, altogether \reff{ndefc} can be violated for a specific problem. 

c) For $\pa X\ne \emptyset$, and in particular for the cases of cylinders and 
caps considered below, we are not aware 
of analytic bifurcation results, although \cite{PP22} presents some 
results for caps in a slightly different setting. 
}\eex\erem 

%\subsection{Biomembrane cylinders}\label{hcyl}
%\subsubsection{Cylinders, general setup}\label{hcyl}
\def\dhome{geomtut/biocyl} 
In the demo {\tt biocyl} we consider a cylindrical 
topology along the $x$--axis with BCs \reff{cBCs}. 
The  equation and BCs thus are 
\begin{subequations}\label{HC}
\hual{
&\Delta H+2H(H^2-K)+2c_0K-2c_0^2H-2\lam_1H=0,\label{Hca} \\
&\text{$X_2^2+X_3^2|_{X_1=\pm 1}=\al^2$,\quad and\quad $\pa_x X_2|_{X_1=\pm 1}
=\pa_x X_3|_{X_1=\pm 1}=0$.} \label{cBC1} 
}
\end{subequations}
For $c_0=0$ we have the explicit family 
\def\Xc{X_{\text{cyl},\al}}
\huga{\label{hcyl1}
\Xc=(x,\al\cos\phi,\al\sin\phi), \ \  x\in[-1,1],\ 
\phi\in [0,2\pi), \text{ with }\lam_1=\frac 1 {4\al^2},
} 
of Helfrich cylinders, and \cite{DDG21} proves various existence 
results for axisymmetric solutions near \reff{hcyl1} 
and in other regimes in the $\al$--$\lam_1$ plane. See also 
\cite{VDM08} for other shapes  with cylindrical topology, which 
fit into our setting by prescribing other contours at $x=\pm 1$ 
instead of the circles of radius $\al$ in \reff{cBC1}. 

%We continue \reff{hcyl1} in $c_0$ and find various surfaces of revolution, and bifurcations to non--axisymmetric branches, see Fig.\,\ref{hf3a} and  Fig.\,\ref{hf3b}. 
The basic setup again consists in setting $X=X_0+uN_0$ (with here $N$ the 
inner normal), and then 
writing \reff{H2} as a system of two second order equations for 
$(u_1,u_2)=(u,H)$, namely 
\begin{subequations}\label{HCd}
\hual{
L u_2+Mf(u_2)+\srot q(X)&=0,\\
Mu_2-H&=0. 
}
\end{subequations}
As before, $L$ is the cotangent Laplacian, $M$ the (Voronoi) mass matrix, 
and $f(u_2)=2u_2(u_2^2-K)-2\lam_1u_2+2c_0K-2c_0^2 u_2$. The mean curvature  
$H=H(u_1)$ is computed as $H=\frac 1 2 \spr{LX,N}$,  
the Gaussian curvature $K=K(u_1)$ is obtained from {\tt discrete\_curvatures}, 
cf.~\reff{dKG1}, 
and $\srot q(u)$ implements the rotational PC with $\srot$ an  
active parameter for branches of non-axisymmetric 
solutions. %, while initially $\srot=0$.  
The reason for the reformulation of \reff{HC} as two second order equations 
\reff{HCd} is that this way we can easily 
implement the {\em two} BCs \reff{cBC1}, 
see {\tt sG.m}. 

The external parameters for \reff{Hca} are $(\al,\lam_1,c_0)$, 
and continuation in any of these yields interesting results. 
We first continue in $c_0$ at fixed $\al=1$, since we believe this is 
the numerically most challlenging case due to the exact solution \reff{hem2}. 
Subsequently we continue in $\al$, and again in $c_0$ at different $\al$, 
and finally in $\lam_1$, 
always starting from the Helfrich cylinder \reff{hcyl1}.  

\taskip
\begin{table}[ht]\caption{
{\sm Scripts in {\tt pde2path/demos/geomtut/biocyl}.} 
\label{hcyltab}}
\centering\vs{-2mm}
{\small 
\begin{tabular}{p{0.1\tew}|p{0.93\tew}} 
{\tt cmds1.m}&Continuation of $\Xc$ in $c_0$, $(\al,\lam_1)=(1,1/4)$, Figs.~\ref{hf3a} and 
\ref{hf3b}\\
{\tt cmds2.m}&Continuation of $\Xc$ in $\al$, $(c_0,\lam_1)=(1/4,1)$, Fig.\,\ref{hf22}.\\
{\tt cmds3.m}&Continuation of ``big'' and ``small''  cylinders $\Xc$ 
in $c_0$, i.e., $\al=1.4$ and $\al=0.6$,  Fig.\,\ref{hf20}. \\
{\tt cmds4.m}&Continuation of $\Xc$ like solutions in $\lam_1$, 
$(\al,c_0)=(1,0.7)$, Fig.\,\ref{hf2}.\\
{\tt cmds4b.m}&Like {\tt cmds4.m} but with $(\al,c_0)=(1,0)$, 
Figs.~\ref{hf2a} and \ref{hf2g}. 
\end{tabular}
}
\end{table}
\teskip

Table \ref{hcyltab} lists the command scripts in {\tt biocyl/}. 
Besides {\tt cylinit.m}, 
{\tt sG.m}, and {\tt getM.m} (see Remark \ref{hstabrem}), 
we then additionally 
have the PC {\tt qfrot.m} and its derivative {\tt qfrotder.m}, 
and the scripts 
{\tt bdmov1.m} and {\tt bdmov2.m} to produce movies of the BDs in 
Fig.\,\ref{hf2}. In {\tt hcylbra.m} we put $A,V,E,\lam_1 A$ and 
the mesh--quality data on the branch.
%, where $E_A=\lam_1 A$ is used in Fig.\,\ref{hf2d}. Finally, 
In {\tt cylinit.m} we set 
{\tt p.sw.Xcont=2}, such that the colors in solution plots 
 mainly give some visual structure to $X$, but 
do not in general give the continuation direction $uN$, cf.~the remarks after 
\reff{Xsw}. 
%For instance, the colors in {\tt l3qr/pt30} and {\tt l4qr/pt34} of Fig.\,\ref{hf2} 
%indicate that the {\em last} Newton step was mainly for phase correction. 
The PC multiplier $|\srot|<10^{-4}$ on all the branches, and 
in all cases the final $u$ plotted is of order $10^{-6}$ or smaller.

\brem\label{hstabrem}{\rm Since $N$ is the inner normal, 
the stability for \reff{HC} refers to the Helfrich flow 
\huga{\label{Hflow}
\dot X=-[\Delta H+2H(H^2-K)+2c_0K-2c_0^2H-2\lam_1H]N, 
} 
with BCs \reff{cBC1}. This is encoded in 
the dynamical mass matrix $\ds\CM=\bpm M&0\\0&0\epm$ in {\tt getM}, where 
$M$ is the 1--component (Voronoi) mass matrix. 
}
\eex
\erem 

%\input{hcylf10a}
%\paragraph{Continuation in $c_0$.}
\subsubsection{Continuation in the spontaneous curvature}
\begin{figure}[ht]
\centering 
\btab{ll}{{\sm (a)}&{\sm (b)}\\[-0mm]
\hs{-1mm}\rb{0mm}{\btab{l}{\ig[height=38mm]{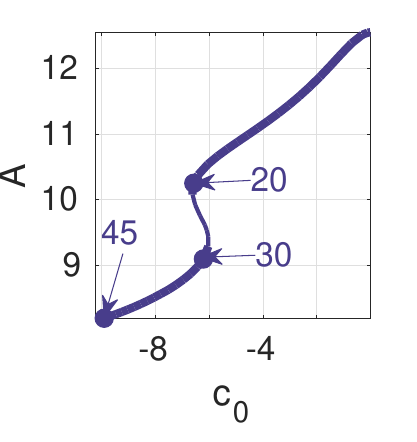}\\
\ig[height=39mm]{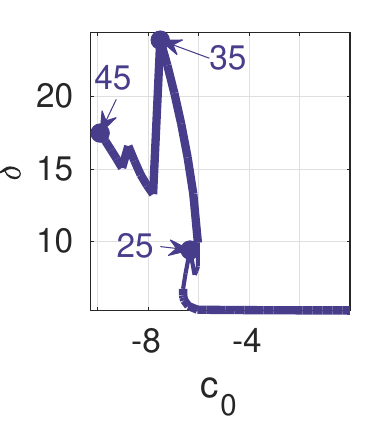}}}
&\hs{-10mm}\btab{l}{\ig[height=44mm]{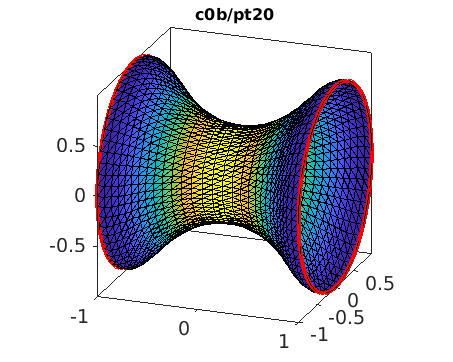}
\hs{-8mm}\ig[height=44mm]{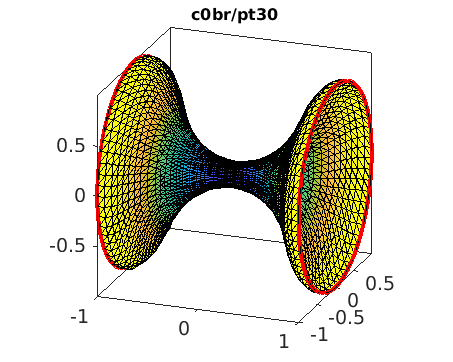}
\hs{-8mm}\ig[height=44mm]{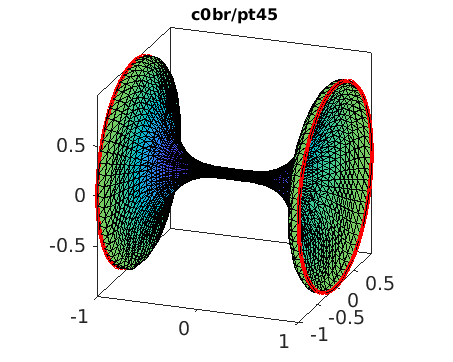}\\[-2mm]
\hs{0mm}\ig[height=42mm]{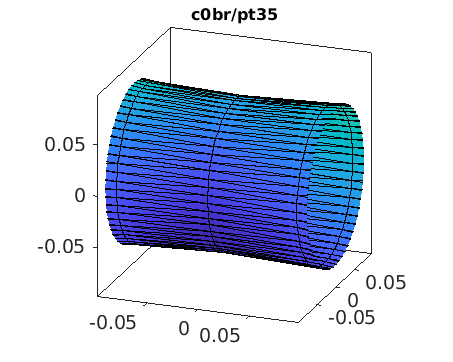}
\hs{-5mm}\ig[height=42mm]{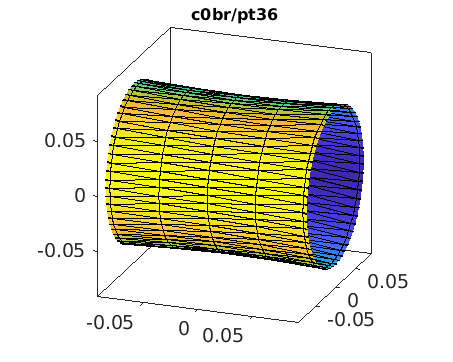}
\hs{-5mm}\ig[height=42mm]{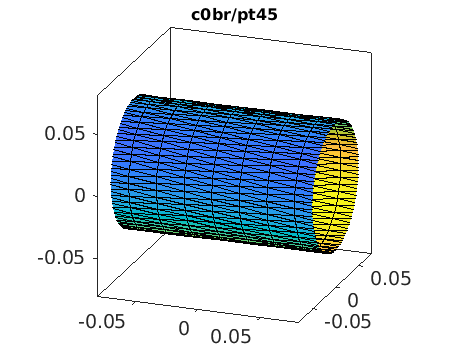}}}
\vs{-4mm}
\caption{{\small Results for \reff{HCd} from {\tt biocyl/cmds1.m}. 
$(\al,\eps)=(1,1/4)$ fixed, continuation to $c_0<0$, with a bistability region near $c_0=-6.25$. Initially uniform mesh of $n_p=1770$ points, and the 
2nd plot in (a) 
indicates that with suitable refinement (at {\tt pt25} and {\tt pt35}, 
to $n_p=3561$) the meshes stay good.   
Samples in (b), with zooms of the necks to illustrate the meshes.  
\label{hf3a}}}
\end{figure}

In {\tt cmds1.m} we fix $\al=1$, and start with the Helfrich cylinder 
\reff{hcyl1}, hence $\lam_1=1/4$. We choose a rather coarse 
uniform initial mesh of {\tt np}=1770 points, and first continue to $c_0<0$, 
yielding the branch {\tt c0b} in Fig.\,\ref{hf3a}. 
The solutions contract in the middle and via two folds produce a 
bistability region around $c_0=-6.25$, but otherwise no 
bifurcations. As the neck thins, we refine the mesh based on 
{\tt e2rsshape1}, cf.~\reff{delmdef}, with {\tt p.sw.rlong=1} 
and combined with {\tt retrigX} to 
avoid obtuse triangles, cf.~Remark \ref{Xcontrem2}. 
The 2nd plot in (a) shows the distortion $\delm$, with refinements at 
{\tt pt25}, {\tt pt35} 
and {\tt pt40}, and the second row in (b) shows rather strong zooms 
into the necks, illustrating the refinement step from {\tt pt35} to 
{\tt pt36}, and the reasonable mesh in the neck at point {\tt pt45}. 

\begin{figure}[ht]
\centering 
\btab{ll}{{\sm (a)}&{\sm (b)}\\[-0mm]
\hs{-2mm}\ig[height=52mm]{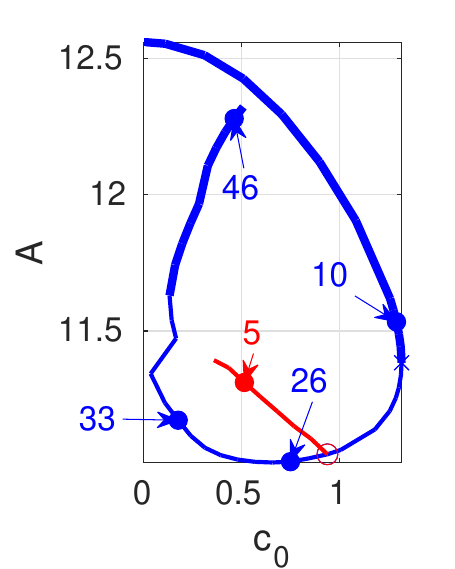}\hs{-0mm}\ig[height=52mm]{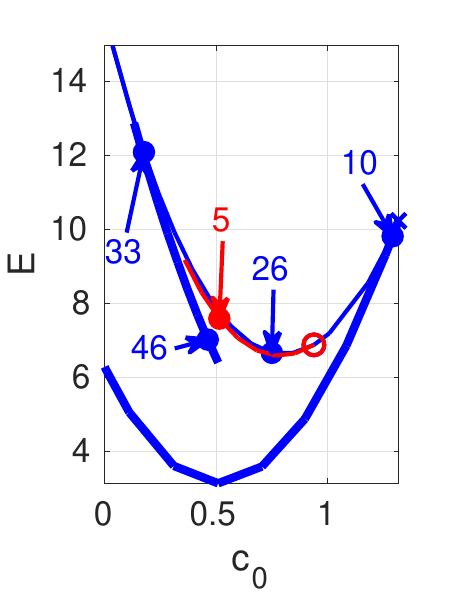}
\hs{-0mm}\ig[height=52mm]{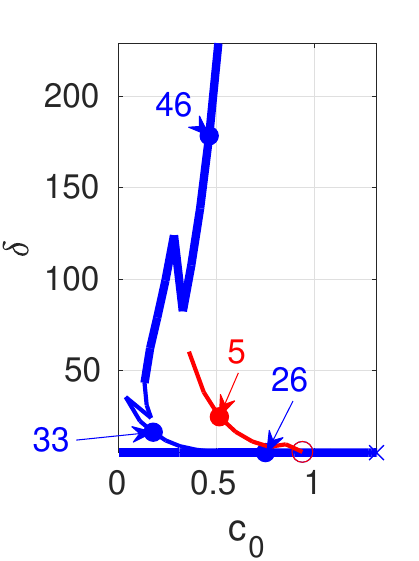}
&\hs{-6mm}\ig[height=52mm]{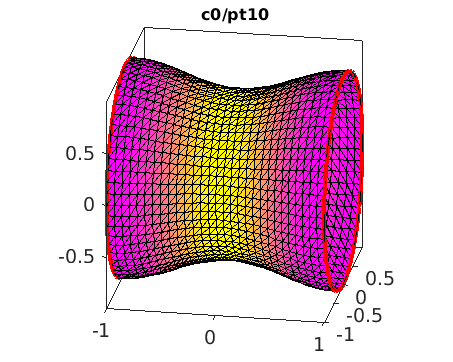}%\hs{-6mm}\ig[height=52mm]{h8/26c}
%\hs{-8mm}\ig[height=42mm]{h8/026}
}
\\[-3mm]
\btab{ll}{{\sm (c)}&{\sm (d)}\\
\hs{-2mm}\ig[height=52mm]{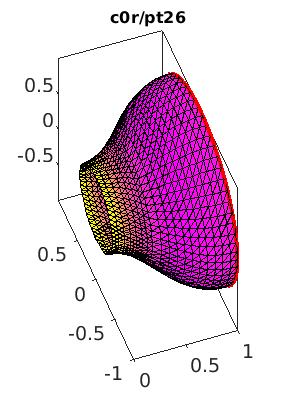}\hs{-4mm}\ig[height=50mm]{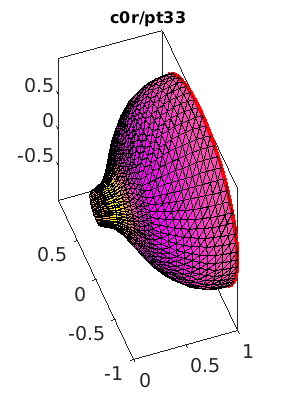}
\hs{-4mm}\ig[height=50mm]{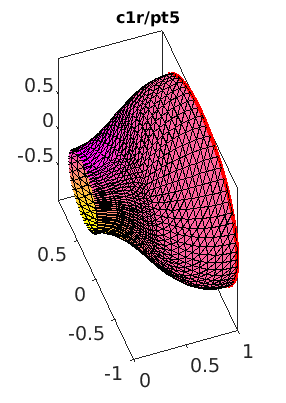}
\hs{-4mm}\ig[height=50mm]{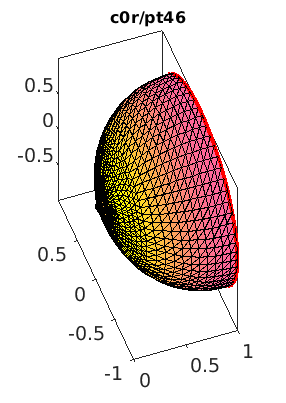}
&\hs{-2mm}\ig[height=50mm]{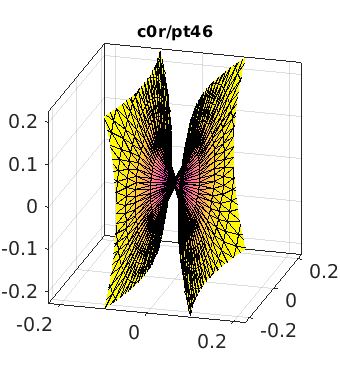}
}
\vs{-4mm}
\caption{{\small Further results from {\tt biocyl/cmds1.m},  
$(\al,\eps)=(1,1/4)$, continuation of Helfrich cylinder 
to $c_0>0$ (blue branch {\tt c0}, {\tt c0r} after mesh refinement), 
with bifurcation to {\tt c1r} (red). (a) BDs of area $A$, energy $E$, and 
mesh distortion $\delm$. (b) A full sample at {\tt c0/pt10},  
different colormap (yellow$>$pink) 
for better visibility of mesh. (c) further samples, cut 
at $x=0$. (d) zoom into solution close to $\Xhem$. 
On {\tt c0}, results after {\tt pt33} 
become somewhat mesh--dependent, but our experiments suggest that the 
branch connects to $\Xhem$ at $c_0=1/2$ (with {\tt c0r/pt46} 
shortly before that connection). 
\label{hf3b}}}
\end{figure}

In Fig.\,\ref{hf3b} we continue to $c_0>0$. This is initially easy, 
but after a first fold near $c_0\approx 1.3$ and then decreasing $c_0$ 
below $c_0\approx 0.15$ it becomes difficult to maintain mesh distortion 
$\del<100$. We give the full result of our continuation 
but remark that the behavior on {\tt c0} ({\tt c0r} after refinement, 
to $n_p=3130$ at {\tt pt46})   
after {\tt pt30}, say, becomes mesh dependent. 
First, however, at $c_0\approx 0.94$, we find a BP to a non--axisymmetric 
branch ({\tt c1}, red), 3rd sample in (c). We then get a second 
fold at $c_0\approx 0.08\pm \eta$, with $\eta\in(-0.02,0.02)$ 
dependent on the mesh. Relatedly, the further continuation becomes 
somewhat non--smooth and quantitatively depends on the mesh, 
but qualitatively we get similar behavior for different meshes, 
namely: The neck becomes thin and short, and the 
solutions seem to approximate a solution 
\huga{\label{hem2}
\text{$\Xhem$ consisting of two 
hemispheres with radius 1, centered at $(x,y,z)=(\pm 1,0,0)$,}
} 
and hence touching at $(0,0,0)$, see {\tt c0r/pt46}, at $c_0=0.47$. 
$\Xhem$ is an exact solution of \reff{HC} with 
$\lam_1=1/4$ at $c_0=1/2$, and given our various experiments with 
mesh refinement we believe that the branch {\tt c0r} connects 
to $\Xhem$. However, the third plot in (a), and the samples in 
(c) (at $c_0\approx 0.47$, shortly before the supposed connection 
to $\Xhem$) show how the mesh quality seriously degrades as we approach 
the supposed connection.

Also, we have, e.g., $\ind(X)=3$ at 
{\tt pt33} (with one unstable direction from the fold at $c_0\approx 1.3$, 
and two from the double BP at $c_0\approx 0.94$), while the 
``almost--$\Xhem$''--solutions near {\tt pt46} have $\ind(X)=0$, i.e., are stable. 
Hence, since $\ind(X)$ should decrease by one at the (supposed) fold between 
{\tt pt33} and {\tt pt46}, we expect another double BP between 
{\tt pt33} and {\tt pt46}. However, the behavior of 
the branch is quantitatively mesh--dependent also near the left fold, 
and while the changes in $\ind(X)$ are detected, the bisection 
loops to localize the fold and/or the BPs do not converge. 
Thus, altogether the behavior after {\tt pt33} {\em is conjectured}. 
On the other hand, for $c_0\in(0.3,0.5)$, say, and in particular 
near {\tt pt46}, the solutions like {\tt pt46} stay stable under 
mesh refinement. In any case, a pinch--off (topological change) 
which should occur after {\tt pt46} as we approach $\Xhem$ 
cannot be resolved with our methods, 
but probably requires the use of some phase field method, see, e.g., 
\cite{DLW06}.%
\footnote{For Neumann BCs \reff{cBCs}, hemispheres 
are highly degenerate in the following sense: Setting up \reff{HC} over 
a {\em single} circle at wlog at $X_1=1$, 
and wlog with $(\al,\lam_1,c_0)=(1,1/4,1/2)$ the single unit hemisphere 
is again an exact solution, but we are not able to continue this in 
$\lam_1$ or $c_0$, i.e., it seems isolated. This is different for BCs \reff{natBCs}, see \S\ref{hcap}.}

The results for {\tt c0r} up to {\tt pt33} are 
mesh independent, including the bifurcation of the 
non-axisymmetric branch {\tt c1} at {\tt c0/bpt1}. 
The mesh then also degrades on the non--axisymmetric branch ({\tt c1}, 
and in a similar way like on the branch {\tt c0} 
but at larger $c_0$), and since next we focus on non--axisymmetric branches in 
a slightly different setting we do not further pursue {\tt c1} here.

\subsubsection{Intermezzo: Other radii}
\begin{figure}[h]
\centering 
\btab{l}{
\hs{-0mm}\ig[height=42mm]{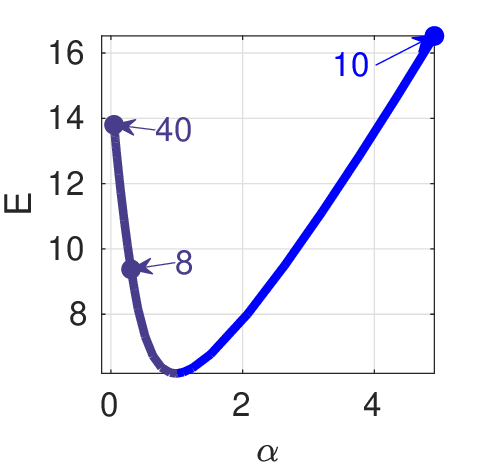}
\hs{-0mm}\ig[height=44mm]{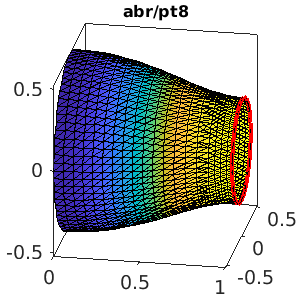}
\hs{-0mm}\ig[height=44mm]{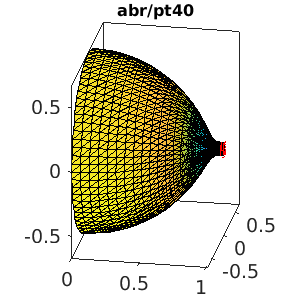}
\hs{-0mm}\ig[height=44mm]{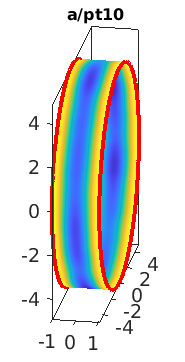}

}
\vs{-4mm}
\caption{{\small Results for \reff{HCd} from {\tt cmds2.m}, continuation 
 in $\al$.
\label{hf22}}}
\end{figure}

%Continuation in $\al$, and continuation at different $\al$}
In {\tt cmds2.m} we continue $\Xc$ in $\al$, with fixed $(c_0,\lam_1)=(0,1/4)$, 
see Fig.\,\ref{hf22}. This %behaves as expected, and 
shows no bifurcations. 
Numerically, the case $\al\searrow 0$ is challenging due to boundary 
layers developing at $x=\pm 1$, see \cite{HCG10}. Alternating continuation 
and mesh adaptation  
based on {\tt e2rsshape1} we can reliably continue to $\al=0.05$ and slightly further. 

\begin{figure}[H]
\centering 
\btab{l}{{\sm (a)}\\
\hs{-0mm}\ig[height=40mm]{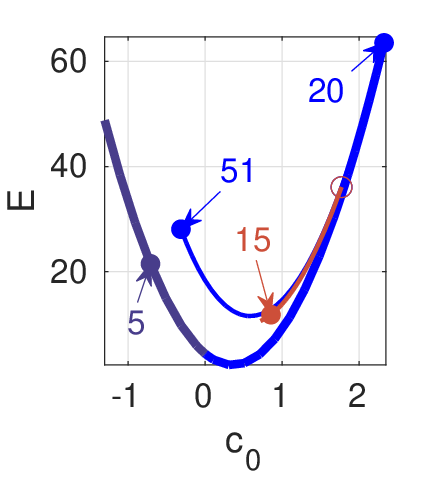}
\hs{3mm}\rb{3mm}{\ig[height=37mm]{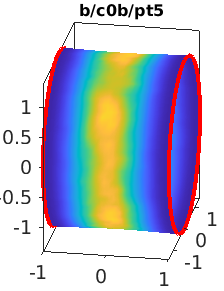}
\hs{1mm}\ig[height=37mm]{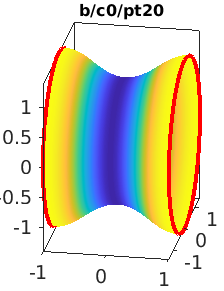}
\hs{1mm}\ig[height=37mm]{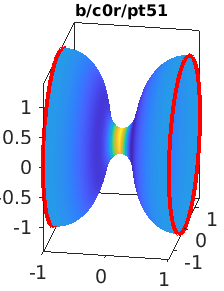}
\hs{1mm}\ig[height=37mm]{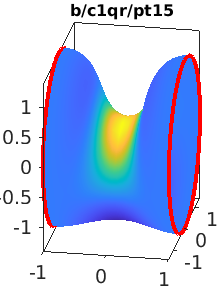}}\\[-2mm]
{\sm (b)}\hs{130mm}{\sm (c)}
\\[-0mm]
\hs{-2mm}\rb{0mm}{\ig[height=35mm]{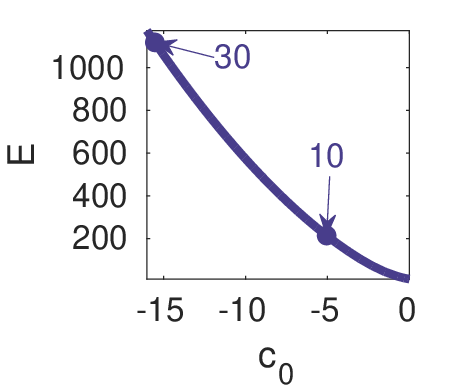}
%\ig[height=40mm]{hcp2/bdsEa}
}
\hs{-0mm}\ig[height=32mm]{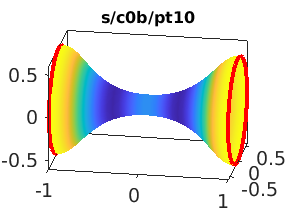}
\hs{-0mm}\ig[height=32mm]{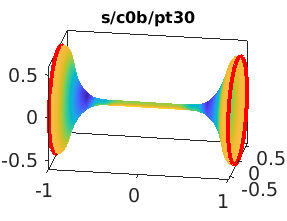}
\ig[height=35mm]{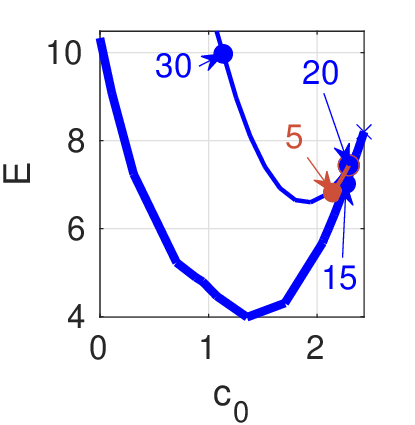}
\\[-3mm]
{\sm (d)}\\
%\ig[height=40mm]{hcp2/bdsEb}
\rb{5mm}{\hs{-0mm}\ig[height=31mm]{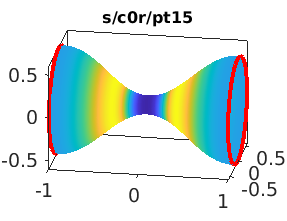}
\hs{-0mm}\ig[height=31mm]{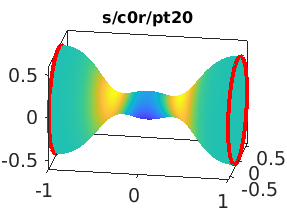}
\hs{-0mm}\ig[height=31mm]{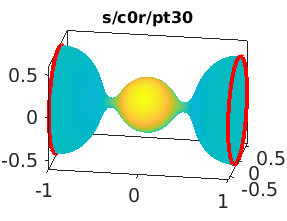}
\hs{-0mm}\ig[height=31mm]{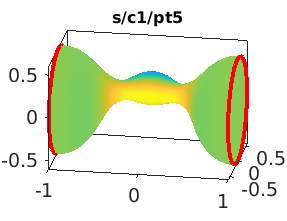}}
}
\vs{-6mm}
\caption{{\small % Results for \reff{HCd} from 
{\tt cmds3.m},  
continuation in $c_0$ for 
$(\al,\lam_1)=(1.4, 1/5.6)$ in (a), 
and $(\al,\lam_1)=(0.6,2.4)$ in (b--d). 
% For the ``big'' cylinder in (a), 
% the continuation to negative $c_0$ proceeds as in Fig.~\ref{hf3a} but 
% without a fold. Also the continuation to positive $c_0$ (branch {\tt c0}(r), blue) is qualitatively similar, but there seems to be no second fold, and 
% the branch continues with decreasing $c_0$ also after {\tt pt51}, 
% again needing strong mesh refinement at the neck, and similar in 
% {\tt c1qr/pt15}. For the ``small'' cylinder in (b--d), the continuation to 
% $c_0<0$ (b) yields long and thin necks, as expected. The continuation 
% to positive $c_0$ now gives a ``pearling'' behavior after the fold 
% near $c_0\approx 2.3$, also visible on the non-axi branch {\tt c1}.   
\label{hf20}}}
\end{figure}

In {\tt cmds3.m} and Fig.\,\ref{hf20} we repeat the continuation in $c_0$ 
from Fig.\,\ref{hf3a} and Fig.\,\ref{hf3b} for different starting cylinders 
$\Xc$, namely $\al=1.4$ in Fig.\,\ref{hf20}(a) and $\al=0.6$ in (b--d). 
The main differences to the case $\al=1$ are as follows: For both, 
$\al=1.4$ and $\al=0.6$, continuing to negative $c_0$ we no longer 
have the two folds and associated bistability range as in Fig.\,\ref{hf3a}, 
and the branches seem to extend to arbitrary large negative $c_0$. 
We refrain from plotting this for $\al=1.4$ in (a), and for 
$\al=0.6$ in (b) we mainly remark that large negative $c_0$ yields 
very long and thin necks, which can be handled with adaptive mesh refinement as 
in Fig.\,\ref{hf3a}. %, with large $E$. %, but never pinch--off. 

Continuing to positive $c_0$, for $\al=1.4$ in (a), the main difference 
to Fig.\,\ref{hf3b} is that the branch {\tt c0r} only has one fold 
(at $c_0\approx 2.3$) and then continues to negative $c_0$, see 
{\tt b/c0r/pt51}. The difference to $\al=1$ is that a solution like 
$\Xhem$ no longer exists, i.e., two hemispheres of radius $\al=1.4$ 
(or any $\al\ne 1$) cannot be near a (genus 1) cylinder type solution. 
Additionally, the larger radius makes the continuation of 
the red non--axisymmetric branch {\tt b/c1qr} easier wrt to mesh handling, 
see the last sample in (a). 
For $\al=0.6$ in (c,d), the smaller radius gives a ``pearling'' behavior 
after the fold on the blue branch, and it seems likely that the branch 
{\tt s/c0r} continues to a solution of type hemisphere--sphere--hemisphere, 
which would give similar problems as discussed for $\Xhem$ in Fig.\,\ref{hf3b}. 
Here we stop the continuation after {\tt pt30} (third sample in (d)), 
as further continuation requires excessive mesh adaptation and the solutions 
are unstable anyway. The pearling shape is also inherited by the 
non-axisymmetric branch {\tt s/c1} (last sample in (d)).

%\input{hcylf21}

%\clearpage

%\paragraph{Continuation in $\lam_1$.} 
\subsubsection{Continuation in surface tension} 
\begin{figure}[h]
\centering 
\btab{ll}{{\sm (a)}&{\sm (b)}\\[-0mm]
\hs{-4mm}\rb{0mm}{\ig[height=40mm]{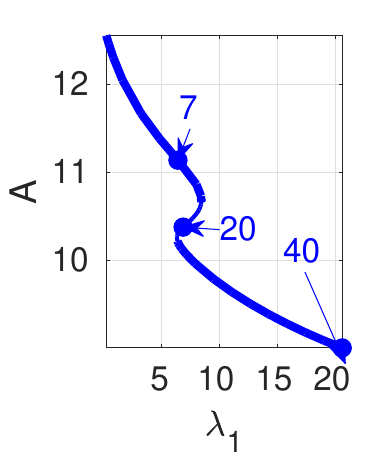}
\ig[height=40mm]{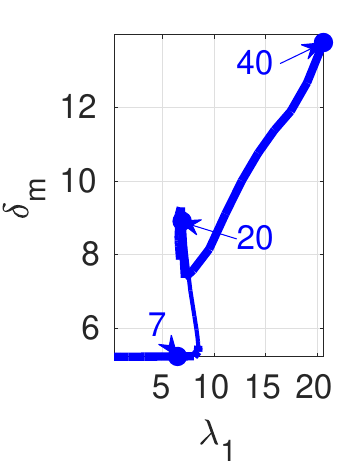}}
&\hs{-0mm}\ig[height=42mm]{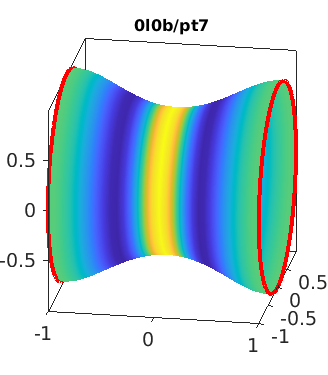}
\hs{-0mm}\ig[height=42mm]{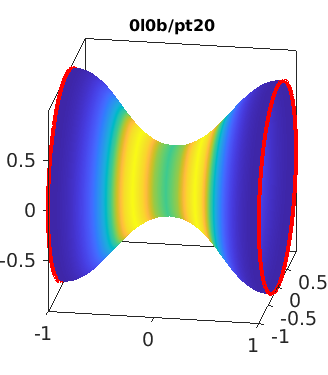}
\hs{-0mm}\ig[height=42mm]{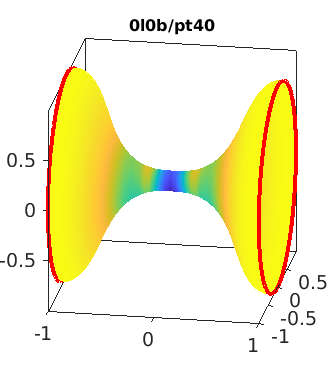}
}
\vs{-4mm}
\caption{{\small %Results for \reff{HCd} from 
{\tt cmds4b.m}, 
continuation of $\Xc$ in $\lam_1$, 
$(c_0,\al)=(0,1)$. 
%Refinement strategy as in Fig.~\ref{hf3a} to resolve the thin neck. 
\label{hf2a}}}
\end{figure}

%Almost all branches (solutions) from Figs.~\ref{hf3a}--\ref{hf20} are axisymmetric, and hence could be 
In {\tt cmds4.m} we return to fixed $\al=1$ and 
continue in $\lam_1$, starting at $\lam_1=0.25$ and 
$c_0=0.7$ corresponding to Fig~\ref{hf3b}, and in {\tt cmds4b.m} 
we repeat this for $c_0=0$, starting from the genuine Helfrich 
cylinder $\Xc$. 
In both cases, for increasing $\lam_1$ the solutions initially 
only slightly change shape, but at larger $\lam_1$ 
(near $\lam_1=8$ for $c_0=0$, and near $\lam_1=4$ 
for $c_0=0.7$) we get an S--shaped bistability region due to 
two consecutive folds, similar to the case of $c_0\approx -6$  
in Fig.\,\ref{hf3a} (with fixed $\lam_1=1/4$).%
\footnote{Such folds were already found in \cite[\S6.3]{doem17} 
for $c_0=0$ in the 1D axisymmetric setting.}
However, we find no bifurcations to non--axissymmetric branches. 
The case of decreasing $\lam_1$ is more interesting, and 
in Fig.\,\ref{hf2a} we only show the basic result for $c_0=0$ and 
large $\lam_1>0$ for completeness.

\begin{figure}[h]
\centering
\btab{l}{{\sm (a)}\hs{50mm}{\sm (b)}\hs{40mm}{\sm (c)}\hs{40mm}{\sm (d)}\\
\ig[height=48mm]{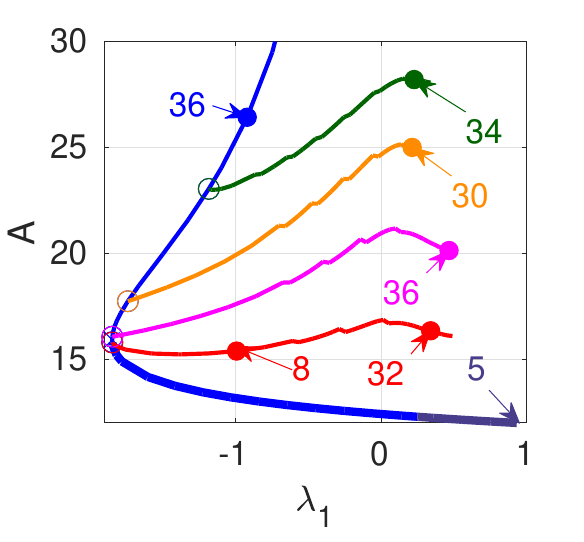}\ig[height=46mm]{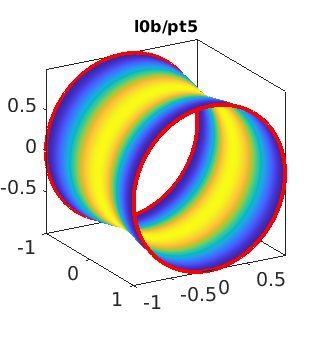}
\hs{-1mm}\ig[height=48mm]{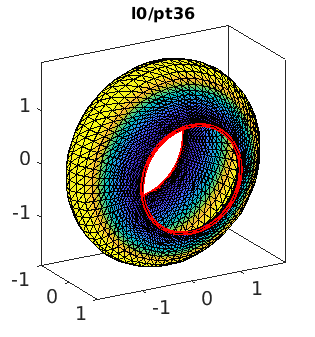}
\hs{-4mm}\ig[height=48mm]{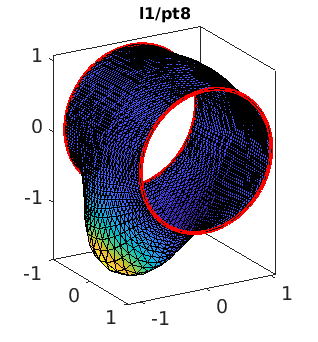}\\[-2mm]
{\sm (e)}\hs{50mm}{\sm (f)}\\
\hs{0mm}\ig[height=52mm]{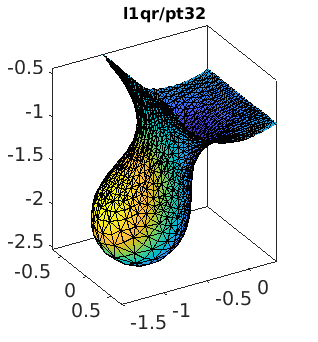}\hs{-4mm}\ig[height=50mm]{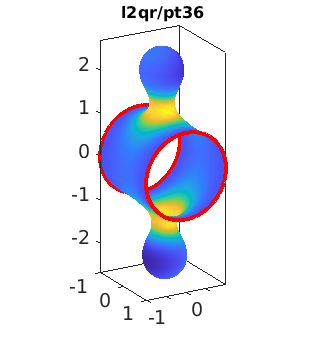}
\hs{-10mm}\ig[height=50mm]{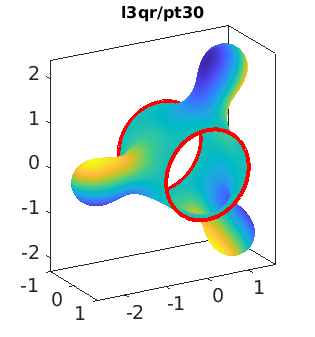}\hs{-2mm}\ig[height=50mm]{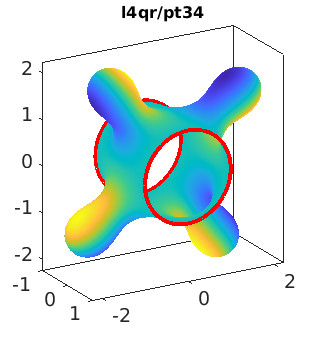}
}
\vs{-4mm}
\caption{{\small %Results from 
{\tt cmds4.m}, 
continuation in $\lam_1$ at $c_0=0.7$. In (a): 
axisymmetric branch  {\tt l0b} (dark blue, sample in (b)), 
{\tt l0} (blue, sample in (c))), and 
four bifurcating branches, {\tt l1} (red, refined to {\tt l1qr}, 
samples in (d) and (e), with zoom), 
{\tt l2} (magenta), {\tt l3} (orange), {\tt l4} (green), samples in (f). 
The kinks of {\tt l1}--{\tt l4} 
are due to adapative mesh refinement every 5th step (after pt10), 
from ${\tt np}=2700$ to ${\tt np}\approx 3600$ at the end of each 
branch.
  \label{hf2}}}
\end{figure}
 
\begin{figure}[h]
\centering
\btab{l}{
\ig[height=57mm]{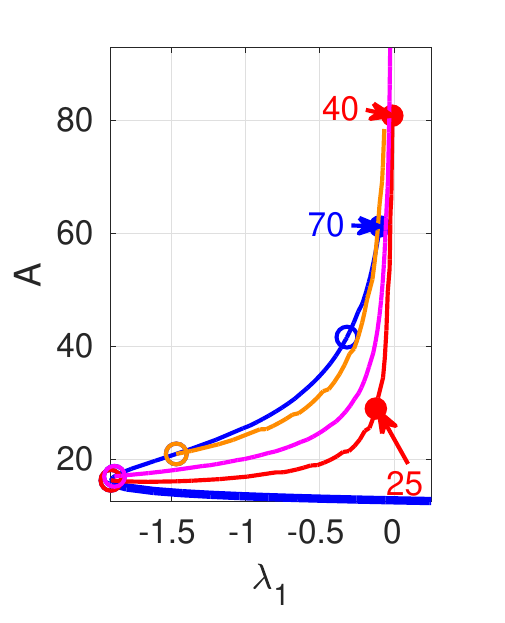}
\hs{-3mm}\ig[height=54mm]{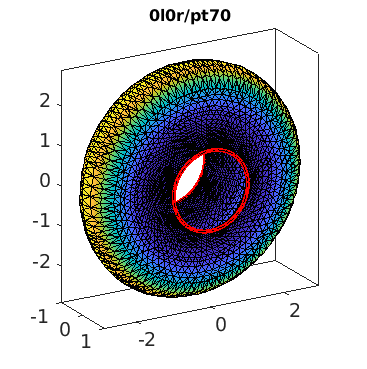}
\hs{-2mm}\rb{5mm}{\ig[height=40mm]{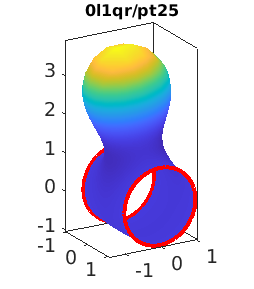}}
\hs{-4mm}\ig[height=57mm]{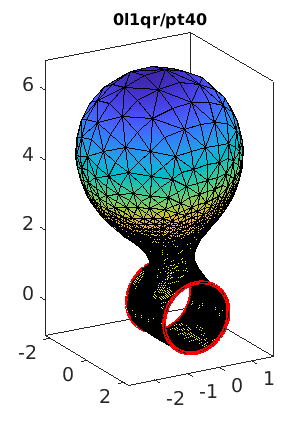}
}
\vs{-4mm}
\caption{{\small %Results from 
{\tt cmds4b.m}, 
continuation in $\lam_1$ at $c_0=0$; 
color code as in Fig.\ref{hf2},  
and, e.g., {\tt np}=4400 at {\tt 0l1qr/pt40}. 
  \label{hf2g}}}
\end{figure}

\begin{figure}[h]
\centering 
\btab{ll}{{\sm (a)}&{\sm (b)}\\[-0mm]
\hs{-2mm}\ig[height=42mm]{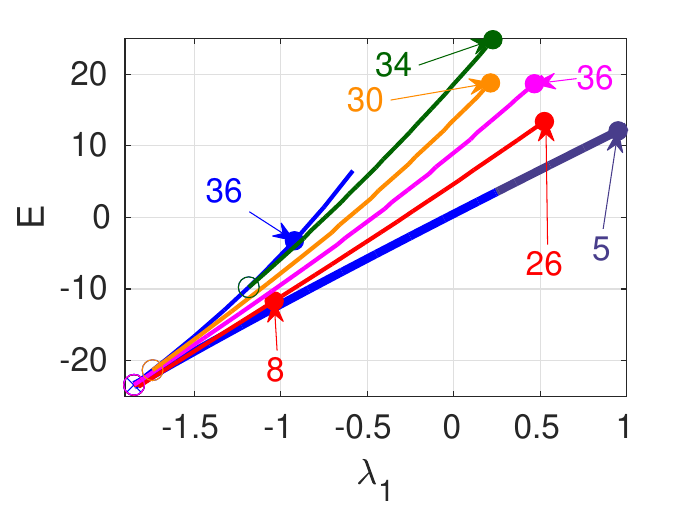}\ig[height=42mm]{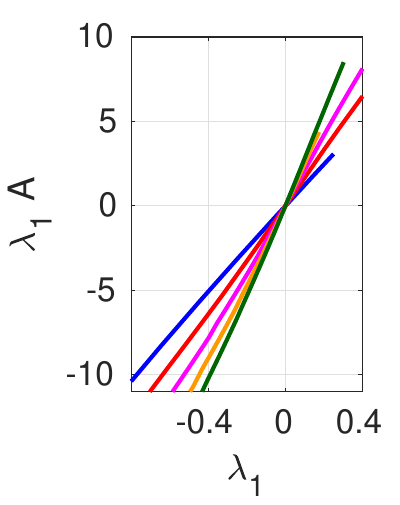}&
\ig[height=42mm]{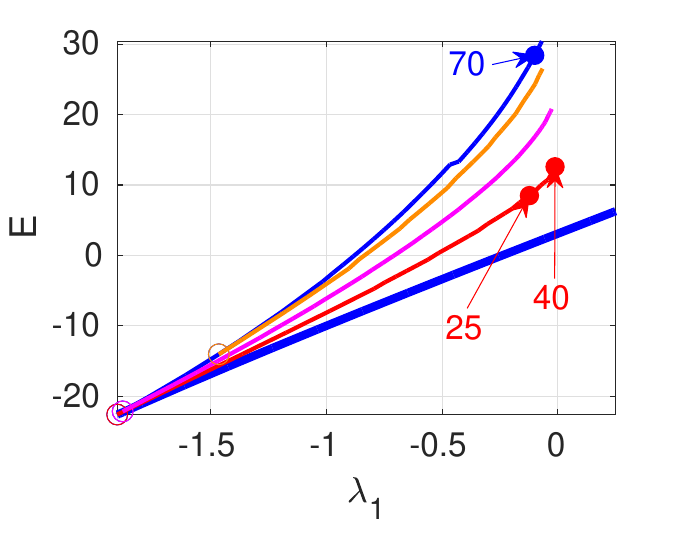}\ig[height=42mm]{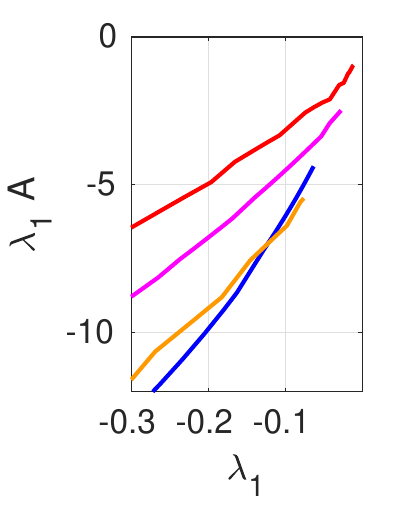}
}
\vs{-4mm}
\caption{{\small Energies $E$ and the part $\lam_1A$ of $E$ along the branches 
of Fig.~\ref{hf2}; (a) $c_0=0.7$, (b) $c_0=0$. 
\label{hf2d}}}
\end{figure}

Even if $\lam_1<0$ is in general unphysical, see Remark \ref{vesrem}(b), 
continuation  
to $\lam_1<0$ yields interesting results, in particular for $c_0>0$. 
For both, $c_0=0.7$ and $c_0=0$, for $\lam_1<0$ the 
solutions start to bulge out in the ``middle'' (near $x=0$), and 
there are folds around $\lam_1=-2$ after which we obtain pronounced 
``tire shapes''. % (Fig.\ref{hf2}(c) and (h)). 
Moreover, we obtain 
bifurcations to $D_m$--symmetric branches, with increasing 
angular wave numbers $m=1,2,3,4,\ldots$. For $c_0=0.7$, 
the bifurcating branches return to $\lam_1>0$, such that 
in Fig.\,\ref{hf2}(a--f) we obtain four new solutions at $\lam_1=0.25$. 
For $c_0=0$, the $D_m$ symmetric branches initially behave similarly, 
but do {\em not} reach $\lam_1>0$ and instead asymptote to $\lam_1=0_-$, 
with arbitrary large $A$, see Fig.\,\ref{hf2g} and again Remark 
\ref{vesrem}(b).  
Therefore we put $c_0=0.7$ first in Fig.\ref{hf2}. 
Nevertheless, for both $c_0=0$ and $c_0=0.7$, the BCs 
do not allow ``large portions of planes'' as solutions and thus 
yield the folds of the blue branches. 
Moreover, $E$ stays bounded below 
(see Fig.\,\ref{hf2d}), and also for $c_0=0$ in Fig.\,\ref{hf2g} 
and Fig.\,\ref{hf2d}(b) it seems that $\lam_1 A\to 0$ 
as $\lam_1\nearrow 0$, i.e., that the area grows slower than $|1/\lam_1|$.

The main issue in both, {\tt cmds4.m} and {\tt cmds4b.m} (which is essentially 
a copy of {\tt cmds4.m}, with a different start (at $c_0=0$) 
of the blue branch and hence directory names starting with 0), is the 
mesh adaptation for the strong budding of the $D_m$ symmetric branches. 
In detail, to compute the branch, e.g., {\tt l1} ($m=1$, 
sample in Fig.\,\ref{hf2}(d)) 
with refinement (sample in (e)) we proceed as follows. 
After branch switching at {\tt l0/bpt1} (double, by rotational 
symmetry, but again we just choose the first of the two kernel vectors) 
we do a few steps without rotational PC, which we then switch on. 
Subsequently, we continue with mesh adaptation based on area 
every five continuation steps. This works very robustly, 
also on the branches with $m=2,3,4$ and similarly for $c_0=0$ in Fig.\,\ref{hf2g}, 
and allows continuation to large buds. 
This shows again that bulging {\em out} (i.e., expanding) 
is typically not problematic wrt meshes. 
The only unbeautiful  effect are the kinks on the $m=1,2,3,4$ 
branches in Fig.\,\ref{hf2}(a), which occur at mesh--refinement. 
To have smoother branches we would need more frequent but less strong 
refinement, but this is clearly not critical.

The branches in Fig.\,\ref{hf2} 
can be continued further by alternating {\tt cont} and {\tt refineX}, 
but this requires some fine--tuning of the 
{\tt cont}--{\tt refineX} loop parameters, and eventually the 
continuation of the branches {\tt l*qr} fails again due to 
bad triangles in the necks, which apparently cannot easily be fixed. 
%by {\tt degcoarsenX}. 
Thus,  to keep 
the script {\tt cmds4.m} simple we stop at the given points. 
%, although the further behavior of these branches (with decreasing $A$) is not clear. 
%Similarly, we stop the continuation at {\tt 0l0r/pt70} as here the main features are already clear. 
Finally we note that the stabilities 
in Fig.\,\ref{hf2} and Fig.\,\ref{hf2g} are as expected, 
namely $\ind(X)=m$ on the $D_m$ branch(es). 

\subsection{Biomembrane caps}\label{hcap}
%\subsubsection{Caps}\label{hcap}
\def\dhome{geomtut/biocaps} 

In the demo {\tt biocaps} we consider BCs \reff{natBCs} for 
the circle $\pa X{=}\{x^2+y^2{=}\al^2\}$ in the $x$--$y$ plane. Thus, 
\begin{subequations}\label{Hcb}
\hual{
&\Delta H+2H(H^2-K)+2c_0K-2c_0^2H-2\lam_1H=0,\label{Hcbp} \\
\intertext{on $X$, and, on $\pa X$,} %=\{X_1^2+X_2^2=\al^2\}$,}
&u=0\text{ (since we keep $\al$ fixed throughout)},\label{cabc1} \\ 
&H-c_0+b\kap_n=0. \label{cabc2} 
}
\end{subequations}
In our experiments we fix $\al=1$ and $\lam_1=1/4$, and first $c_0=1/2$ 
and vary $b$, and want to start with the upper unit hemisphere. Then $H=K=1$ 
(choosing the inner normal for the hemisphere) and 
$\kap_n=1$ and hence \reff{cabc2} requires $b=-1/2$. The files in {\tt biocaps} 
are as usual, with the main adaption in {\tt hcapbra.m} to compute 
$\int_X K\dd S$ for the energy $E$, see Remark \ref{hcaprem}c). 

\begin{figure}[h]
\centering 
\btab{l}{{\sm (a)}\\
\hs{-2mm}\ig[width=48mm]{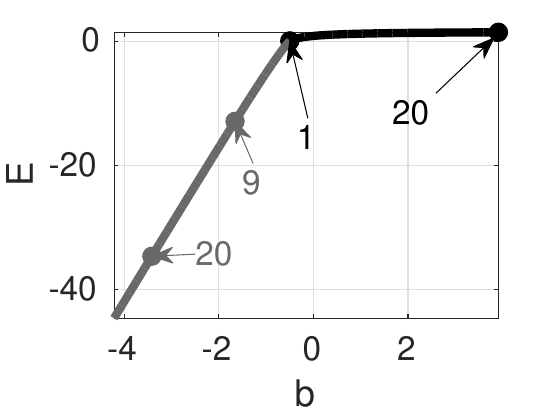}
\rb{5mm}{\hs{-2mm}\ig[width=45mm]{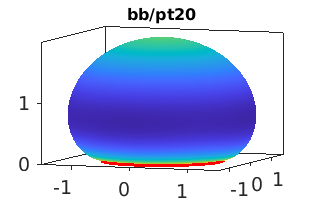}
\hs{-1mm}\ig[width=45mm]{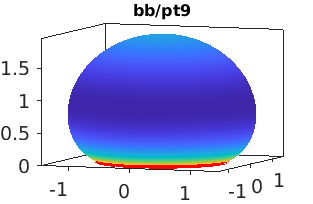}
\hs{-1mm}\ig[width=45mm]{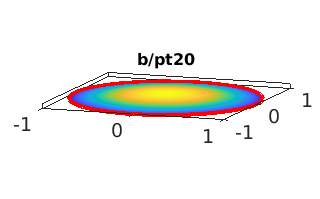}}\\
{\sm (b)}\\
\hs{-2mm}\ig[width=48mm]{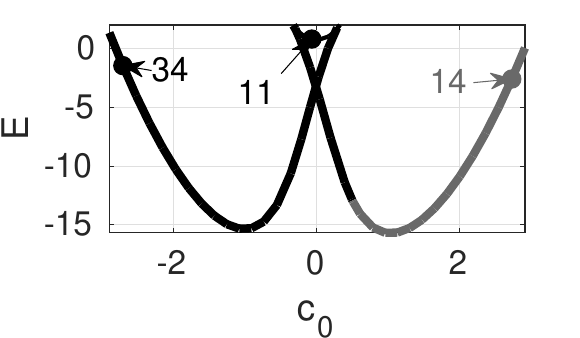}
\rb{5mm}{\hs{-0mm}\ig[width=44mm]{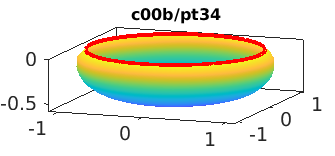}
\hs{-1mm}\ig[width=44mm]{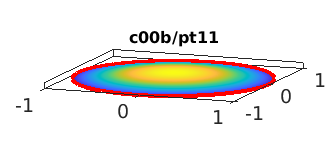}
\hs{1mm}\ig[width=44mm]{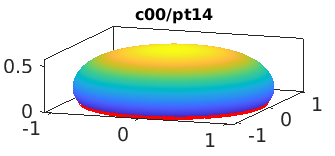}}
}

\vs{-4mm}
\caption{{\small Initial results for \reff{H2},\reff{natBCs} from {\tt biocaps/cmds1.m}. In (a) 
we continue in $b$ starting at {\tt b/pt1} from the unit hemisphere 
with $(\al,\lam_1,c_0,b)=(1,1/4,1/2,-1/2)$, 
to increasing $b$ (branch {\tt b}, black) and decreasing $b$ (branch 
{\tt bb}, grey). On {\tt b} we go to the flat disk (last sample), 
while on {\tt bb} the hemisphere bulges out. This is mainly intended 
for later continuation in $c_0$, and in (b) we do so starting from {\tt bb/pt9} 
at $b\approx -1.66$. This gives the double well shape for $E$, with a short 
unstable segment between the two folds at $c_0\approx \pm 0.33$. 
See Fig.~\ref{hcaf2} for the cases of $b\approx -3.4$ ({\tt bb/pt20}) and 
$b\approx -4$. 
\label{hcaf1}}}
\end{figure}

\begin{figure}[H]
\centering 
\btab{l}{{\sm (a)}\hs{45mm}{\sm (b)}\\[-0mm]
\hs{-6mm}
\btab{l}{\ig[width=52mm]{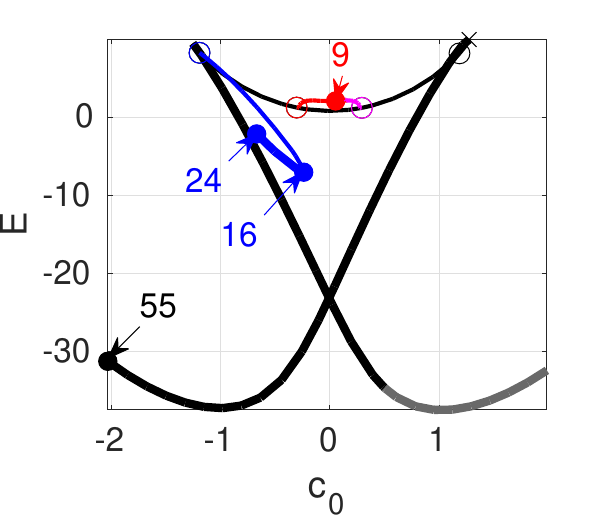}\\
\ig[width=52mm]{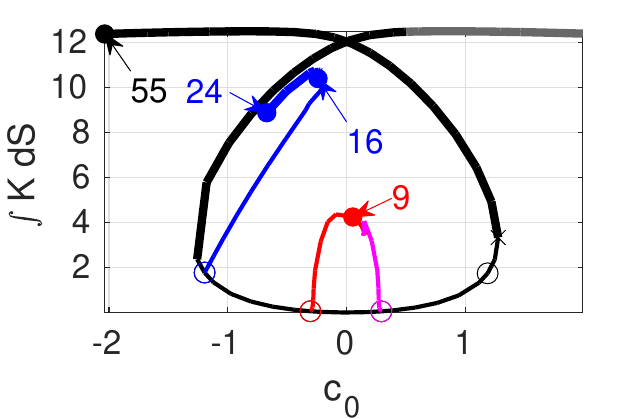}}
\rb{0mm}{\btab{ll}{
\hs{-5mm}\ig[height=28mm]{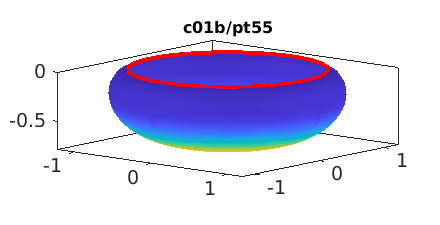}
&\hs{-5mm}\ig[height=28mm]{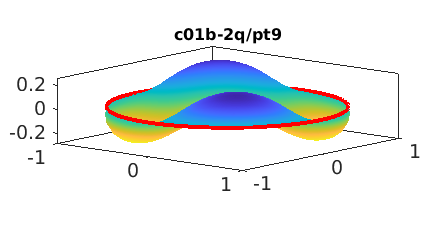}\\[-3mm]
\hs{-5mm}\ig[height=56mm]{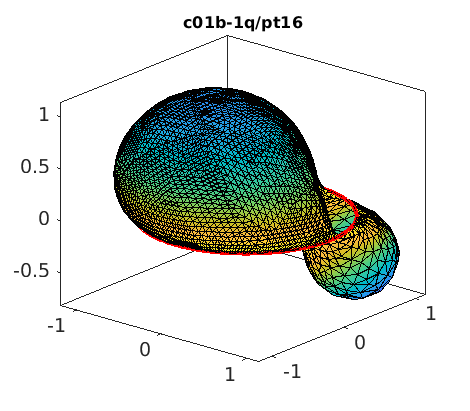}
&\hs{-5mm}\rb{0mm}{\ig[height=56mm]{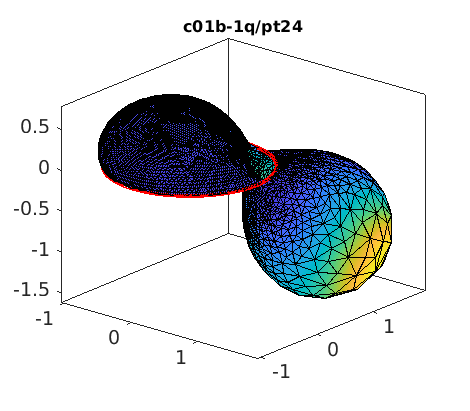}}}}\\[-4mm]
{\sm (c)}\\[-0mm]
%\hs{-2mm}\ig[width=44mm]{hca/bdE4}
\hs{-2mm}\ig[width=46mm]{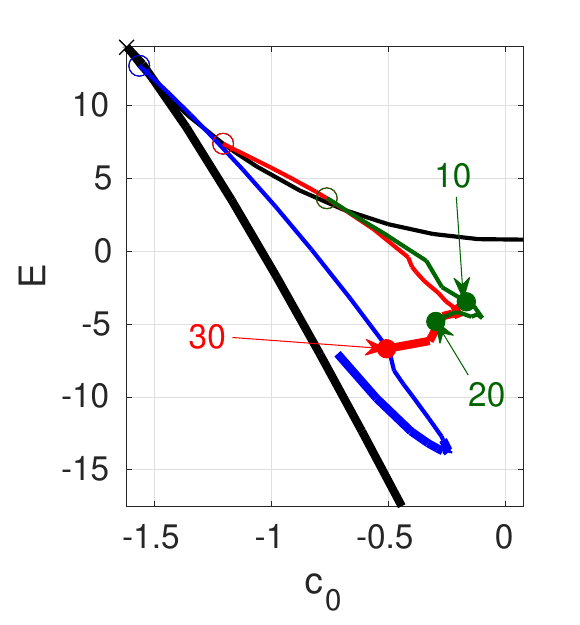}
\hs{2mm}\rb{5mm}{%\ig[width=44mm]{hca/c1b-1-26}
\hs{-3mm}\ig[width=46mm]{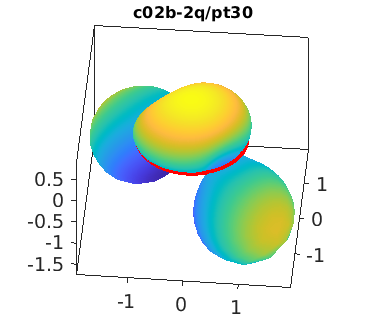}
\hs{-3mm}\ig[width=46mm]{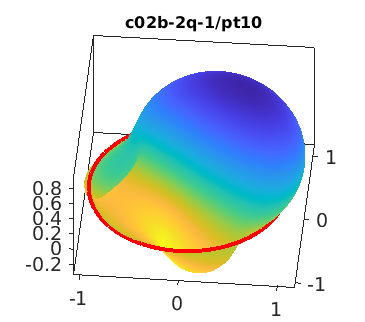}
\hs{-3mm}\ig[width=46mm]{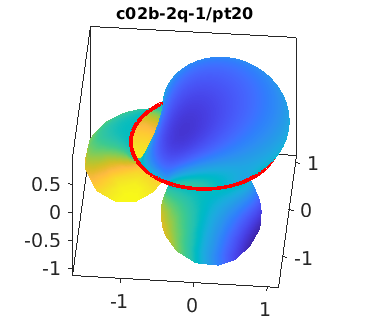}
}
}
\vs{-4mm}
\caption{{\small 
Further results from {\tt cmds1.m}. (a,b) 
Continuation in $c_0$ from {\tt bb/pt20}, $(\al,\lam_1,b)=(1,0.25,-3.4)$, 
starting from $c_0=0.5$, branches {\tt c01b} (black, to decreasing $c_0$) and 
{\tt c01} (grey, to increasing $c_0$}). There are two BPs on the unstable 
part of {\tt c01b} for $c_0<0$, and the symmetric BPs for $c_0>0$. 
The blue branch {\tt c01b-1q} has azimuthal wave number $m=1$ 
and is stable after its fold. 
The red branch {\tt c01b-2q} has $m=2$ 
and connects to the symmetric BP at $c_0>0$. The 2nd plot in (a) 
shows where the part $b\int K\dd S$ of $E$ becomes dominant, 
taking into account the rather large $|b|$. 
(c) Similarly starting at {\tt bb/pt24} with $b\approx -4$; zoom of BD 
near upper left fold of the branch {\tt c02b} (black) 
similar to {\tt c01b} from (a). The blue branch is 
qualitatively as in (a),  but now the $m=2$ branch {\tt c02b-2q} (red) 
also folds back giving stable solutions, and there is a secondary BP on it, 
giving the green branch {\tt c02b-2q-1}. 
\label{hcaf2}}
\end{figure}

Fig.\ref{hcaf1}(a) shows the continuation in $b$. %, starting at {\tt pt1}. 
This is mainly intended for subsequent continuation in $c_0$ at 
negative $b$, and (b) shows the case of $b\approx -1.66$. The problem 
is symmetric under $(c_0,X_3)\mapsto -(c_0,X_3)$, and in particular at 
$c_0=0$ we have the flat disk as an exact solution (for any $b$). 
See {\tt c00b/pt11} for a nearby solution with $c_0\approx  -0.07$, 
which lies between two folds with exchange of stability. This unstable 
part will feature interesting bifurcations to non--axisymmetric 
branches at more negative $b$, 
see Fig.\,\ref{hcaf2}, while the remainder(s) of the axisymmetric 
branches are all stable, with the samples {\tt c00b/pt34} and {\tt c00/pt14} 
in  \ref{hcaf1}(b) showing the typical behavior at strongly negative or 
positive $c_0$, respectively.

\brem\label{hcaprem}{\rm 
The main  numerical challenges and used tricks for \reff{Hcb} 
are: 

a) To have a mesh as regular as possible for the initial hemisphere 
we use a classical subdivision and projection algorithm, 
see {\tt subdiv\_hemsphere.m}, modified from {\tt subdivided\_sphere.m} 
from the \gptool. Subsequently, we do one mesh--refinement with 
{\tt e2rsbdry} as selector, as a good resolution near $\pa X$ turns 
out helpful later. The initial mesh then has {\tt np=2245} nodes, 
which for instance at {\tt c01b-1q/pt24} is refined (by triangle 
area, i.e., {\tt e2rsA}, followed by {\tt retrigX}) to {\tt np=4475}. 
The mesh quality in all our solutions stays quite good, i.e., 
$\delm<20$ for all solutions, and mostly $\delm<10$. 

b) The boundary $\ga=\pa X$ is parameterized by arclength as 
$\ga(\phi)=\al(\cos(\phi/\al),\sin(\phi/\al),0)$. Then 
$\kap=\ga''=-\ga/\al$ and the normal curvature on 
$\pa X$ reads  $\ds \kap_n=-\frac 1 \al \spr{N,X}$, which 
is used to implement the BCs \reff{cabc2}. 

c) The ``integral'' {\tt sum(K)} over the 
discrete Gaussian curvature ${\tt K}$ always evaluates to $2\pi\chi(X)$, 
cf.~Footnote \ref{GBfoot}. 
%i.e., the definition of {\tt K} (necessarily) ignores the geodesic curvature $\kap_g$, i.e., assumes $\kap_g=0$. 
Thus we once more use 
Gauss-Bonnet $\ds\int_X K\dd S=2\pi\chi(X)-\int_{\pa X}\kap_g\dd s$, where 
$\kap_g=\sign(N_3)\frac 1 \al \|N\times \ga\|$, 
to compute the energy $E$, 
and in {\tt hcapbra.m} we evaluate $\ds\int_{\pa X}\kap_g\dd s$ 
by a trapezoidal rule. 
Another alternative is to zero out 
$K|_{\pa X}$ before evaluating  {\tt sum(K)}, and for comparison 
we also put that into {\tt hcapbra.m}. In our computations we find a 
relative error $<0.01$ (usually much smaller) between the two methods, 
when doing the initial refinement step near $\pa X$ from (a). 
}\eex\erem 

In Fig.\,\ref{hcaf2} we repeat the continuation in $c_0$ from 
Fig.\,\ref{hcaf1}(b) at more negative $b$, namely $b\approx -3.4$ in (a) 
and $b\approx -4$ in (b). For lower $b$, the unstable part of the 
$c_0$ continuation expands, and we find two (or more, for even lower $b$) 
BPs between the left fold and $c_0=0$, with azimuthal wave numbers 
$m=1$ and $m=2$. As before, these bifurcations are double by $S^1$ symmetry, 
and to continue the bifurcating branches we set the usual rotational PC 
after two steps. The blue $m=1$ branch then behaves similarly 
in (a) and (b), i.e., it becomes stable after a fold at $c_0\approx -0.2$ 
($b=-3.4$) resp.~$c_0\approx -0.27$ ($b=-4$). 
%, where also the further continuation becomes more difficult wrt meshes). 
However, the $m=2$ branch behaves differently: For $b=-3.4$ it connects 
to the symmetric BP at $c_0>0$. For $b=-4$, the red branch {\tt c02b-2q} 
first shows 
a secondary BP to a branch ({\tt c02b-2q-1}, green) with broken 
$\Z_2$ symmetry, and then shows a fold at 
$c_0\approx -0.21$ where it becomes stable. 
The branch {\tt c02b-2q-1} also shows a fold, at $c_0\approx -0.11$, after 
which however one unstable eigenvalue remains, i.e., $\ind(X)=1$ at, e.g., 
{\tt c02b-2q-1/pt20} (last sample in (c)). The somewhat non--smooth shape 
of the red and green branches is due to repeated 
and heavy mesh refinement, to, e.g., $n_p=5560$ at {\tt c02b-2/pt30}.

\begin{figure}[h]
\centering 
\btab{l}{{\sm (a)}\\[0mm]
\hs{-2mm}\ig[width=40mm]{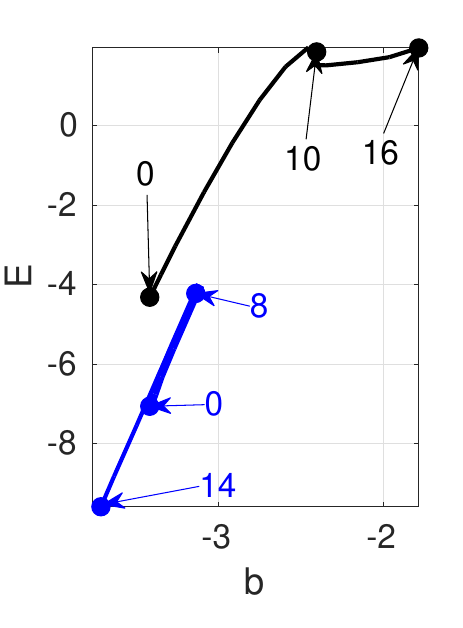}
\rb{25mm}{\btab{l}{\ig[width=45mm]{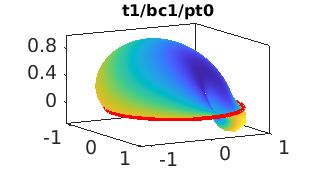}
\ig[width=45mm]{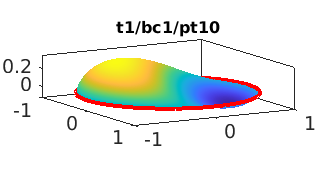}
\ig[width=45mm]{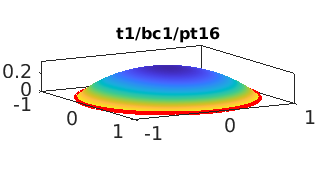}\\
\ig[width=45mm]{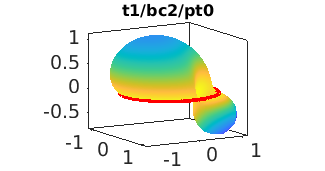}\ig[width=45mm]{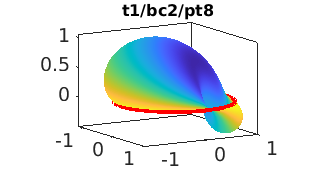}
\ig[width=45mm]{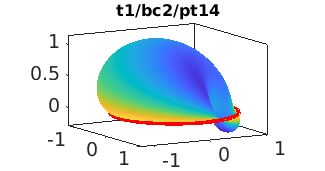}}}\\[-4mm]
{\sm (b)}\\
\hs{-2mm}\ig[width=40mm]{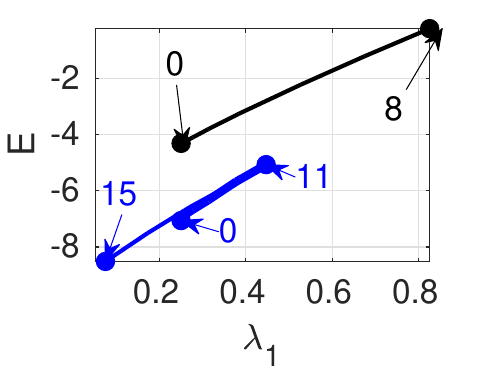}
\hs{2mm}\rb{3mm}{\ig[width=45mm]{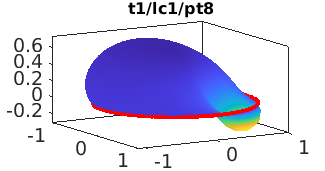}
\ig[width=45mm]{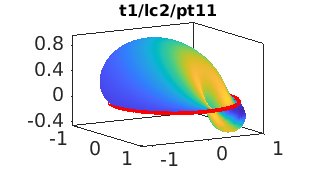}
\ig[width=45mm]{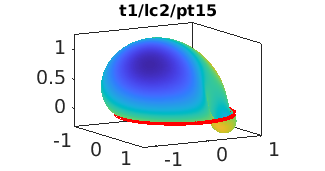}}
}

\vs{-2mm}
\caption{{\small {\tt cmds2.m}, experiments with continuation of 
$m=1$ solutions from Fig.~\ref{hcaf2} in $b$ (a), and in $\lam_1$ (b). 
The starting points in (a) are from $c\approx -0.35$ for the black branch 
{\tt bc1}, and from $c\approx -0.24$ for the blue branch {\tt bc2} (this is 
the same solution as {\tt c01b-1q/pt16}).  
The continuation in $\lam_1$ in (b) has the same 
starting points and the same colors. The unstable branch (now 
with fixed $(c_0,b)=(-0.35,-3.4)$) continues to large $\lam_1$, 
but the blue branch with fixed  $(c_0,b)=(-0.24,-3.4)$ has a fold at 
$\lam_1\approx 0.44$. 
\label{hcaf8}}}
\end{figure}

In Fig.\,\ref{hcaf8} we continue some $m=1$ solutions from the blue 
branch in Fig.\,\ref{hcaf2}(b) in $b$ and in $\lam_1$, with fixed 
$c_0$. The unstable black solutions with $c_0\approx -0.35$ in (a) 
continue to $b\approx -1.8$, where the branch bifurcates from 
the axisymmetric branch. This is in contrast to Fig.\,\ref{hcaf1} 
where for the continuation in $b$ at $c_0=1/2$ no bifurcations from 
the axisymmetric branch were found. 
The blue branch in Fig.\,\ref{hcaf8}(a) with $c_0\approx -0.24$ 
shows a fold at $b\approx -3.11$. 
From these and similar experiments, $m=1$ solutions do not seem to 
exist for $b$ somewhat larger than $-2$, and no stable 
ones for $b$ larger than $-3$. In Fig.\,\ref{hcaf8}(b), the black unstable branch 
continues to large $\lam_1$, but the stable blue branch again shows a 
fold, at moderate $\lam_1\approx 0.44$. 
This indicates that solutions of this type also do 
not continue to ``large'' $\lam_1$. 
Nevertheless, 
while the physical (biological) significance of the parameter regimes 
and solutions in Figs.~\ref{hcaf1}--\ref{hcaf8} certainly 
needs to be discussed, 
mathematically we have obtained some stable 
non--axisymmetric solutions.

\section{Summary and outlook}\label{sosec} 
We explained the basic setup of the \pdep\ extension library {\tt Xcont} 
for continuation of 2D submanifolds $X$ (surfaces) of $\R^3$, 
and gave a number of examples. These were partly introductory,  
e.g.~the spherical caps in \S\ref{scsec}, partly classical, e.g., 
the Enneper minimal surface in \S\ref{ensec}, the nodoids in 
\S\ref{nodnumDBC} 
and \S\ref{nodnumpBC}, and the TPS in \S\ref{tpssec}, 
and partly rather specific such as the Plateau problem in \S\ref{bc2sec} 
or the 4th order Helfrich type cylinders and caps in 
\S\ref{wsec} and \S\ref{hcap}. Besides 
\cite{bru18}, and to some extent \cite{Bra96}, 
there seem to be few numerical continuation 
and bifurcation experiments for such geometric problems for 2D surfaces, 
i.e., without imposing some axial symmetry, 
and we are not aware of a general %and user--friendly 
software for such tasks. 

The basic setup for all our problems (except \reff{H2}, which is of 
4th order) is similar: We consider CMC surfaces, which mainly differ 
wrt constraints and/or boundary conditions. Along the way we explained 
a number of techniques/tricks which we expect to be crucial in 
many applications. 
%\bci 
%\item 
A major problem for 
continuation (over longer parameter regimes) is the mesh handling 
as $X$ changes and hence the mesh distorts. 
We explained how this (often) can be abated via moving of mesh points 
({\tt moveX}), refinement ({\tt refineX} which sometimes should be 
combined with re--triangulation by {\tt retrigX}) and 
coarsening ({\tt coarsenX}), 
and coarsening of degenerate triangles ({\tt degcoarsenX}), although 
the choice of the parameters controling these functions often 
requires some trial and error.  In any case, $X$ bulging out (increasing area) 
is usually 
harmless, but bulging in (the development of necks) is more challenging. 
%Additionally,  restarting as for the nodoids in \S\ref{nodnumDBC} is an option, 
%if the desired branch can be approximated. 
%\item Symmetries, leading to nontrivial kernels (and higher 
%multiplicities of BPs) can be handled with suitable phase conditions. 
%\eci 

This is a first step. 
With the demos we hope to give a pool of applications which 
users can use as templates for their own problems, and we are 
curious what other applications users will consider, and of course are happy  
to help if problems occur. As indicated above, 
our own further research, to be presented elsewhere, mostly since bifurcation 
diagrams become complicated, but also since some of the numerics 
require further tricks, includes: 
\bci 
\item Further classical minimal surfaces (and CMC companions) such 
as Schwarz H 
and Scherk surfaces (surface families); 
\item Bifurcations of {\em closed} vesicles in Helfrich type problems. 
\item Coupling of membrane curvature and reaction--diffusion equations 
for proteins as in, e.g., \cite{TN20}. 
\eci

\appendix 
%\section{Spheres, hemispheres, and VPMCF}\label{hssec0}
\section{Spheres, hemispheres, VPMCF, and an alternative setup}\label{hssec0}
\subsection{Spheres}\label{spsec}
The demo {\tt spheres}, containing only the (somewhat minimally necessary) 
files {\tt sphereinit.m}, {\tt sG.m}, {\tt getM.m},  and 
{\tt cmds1.m}, 
is mostly meant to illustrate volume preserving mean curvature 
flow (VPMCF) near  spheres, see Fig.\,\ref{sf1}.%
\footnote{Additionally, the demo contains {\tt convtest.m} for convergence 
tests, see Fig.\,\ref{nf1}.} 
We refer to {\tt sphereinit.m} and 
{\tt cmds1.m} for comments (and to {\tt geomflow.m} and 
{\tt vpmcff.m} from {\tt libs/Xcont} for the VPMCF) 
and here only note: 
\begin{figure}[ht] %H]
\centering
\btab{l}{{\sm (a)}\hs{80mm}{\sm (b)}\\
\hs{0mm}\ig[width=0.21\tew]{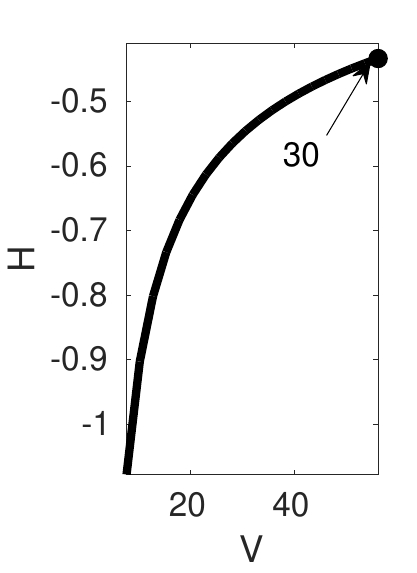}
%\hs{2mm}\ig[width=0.30\tew]{spics/H30}\ig[width=0.30\tew]{spics/K30}
\hs{4mm}\ig[width=0.18\tew]{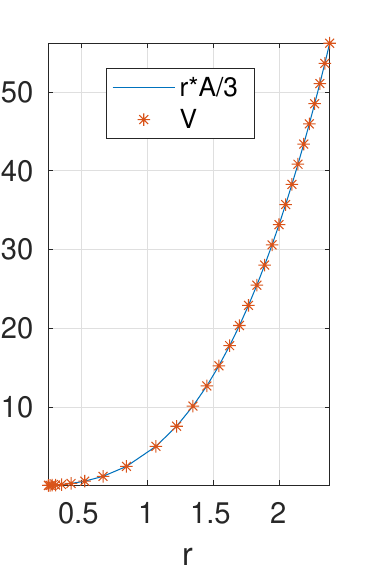}
\hs{8mm}\rb{0mm}{\ig[width=0.27\tew]{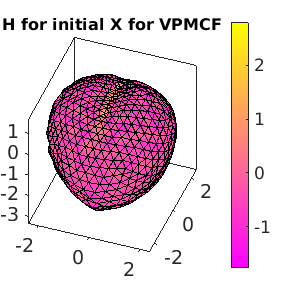}}\\[-3mm]
{\sm (c)}\\[-3mm]
\hs{0mm}\ig[width=0.25\tew]{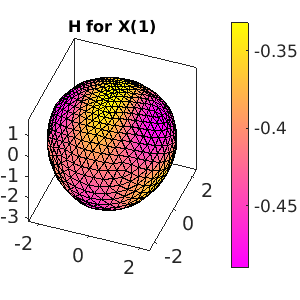}
\hs{0mm}\ig[width=0.27\tew]{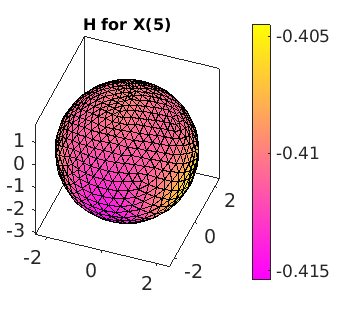}
\hs{3mm}\rb{0mm}{\ig[width=0.13\tew]{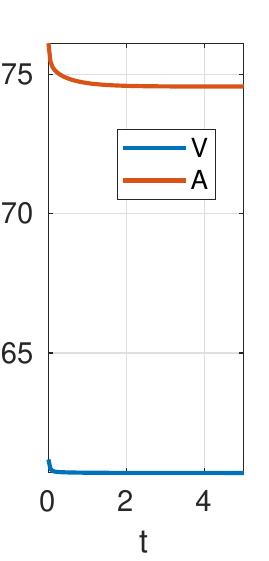}}
}\
\vs{-3mm}
\caption{\small Results from {\tt sphere/cmds1.m}. (a) Continuation 
of a sphere in $V$, with comparison of  $A$ and $V$. 
(b) An IC for VPMCF from a perturbation of {\tt S3/pt30}, 
with (c) solutions at $t=1$ and $t=5$, and a time--series for $V,A$; 
$V$ is conserved to within $0.5\%$. 
\label{sf1}}
\end{figure}
\bcen 
\item The comparison between $rA/3$ and $V$ in Fig.\,\ref{sf1}(a) shows 
a very small error which indicates that the solutions are good approximations 
of spheres. 
\item For convex closed initial $X$ (meaning that $X=\pa\Om$ for 
a convex domain $\Om\subset\R^3$), the VPMCF converges to a sphere. 
See also, e.g., \cite{ESi98} for theoretical background. 
This also holds for ``slightly'' non--convex initial $X$ as 
in Fig.\,\ref{sf1}(b).% 
\footnote{In detail, $X(0)$ here is obtained as $X(0)=S_{r_0}+0.4(\sin(\vt)(|x|-{r_0})+\xi)N$, 
where $\vt$ is the azimuth, and $\xi=0.2({\tt rand}-0.5)$ with ${\tt rand}\in[0,1]$ a \mlab\ 
random variable on each node.} 
\item  
Our (explicit Euler) implementation of the VPMCF does not 
conserve $V$ exactly, but with ``reasonable accuracy'', i.e.: 
Even for quite ``non--spherical'' initial $X(0)$, the error 
$|1-V_\infty/V_0|{<}0.01$, where $V_\infty=\lim_{t\to\infty}V(t)<V_0$ 
in all our tests, i.e., $V_\infty$ is always slightly smaller than 
$V_0$.%
\footnote{This ``volume--accuracy'' of {\tt geomflow} for VPMCF 
depends weakly on the Euler--stepsize {\tt dt}, and more strongly 
on the fineness of the discretization of $X$. This can be checked 
by changing {\tt sw} in {\tt spheres/cmds1.m} for initializing $X$.} 
\item During continuation, the position of {\tt X} is not fixed, i.e., 
we have the threefold (in the discrete setting approximate) 
translational invariance in $x,y,z$, and hence always a three-dimensional 
(approximate) kernel. This could be removed by suitable translational 
PCs (see the demo {\tt hemispheres}). However, the approximate 
kernel is not a problem for the continuation here since 
in the Newton loops the right hand sides are orthogonal to this kernel.  
As here we are mostly interested in the VPMCF, for which the 
translational invariances are irrelevant, we refrain from 
these PCs to keep the demo sleek and simple. % as possible. 
\ecen

\subsection{Hemispheres}\label{hssec}
\begin{figure}[ht] %H]
\centering
\btab{l}{{\sm (a)}\hs{80mm}{\sm (b)}\hs{30mm}{\sm (c)}\\
\hs{0mm}\ig[width=0.2\tew]{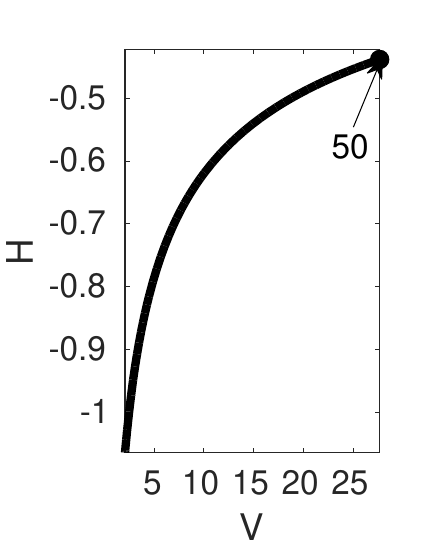}
\hs{-3mm}\rb{0mm}{\ig[width=0.28\tew]{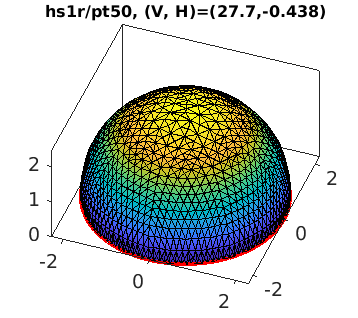}}
\hs{2mm}\ig[width=0.18\tew]{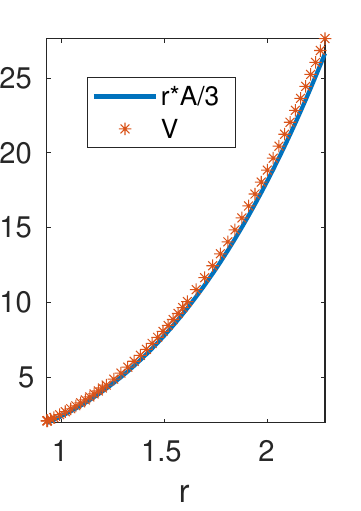}
\ig[width=0.3\tew]{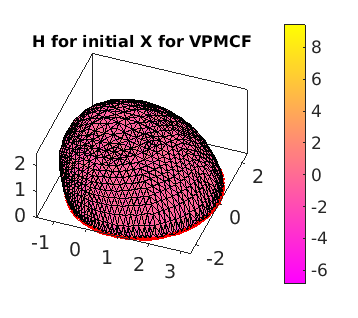}\\[-2mm]
{\sm (d)}\\
\hs{-0mm}\ig[width=0.3\tew]{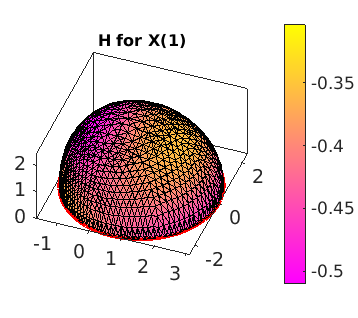}
\hs{-0mm}\ig[width=0.3\tew]{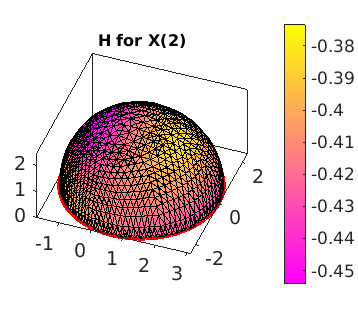}
\hs{8mm}\rb{0mm}{\ig[width=0.18\tew]{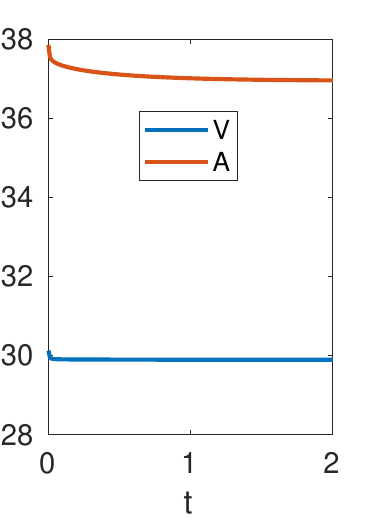}}
}
\vs{-3mm}
\caption{\small Results from {\tt hemisphere/cmds1.m}. (a) Continuation 
of hemisphere in $V$, with sample at end. 
(b) Comparison of $A$ and $V$ for (a). 
(c) A perturbation of {\tt hsr1/pt50} as IC for VPMCF. 
(d) Solutions at $t=1$ and $t=2$, and time--series for $V,A$ 
from (c,d); $V$ conserved to within $0.8\%$. 
\label{sf2}}
\end{figure}

In the demo {\tt hemispheres} we continue  in volume $V$ hemispheres $X$ 
sitting orthogonally on the $z=0$--plane, i.e., 
\huga{\label{sNBC}
\text{$\pa_r X_3=0$ where $r=\sqrt{x^2+y^2}$,}
}  
and test VPMCF for perturbations of such hemispheres, see Fig.\,\ref{sf2}. 
Additionally we use the PCs 
\huga{\label{hspc}
\int_X N_1\dd S=\int_X N_2\dd S=0, \quad X=X_0+uN_0, 
} 
to fix the translational invariance. 

\taskip
\begin{table}[ht]\caption{
Files in {\tt pde2path/demos/geomtut/hemispheres}. 
\label{hstab1}}
\centering\vs{-2mm}
{\small 
\begin{tabular}{l|l}%p{0.2\tew}|p{0.82\tew}} 
%file&remarks\\ \hline
{\tt cmds1.m}&continuation in $V$, and VPMCF flow test. \\
{\tt hsinit.m, sGhs.m}&init and rhs\\
{\tt qf2.m, qjac2.m}&PC \reff{hspc}, and derivative \\
{\tt getN.m}&mod of getN, correction at $z=0$\\
{\tt diskpdeo2.m}&mod of (\pdep--)default  {\tt diskpdeo2.m} to have 
a finer mesh at $r=1$. 
\end{tabular}
}
\end{table}
\teskip

Compared to the spheres in \S\ref{spsec} this requires a few more files, 
listed in Table \ref{hstab1}. The PCs \reff{hspc} (and $u$--derivatives) 
are implemented in {\tt qf2, qjac2}, and the rhs $G(u)=H-H_0$ in 
{\tt sGhs} is augmented to $\wt G(u)=H-H_0+s_x N_1+s_y N_2$ 
with multipliers $s_x,s_y$. These stay $\CO(10^{-6})$ during continuation, 
and the only effect of the construction is that 
the 2D kernel of translational invariance of $\pa_u G$  is removed from 
the linearization of the extended system $(\wt G,q)$. See the end of 
{\tt cmds1.m} for further comments.%
\footnote{In {\tt cmds1.m} we also monitor the ``positions'' $(x_0,y_0)$ 
of $X$ defined as $x_0=\int_X X_1\dd S,\, y_0=\int_X X_2\dd S$. These behave 
as expected, namely: very small drifts for continuation without PCs, 
but no drifts for continuation with the PCs \reff{hspc}.} 

The BCs \reff{sNBC} allow motion of $\pa X$ in 
the ``support--plane'' $z=0$, and 
are implemented as %{\tt sGhs} as
$$\text{ {\tt r(idx)=grXz*(X(idx,1).\^{}2+X(idx,2).\^{}2)}($\stackrel!=0$)}$$ 
in {\tt sGhs}, where as usual {\tt idx} are the boundary indices, and  
{\tt grXz} is the $z$--derivative (operator) 
on $X$, as before obtained from {\tt grX=grad(X,p.tri)} 
and interpolation to the nodes. This forces the $x,y$--coordinates of 
the points on $X$ directly above the $z=0$ layer to also 
fulfill $x^2+y^2=1$, i.e., we obtain a ``cylindrical socket''  
for the hemispheres. To mitigate this effect, in {\tt hsinit} we 
initialize with a somewhat specialized mesh over the unit disk $D$ 
with higher density towards $\pa D$, which is then mapped to the 
unit hemisphere via $z=\sqrt{1-x^2-y^2}$. Nevertheless, after continuation 
to larger $V$, which is combined with some mesh--adaptation by area, 
we obtain a mismatch between $rA/3$ and $V$, 
see Fig.\,\ref{sf2}(a,b), and compare Fig.\,\ref{sf1}(b).  

In Fig.\,\ref{sf2}(c,d) we give an example of VPMCF from a perturbation 
of {\tt hs1r/pt50}, here of the form $X|_{t=0}=X+0.4(\cos(\vt)(\max(z)-z
+0.1({\rm rand}-0.5))N$, where $\vt$ is the angle in the $z=0$ plane. 
We {\em do not use any BCs}  
in the rhs {\tt vpmcff.m}. Instead, we use the correction 
\huga{\label{hsNc}
\text{{\tt N(p.idx,3)=0; N=normalizerow(N);}} 
}
in {\tt getN.m}, to let $N$ lie in the $x$--$y$--plane, cf.~Remark \ref{Prem1}(b) 
for a similar trick. 
Without \reff{hsNc}, the {\em continuation} of hemispheres works as before, 
but the solutions of the VPMCF start to (slowly) lift off the $z=0$ plane, 
and then (quickly) evolve towards a planar $X$ such that also $V\to 0$. 
With \reff{hsNc}, the solutions flow back to hemispheres, 
and $V$ is again conserved up to $0.5\%$, and ``after convergence''  
(e.g., from $X|_{t=2}$) we can start continuation again, showing 
consistency. 
%In summary, the demo {\tt hemispheres} 

%\input{hs}

\subsection{Spherical caps via 2D finite elements }\label{sc2sec}
\def\dhome{geomtut/spcap2} 
For the sake of completeness and possible generalization, 
we show how the spherical caps with DBCs 
can be treated in a classical FEM setting. 
Let $\Om\subset\R^2$ be a bounded domain and 
$X_0$ be a surface with parametrization $\phi_0:\Omega\ra\R^3$, 
and as before define a new surface via 
$X=X_0+u\,N_0$, $u:\Om\to\R$.  Then \reff{CMC} reads  
\huga{\label{feq2}
	G(u,H,V)=\begin{pmatrix}
		H(X)-H\\
	V(X)-V
	\end{pmatrix}=0, \text{ with $u|_{\pa \Om}=0$}. 
}
This is a quasilinear elliptic equation for $u:\Om\to\R$, 
and after solving \reff{feq2} we can again update $X_0=X_0+uN_0$ and repeat. 
To discretize \reff{feq2} we now use the FEM 
in $\Om$. The main issue is how to compute $H$ (and $A$ and $V$ and 
similar quantities) without the \gptool, 
and Table \ref{sctab1} lists the pertinent files. 

\taskip
\begin{table}[ht]\caption{
Overview of files in {\tt geomtut/spcap2}. \label{sctab1}}
\centering\vs{-2mm}
{\small 
\begin{tabular}{p{0.16\tew}|p{0.84\tew}} 
{\tt cmds1}&Main script; initialization, continuation, plotting.\\
{\tt spcapinit}&Initialization, rather standard, except for {\tt p.X=[x,y,0*x]} for the initialization of $X$.\\
{\tt oosetfemops}&Setting mass matrix $M$, and first order differentiation 
matrices which are used to compute $H$ in {\tt sG} (via {\tt getmeancurv}).\\
{\tt sG, qV}&rhs and volume constraint.\\ 
{\tt getN}&compute normal vector. \\
{\tt getA, getV}&compute $A(X)$ and $V(X)$, overload of {\tt Xcont} functions.\\
{\tt getmeancurv}&$H$ from \reff{h1}, see also  {\tt get1ff, get2ff} for 
1st and 2nd fundamental forms.\\
{\tt e2rs}&ElementToRefineSelector function, based on triangle areas.\\
{\tt cmcbra, pplot}&like {\tt Xcont/cmcbra} and {\tt Xcont/pplot}, 
but overloaded since {\tt p.tri} not present here. \\
{\tt getGupde}&overload of library function to deal with larger bandwidth. 
\end{tabular}
}
\end{table}
\teskip

To implement $H(X)$, $X=X_0+uN$, we here directly use the definition 
\hugast{%\label{h1}
H=\frac 1 2 \frac{h_{11}g_{22}-2h_{12}g_{12}+h_{22}g_{11}}{g_{11}g_{22}-g_{12}^2}
}
of the mean curvature 
based on the fundamental forms of $X$. This is brute force
and in particular neither confirms to a weak (FEM) formulation
 nor allows simple Jacobians, see Remark \ref{hrem1}. 
%Similarly, using $\ds |V|=\int_V 1\dd x=\frac 1 3\int_V \mdiv x \dd x=\frac 1 3\left|\int_{\pa V} \spr{x,n}\dd S\right|$ and using that $\spr{x,n}=0$ in the $x$--$y$ plane, we can implement {\tt getV, qf} and {\tt qVjac}.  

\brem\label{hrem1}
{\rm a) In {\tt getN} and {\tt get1ff} we compute the 
{\em nodal} values for $X_x$ and $X_y$ (via weighted averages of 
the adjacent triangles), and the corresponding 
nodal values of $N$ and $g_{ij}$; this is needed here as 
$X_x$ and $X_y$ appear nonlinearly in $N$, and similarly $g_{ij}$ 
appear nonlinearly in $H$. On the other hand, the second derivatives 
$h_{ij}$ only appear linearly in $H$, and hence we 
take the second derivatives in {\tt get2ff} using the {\em element} 
differentiation matrices ${\tt Kx}$ and ${\tt Ky}$. Thus, ${\tt H}$ 
in {\tt r=-2*H+p.mat.M*(H0*ones(p.nu,1))} in {\tt sG} is (an approximation of) 
the element wise mean curvature, and hence ${\tt H0}$ must also be multiplied 
by the mass matrix {\tt p.mat.M}. 

b) Since $G$ as implemented uses products 
of differentiation matrices, the associated Jacobian $\pa_u G$ has 
more bandwidth than before, i.e., the sparsity pattern $S$ of $\pa_u G$ 
is that of $M^2$, rather than that of $M$, and thus we use 
${\tt S=M^2>0}$ in {\tt getGupde}. 
%However, the main bottleneck for numerical Jacobians is $\pa_u q$, for which we therefore provide a function. 
}\eex \erem

\hulst{}{\dhome/spcapinit.m}
\hulst{}{\dhome/oosetfemops.m}
\hulst{}{\dhome/sG.m}
\hulst{}{\dhome/getmeancurv.m}
%\hulst{}{\dhome/get1ff.m} \hulst{}{\dhome/get2ff.m}
\hulst{caption={{\small {\tt spcapinit, oosetfemops, sG, getmeancurv} and {\tt getV} from {\tt spcap1}. See sources for, e.g., 
{\tt get1ff}, {\tt get2ff}, and the Element2RefineSelector function {\tt e2rs} used for mesh adaptation. 
%This implementation is 
}},
label=spl1}{\dhome/getV.m}

Despite the caveats in Remark \ref{hrem1}, we can now set up a 
simple script (Listing \ref{spl2}) and produce a continuation 
diagram and sample plots fully analogous to \S\ref{scsec}. 
For mesh--refinement, to select 
triangles to refine we use the areas on $X$, see {\tt e2rs}, 
but otherwise the adaptive mesh refinement works as usual 
in the legacy setting of \pdep\ via {\tt oomeshada}. 

%\hulst{linerange=3-16}{\dhome/cmds1.m}
%\hulst{caption={{\small Short script {\tt cmds1.m} from {\tt geomtut/spcap}, and Element2RefineSelector function {\tt e2rs} used for mesh adaptation.}},label=spl2}{\dhome/e2rs.m}
\hulst{linerange=3-13,caption={{\small Short script {\tt cmds1.m} from {\tt geomtut/spcap}.}},label=spl2}{\dhome/cmds1.m}

As already said, the main advantage of this setup is simplicity in the 
sense that 
the function {\tt getmeancurv} (based on {\tt get1ff} and {\tt get2ff}) 
is a direct translation of the differential geometric definition. 
Moreover,  we can work with fixed preassembled differentiation 
matrices (independent of $X$) as long as the mesh (in $\Om$) is fixed. 
The main disadvantage is that this 
implementation of \reff{h1} is not a weak form but mixes 
FEM and FD differentiation matrices. 
%, and it does not allow for simple analytical expressions for Jacobians. 

%\input{spcapleg}

%\input{appHU}

\renewcommand{\refname}{References}
\taskip%\renewcommand{\arraystretch}{1.05}\renewcommand{\baselinestretch}{1}
\small
\bibliographystyle{alpha}
%\bibliographystyle{plain}
%\bibliography{/hh/hubib}
\newcommand{\etalchar}[1]{$^{#1}$}

\end{document}